\pdfoutput=1

\newif\ifjournal
\newif\iffull
\newif\ifau
\newif\ifshort
\ifdefined\journalbuild \journaltrue \else \fulltrue \fi
\ifdefined\aubuild \autrue \fi
\ifdefined\shortbuild \shorttrue \fi

\documentclass[11pt]{amsart}

\usepackage[T1]{fontenc}
\usepackage[utf8]{inputenc}
\usepackage{lmodern}
\usepackage{microtype}
\usepackage{amsmath,amssymb,amsthm}
\usepackage{mathtools}
\usepackage{graphicx}
\usepackage{booktabs}
\usepackage{longtable}
\usepackage{array}
\usepackage{enumitem}
\usepackage{xcolor}
\usepackage{newunicodechar}
\usepackage{tikz}
\usepackage[colorlinks=true,linkcolor=blue!50!black,citecolor=blue!50!black,urlcolor=blue!50!black]{hyperref}
\hypersetup{
  pdftitle={No free lattice contains an infinite ascending chain of copies of FL(3)},
  pdfauthor={William DeMeo and Ralph Freese},
  pdfkeywords={free lattice, sublattice, ascending chain, canonical form, self-embedding, Whitman's solution, Agda, machine-checked proof},
  pdfsubject={Mathematics Subject Classification (2020): Primary 06B25; Secondary 03B35},
}

\makeatletter
\newenvironment{expertread}
  {\expert@open\endgroup}
  {\begingroup\def\@currenvir{expertread}\expert@close}
\newenvironment{cutcandidate}
  {\cut@open\endgroup}
  {\begingroup\def\@currenvir{cutcandidate}\cut@close}
\newenvironment{fullonly}
  {\full@open\endgroup}
  {\begingroup\def\@currenvir{fullonly}\full@close}
\newenvironment{furthercut}
  {\further@open\endgroup}
  {\begingroup\def\@currenvir{furthercut}\further@close}
\newif\ifexpert@pending
\newif\ifcut@pending
\newif\iffull@pending
\newif\iffurther@pending
\newskip\expert@skip
\newcount\cut@count
\newcount\full@count
\newcount\further@count
\def\marks@point#1#2#3#4#5{%
  \expandafter\xdef\csname #1@#2@#3@page\endcsname{#4}%
  \expandafter\xdef\csname #1@#2@#3@y\endcsname{#5}}
\def\cut@point{\marks@point{cut}}
\let\cut@open\par
\let\cut@close\par
\let\full@open\par
\let\full@close\par
\let\further@open\par
\let\further@close\par
  \let\expert@open\par
  \let\expert@close\par
\makeatother

\setlist{nosep}
\theoremstyle{plain}
\newtheorem{theorem}{Theorem}[section]
\newtheorem{lemma}[theorem]{Lemma}
\newtheorem{corollary}[theorem]{Corollary}
\newtheorem{proposition}[theorem]{Proposition}
\theoremstyle{definition}
\newtheorem{definition}[theorem]{Definition}
\newtheorem{question}[theorem]{Question}

\theoremstyle{remark} %
\newtheorem{remark}[theorem]{Remark}
\newtheorem{computation}[theorem]{Computation}

\theoremstyle{plain}
\newtheorem{letteredthm}{Theorem}
\newtheorem{letteredcor}{Corollary}
\newtheorem{letteredlem}{Lemma}

\makeatletter
\def\lettered@setup{%
  \let\@upn\@firstofone
  \let\lettered@label\label
  \def\label##1{\ifvmode\AddToHookNext{para/begin}{\lettered@label{##1}}\else\lettered@label{##1}\fi}}
\newenvironment{thmA}{\lettered@setup\letteredthm}{\endletteredthm}
\newenvironment{thmB}{\lettered@setup\letteredthm}{\endletteredthm}
\newenvironment{thmC}{\lettered@setup\letteredthm}{\endletteredthm}
\newenvironment{thmCfour}{\lettered@setup\letteredthm}{\endletteredthm}
\newenvironment{thmCthree}{\lettered@setup\letteredthm}{\endletteredthm}
\newenvironment{thmCcoh}{\lettered@setup\letteredthm}{\endletteredthm}
\newenvironment{thmD}{\lettered@setup\letteredthm}{\endletteredthm}
\newenvironment{thmE}{\lettered@setup\letteredthm}{\endletteredthm}
\newenvironment{thmF}{\lettered@setup\letteredthm}{\endletteredthm}
\newenvironment{thmG}{\lettered@setup\letteredthm}{\endletteredthm}
\newenvironment{thmH}{\lettered@setup\letteredthm}{\endletteredthm}
\newenvironment{corAone}{\lettered@setup\letteredcor}{\endletteredcor}
\newenvironment{corAtwo}{\lettered@setup\letteredcor}{\endletteredcor}
\newenvironment{corAthree}{\lettered@setup\letteredcor}{\endletteredcor}
\newenvironment{corAfour}{\lettered@setup\letteredcor}{\endletteredcor}
\newenvironment{corEone}{\lettered@setup\letteredcor}{\endletteredcor}
\newenvironment{lemone}{\lettered@setup\letteredlem}{\endletteredlem}
\newenvironment{lemtwo}{\lettered@setup\letteredlem}{\endletteredlem}
\newenvironment{lemthree}{\lettered@setup\letteredlem}{\endletteredlem}
\newenvironment{lemfour}{\lettered@setup\letteredlem}{\endletteredlem}
\newenvironment{lemfive}{\lettered@setup\letteredlem}{\endletteredlem}
\newcommand{\formofC}[1]{\@ifundefined{r@thm:C-#1}{C}{\hyperref[thm:C-#1]{C}}}
\makeatother

\newcommand{\FL}{\mathbf{FL}}
\newcommand{\gen}[1]{\langle #1 \rangle}
\newcommand{\Fix}{\operatorname{Fix}}
\newcommand{\Per}{\operatorname{Per}}
\newcommand{\cx}{\kappa}
\newcommand{\agda}[1]{\texttt{#1}}
\newcommand{\pr}{^{\prime}}
\newcommand{\W}{\mathrm{(W)}}
\newcommand{\ba}[1]{\bar{#1}}
\newcommand{\agdadot}{\llap{\raisebox{1.15ex}[0pt][0pt]{\scalebox{0.55}{\ensuremath{\bullet}}}\kern0.28em}}

\usetikzlibrary{arrows.meta,calc,fit,backgrounds}
\tikzset{
  >={Stealth[length=4pt,width=3pt]},
  every picture/.append style={line width=0.4pt, font=\small},
  every label/.append style={inner sep=1.5pt, font=\small},
  elt/.style={circle, draw, fill=black, inner sep=0pt, minimum size=3.4pt},
  fixed/.style={elt},
  moved/.style={elt, fill=white},
  triple/.style={double, double distance=0.9pt, minimum size=4.2pt},
  upper/.style={},
  lower/.style={},
  term/.style={inner sep=2pt},
  prime/.style={term, draw, rounded corners=2pt, dashed, fill=white},
  leaf/.style={prime},
  le/.style={},
  maps/.style={->, densely dashed, shorten <=2pt, shorten >=2pt},
  uses/.style={->, shorten <=1pt, shorten >=1pt},
  region/.style={fill=black!7, draw=none},
  link/.style={draw, rounded corners=3pt},
  result/.style={draw, rounded corners=2pt, inner sep=2.5pt, align=center},
  narrow/.style={result, densely dashed},
  module/.style={draw, rounded corners=1.5pt, inner sep=2pt, font=\scriptsize, fill=white},
  ours/.style={module, fill=black!12},
  example/.style={ours, densely dashed},
  imports/.style={->, >={Stealth[length=3pt,width=2.2pt]}, shorten <=0.5pt, shorten >=0.5pt},
}

\newunicodechar{·}{\textperiodcentered}
\newunicodechar{ʷ}{\textsuperscript{w}}
\newunicodechar{ˢ}{\textsuperscript{s}}
\newunicodechar{̇}{\agdadot}
\newunicodechar{Θ}{\ensuremath{\Theta}}
\newunicodechar{σ}{\ensuremath{\sigma}}
\newunicodechar{τ}{\ensuremath{\tau}}
\newunicodechar{ω}{\ensuremath{\omega}}
\newunicodechar{ᵖ}{\textsuperscript{p}}
\newunicodechar{ᶜ}{\textsuperscript{c}}
\newunicodechar{₁}{\textsubscript{1}}
\newunicodechar{₃}{\textsubscript{3}}
\newunicodechar{₄}{\textsubscript{4}}
\newunicodechar{ℊ}{\textit{g}}
\newunicodechar{→}{\ensuremath{\to}}
\newunicodechar{⇔}{\ensuremath{\Leftrightarrow}}
\newunicodechar{∧}{\ensuremath{\wedge}}
\newunicodechar{∨}{\ensuremath{\vee}}
\newunicodechar{≈}{\ensuremath{\approx}}
\newunicodechar{≤}{\ensuremath{\leq}}
\newunicodechar{⊥}{\ensuremath{\bot}}
\newunicodechar{⨟}{\ensuremath{\mathbin{;}}}

\begin{document}

\title[No infinite ascending chain of copies of \texorpdfstring{$\FL(3)$}{FL(3)}]{No free lattice contains an infinite ascending chain of copies of \texorpdfstring{$\FL(3)$}{FL(3)}}

\author{William DeMeo}
\address{Input Output, Boulder, CO 80302, USA}
\email{williamdemeo@gmail.com}
\urladdr{https://orcid.org/0000-0003-1832-5690}

\author{Ralph Freese}
\address{Department of Mathematics, University of Hawaii, Honolulu, HI 96822, USA}
\email{ralph@math.hawaii.edu}
\urladdr{https://orcid.org/0000-0002-3441-0687}

\dedicatory{In memory of R.~A.~Dean (1924--2022)}

\begin{abstract}
No free lattice contains an infinite strictly ascending chain of sublattices isomorphic to $\FL(3)$, the free lattice on three generators; this answers a question that R.~A.~Dean asked and the second author's 1997 notes on the problem left open.  Every such chain would yield, inside $\FL(3)$, a coherent chain of proper self-embeddings $\Theta_0, \Theta_1, \Theta_2, \dots$ with $\Theta_k \circ \Theta_j = \Theta_j$ for $j < k$, so that each later embedding fixes the image of each earlier one pointwise.  Every proper self-embedding of $\FL(3)$ whose fixed sublattice contains a free generating triple fixes one of the three generators, say $z$, because a power of it or of its dual lies above the identity, and a self-embedding that moves every generator and lies above the identity fixes at most eleven elements.  Up to duality, a proper self-embedding fixing the image of such a $\Theta_1$ pointwise fixes only elements above $z \wedge (x \vee y)$, where no free generating triple containing $z$ can lie.  We also verify the one-step example of the 1997 notes, prove that no step toward such a chain is self-similar, and show that the configuration the argument rules out at length three does occur at length two.

Without freeness the answer reverses: $\FL(3)$ contains an infinite strictly ascending chain of $3$-generated sublattices, so the lattice analogue of Takahasi's chain condition for free groups fails.  The theorem itself is checked by the Agda proof assistant under its safe flags, for free lattices over generators with decidable equality, and its formal proof uses no classical principle: finite searches and an induction on a finite set replace the infinite pigeonhole principle by which a chain yields a coherent chain.
\end{abstract}

\maketitle
\tableofcontents

\section{The problem, its history, and the theorem}
\label{sec:problem}

\begin{expertread}
Let $\FL(X)$ be the free lattice on a set $X$ of generators, $\FL(n)$ the free lattice on $n$ generators, often labeled $x_1, \dots, x_n$, and $\FL(3) = \FL(\{x, y, z\})$.  A \emph{copy of $\FL(3)$} in a lattice $L$ is a sublattice of $L$ isomorphic to $\FL(3)$.  R.~A.~Dean asked the following question in an email shortly after the publication, in 1995, of the monograph \emph{Free Lattices} \cite{FJN} by Freese, Je\v{z}ek and Nation:

\begin{quote}
\textbf{Dean's problem}.  Is there a free lattice $\FL(X)$ with sublattices $L_0 \subsetneq L_1 \subsetneq L_2 \subsetneq \cdots$, every $L_k$ isomorphic to $\FL(3)$?
\end{quote}

Dean had retired from Caltech in 1987, and the second author, at work with Nation on free lattices, had tried to reach him and found no address.  The book reached him, and his email said that of the outstanding problems on free lattices he knew of, only this one was still open.

The second author's survey of free lattices, Appendix~G of Gr\"atzer's \emph{General Lattice Theory} \cite[p.~578]{FreeseAppendixG}, records the problem as follows: ``In a recent communication, R.~A.~Dean informed me of an old problem in free lattice theory: can a free lattice have an ascending chain of sublattices all isomorphic to $\FL(3)$?  This problem did not appear in print and I was unaware of it.''\footnote{The printed text writes $F(3)$, the book's notation, for $\FL(3)$; the revised version of the appendix \cite{FreeseAppendixGRevised} writes $\FL(3)$.}

The second author also wrote ``Notes on Dean's problem'' \cite{Freese1997}; these include an example of one step of such a chain, show that the obvious attack on the problem fails, and provide an ``almost proof'' that an ascending chain cannot exist.  As far as we know, Appendix~G of Gr\"atzer's book is the only place the problem appears in print.  This paper answers the question.\footnote{In the time between the first and second versions of this manuscript, the authors learned of a preprint by M.~Lupini \cite{Lupini2026}, dated 19 October 2026 and posted on 24 October 2026.  Its Theorem~1.1 states that, for every finite $r$, the free sublattices of rank at most $r$ of a free lattice satisfy the ascending chain condition; for $r = 3$ this is Theorem~\ref{thm:H}.  It also contains an ascending chain of $3$-generated sublattices of $\FL(3)$ built by the recursion of Theorem~\ref{thm:ladder} (its Theorem~1.2), a form of Theorem~\ref{thm:A} and Corollary~\ref{cor:A4} for the endomorphisms of a free lattice of finite rank (its Theorem~1.3), and the self-embedding $\theta_1$ of Theorem~\ref{thm:D}, up to the names of the generators, with $\Fix(\theta_1) = \gen{p, q, z}$ (its Theorem~6.1).  The present paper was written with no knowledge of Lupini's work.}
\end{expertread}

\begin{thmH}[Dean's problem]
No free lattice contains an infinite strictly ascending chain of sublattices isomorphic to $\FL(3)$.
\end{thmH}

\subsection{The shape of the question}
\label{sec:shape}

The sublattices of free lattices have been studied since Whitman's two papers \cite{Whitman1941, Whitman1942}, which solved the word problem for $\FL(X)$, found the canonical form of its elements, and showed that $\FL(3)$ contains a copy of $\FL(\omega)$, the free lattice on countably many generators; the monograph \cite{FJN}, which we cite as FJN and by theorem number throughout, tells the history in its introduction.

Dean, J\'onsson, Galvin and Kiefer took the subject up in the 1950s and 1960s, with the finite sublattices as a main interest.  J\'onsson and Kiefer proved that free lattices, hence all their sublattices, are semidistributive, and that a finite lattice satisfying Whitman's condition $\W$ (Section~\ref{sec:whitman}) has breadth at most four, so that the finite sublattices of a free lattice satisfy a nontrivial lattice equation \cite{JonssonKiefer1962} (FJN, Theorems~1.21 and~1.30 and Corollary~1.31).  Dean's paper on sublattices of free lattices \cite{Dean1961} shows, among other things, that the lattice freely generated by any countable ordered set embeds in $\FL(3)$ (FJN, Chapter~V, Section~2).  J\'onsson conjectured that a finite lattice embeds in a free lattice if and only if it is semidistributive and satisfies $\W$, and Nation proved it \cite{Nation1982} (FJN, Theorem~5.55).

So the finite sublattices of a free lattice are characterized.  Dean's question concerns the infinite ones, specifically $\FL(3)$, which contains a copy of every countable one; it asks how the copies of $\FL(3)$ inside a free lattice can sit inside one another.  Read inside $\FL(3)$ (Theorem~\ref{thm:C}), the answer is a statement about its self-embeddings: no infinite sequence of proper self-embeddings of $\FL(3)$ has each later one fixing the images of all earlier ones pointwise (Corollary~\ref{cor:dean-iff-coherent}).

Three classical facts, all from FJN, fix the question's shape.

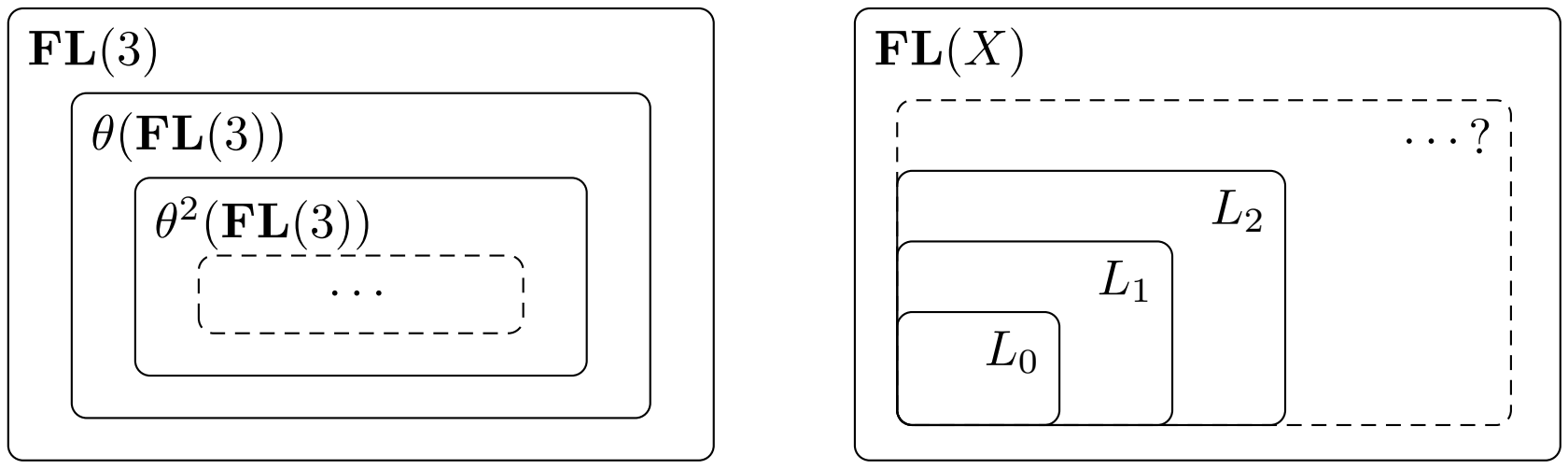
\begin{figure}[ht]
\centering
\begin{tikzpicture}[x=10mm, y=10mm]
  \foreach \k/\name in {0/{$\FL(3)$}, 1/{$\theta(\FL(3))$}, 2/{$\theta^2(\FL(3))$}} {
    \draw[link] (0.45*\k, 0.3*\k) rectangle (5-0.45*\k, 3.2-0.6*\k);
    \node[anchor=north west] at (0.45*\k, 3.2-0.6*\k) {\name};
  }
  \draw[link, densely dashed] (1.35,0.9) rectangle (3.65,1.45);
  \node at (2.5,1.175) {$\cdots$};
  \begin{scope}[xshift=60mm]
    \draw[link] (0,0) rectangle (5,3.2);
    \node[anchor=north west] at (0,3.2) {$\FL(X)$};
    \foreach \k in {0,1,2} {
      \draw[link] (0.3,0.25) rectangle (1.45+0.8*\k, 1.05+0.5*\k);
      \node[anchor=north east] at (1.45+0.8*\k, 1.05+0.5*\k) {$L_{\k}$};
    }
    \draw[link, densely dashed] (0.3,0.25) rectangle (4.65,2.55);
    \node[anchor=north east] at (4.65,2.55) {$\cdots\,?$};
  \end{scope}
\end{tikzpicture}
\caption{Left, the descending chain $\FL(3) \supsetneq \theta(\FL(3)) \supsetneq \theta^2(\FL(3)) \supsetneq \cdots$ that a proper self-embedding $\theta$ of $\FL(3)$ gives; right, the ascending chain $L_0 \subsetneq L_1 \subsetneq L_2 \subsetneq \cdots$ of copies of $\FL(3)$ in a free lattice that Dean asks for.}
\label{fig:chains}
\end{figure}

\begin{enumerate}[label=(\roman*)]
\item Whitman proved that $\FL(\omega)$, the free lattice on countably many generators, embeds in $\FL(3)$ (FJN, Theorem~1.28).  Hence if one free lattice has such a chain then every $\FL(X)$ with $|X| \geq 3$ does, and $\FL(3)$ is the only case.  (For $|X| \leq 2$ the free lattice is finite.)
\item A \emph{descending} chain is easy.  $\FL(3)$ has proper sublattices isomorphic to itself (the elements $x \vee (y \wedge z)$, $y \vee (z \wedge x)$, $z \vee (x \wedge y)$ generate one, Section~\ref{sec:inside}), and iterating a self-embedding $\theta$ with such an image gives $\FL(3) \supsetneq \theta(\FL(3)) \supsetneq \theta^2(\FL(3)) \supsetneq \cdots$.  Figure~\ref{fig:chains} shows the two directions.
\item A copy of $\FL(3)$ inside $\FL(X)$ determines its triple of free generators up to order, so a copy and its free triple are the same datum: the free generators of a free lattice are its elements that are both join prime and meet prime, and the automorphisms of $\FL(X)$ are the permutations of $X$ (FJN, Lemma~1.4 and Corollaries~1.5 and~1.6).
\end{enumerate}

The obvious attack fails.  Suppose $\gen{a, b, c} \subseteq \gen{a\pr, b\pr, c\pr}$, where $\gen{S}$ denotes the sublattice generated by $S$ and both triples freely generate copies of $\FL(3)$.  Then $a$, $b$, $c$ are lattice terms in $a\pr$, $b\pr$, $c\pr$, and one might hope that $a\pr$, $b\pr$, $c\pr$ must occur as subterms of $a$, $b$, $c$, which would force a chain to stabilize, since a term has finitely many subterms.  The one-step example of the 1997 notes (verified here in Section~\ref{sec:freese-example}) shows that hope is false: each new generator in the step has complexity~$4$ and each old one complexity~$3$, and the notes add that the new generators may be chosen as complex as one likes.

\emph{Finite} ascending chains of every length exist trivially.  (Read the chain of~(ii) upward, starting from $\theta^n(\FL(3))$.)  The content of the question is the infinite chain.

\subsection{The shape of the proof}
\label{sec:proof-shape}

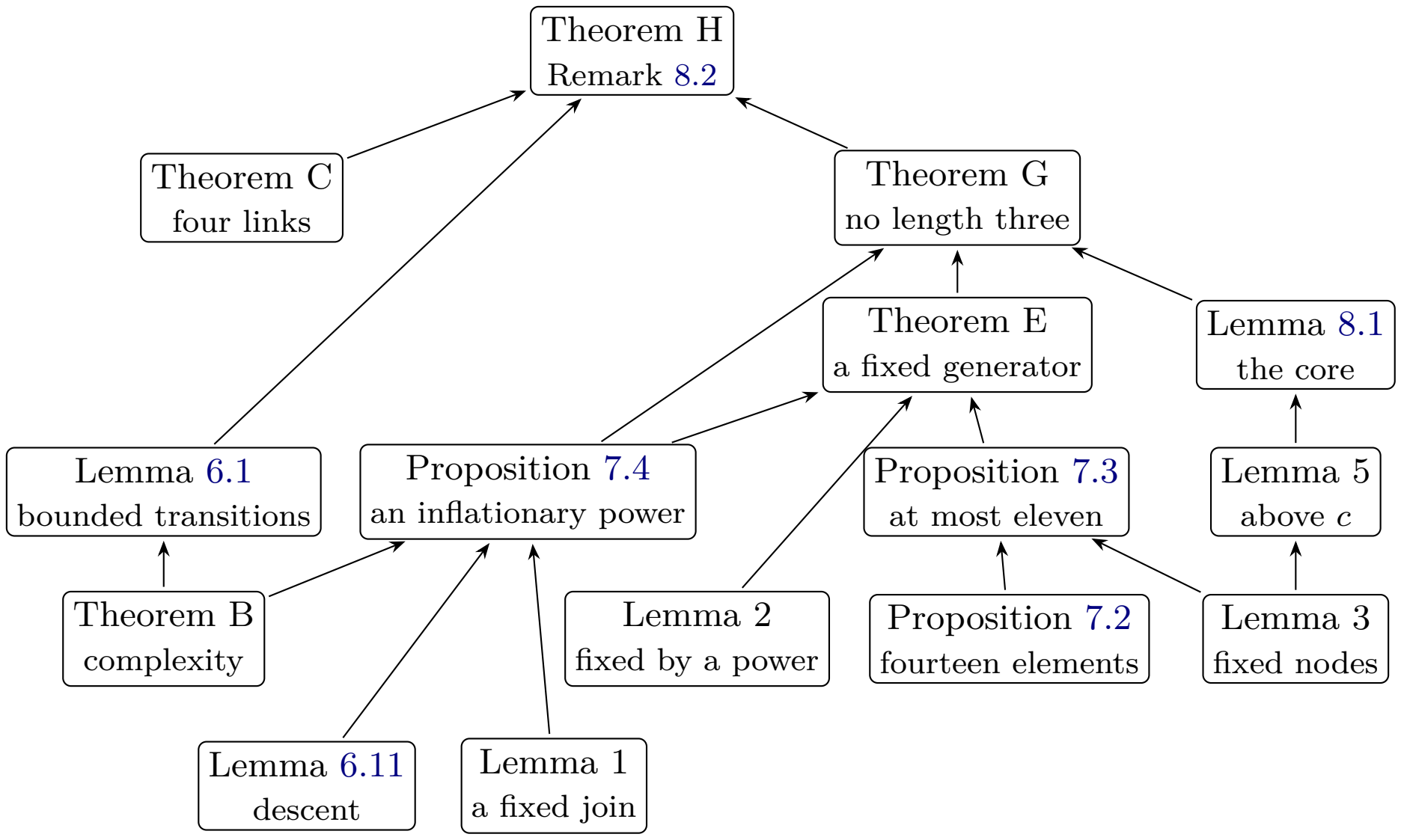
\begin{figure}[ht]
\centering
\begin{tikzpicture}[x=11.5mm, y=13mm, every node/.append style={font=\footnotesize}]
  \node[result] (descent) at (-3.3,0) {Lemma~\ref{lem:descent}\\{\scriptsize descent}};
  \node[result] (L1) at (-1.4,0) {Lemma 1\\{\scriptsize a fixed join}};
  \node[result] (B) at (-4.4,1) {Theorem B\\{\scriptsize complexity}};
  \node[result] (L2) at (-0.3,1) {Lemma 2\\{\scriptsize fixed by a power}};
  \node[result] (fourteen) at (2.1,1) {Proposition~\ref{prop:fourteen}\\{\scriptsize fourteen elements}};
  \node[result] (L3) at (4.3,1) {Lemma 3\\{\scriptsize fixed nodes}};
  \node[result] (bound) at (-4.4,2) {Lemma~\ref{lem:bounded-transitions}\\{\scriptsize bounded transitions}};
  \node[result] (power) at (-1.6,2) {Proposition~\ref{prop:inflationary-power}\\{\scriptsize an inflationary power}};
  \node[result] (eleven) at (2,2) {Proposition~\ref{prop:eleven}\\{\scriptsize at most eleven}};
  \node[result] (L5) at (4.3,2) {Lemma 5\\{\scriptsize above $c$}};
  \node[result] (E) at (1.7,3) {Theorem E\\{\scriptsize a fixed generator}};
  \node[result] (core) at (4.3,3) {Lemma~\ref{lem:core}\\{\scriptsize the core}};
  \node[result] (C4) at (-3.8,4) {Theorem C\\{\scriptsize four links}};
  \node[result] (G) at (1.7,4) {Theorem G\\{\scriptsize no length three}};
  \node[result] (H) at (-0.8,5) {Theorem H\\{\scriptsize Remark~\ref{rem:constructive}}};
  \draw[uses] (B) -- (bound);
  \draw[uses] (B) -- (power);
  \draw[uses] (descent) -- (power);
  \draw[uses] (L1) -- (power);
  \draw[uses] (fourteen) -- (eleven);
  \draw[uses] (L3) -- (eleven);
  \draw[uses] (L3) -- (L5);
  \draw[uses] (power) -- (E);
  \draw[uses] (L2) -- (E);
  \draw[uses] (eleven) -- (E);
  \draw[uses] (L5) -- (core);
  \draw[uses] (power) -- (G);
  \draw[uses] (E) -- (G);
  \draw[uses] (core) -- (G);
  \draw[uses] (G) -- (H);
  \draw[uses] (C4) -- (H);
  \draw[uses] (bound) -- (H);
\end{tikzpicture}
\caption{The results on the way to Theorem~\ref{thm:H}: an arrow runs from a result to each proof that cites it, and the proof of Theorem~\ref{thm:H} is the constructive one of Remark~\ref{rem:constructive}.  The results of Section~\ref{sec:preliminaries}, which every part uses, are left out, and so is Lemma~\ref{lem:two-generators}, that a self-embedding of $\FL(3)$ fixing two generators is the identity, which the proofs of Lemma~\ref{lem:2} and Lemma~\ref{lem:core} cite.}
\label{fig:map}
\end{figure}

Fix a chain $L_0 \subsetneq L_1 \subsetneq \cdots$ in $\FL(X)$ with isomorphisms $\varphi_k \colon \FL(3) \to L_k$.  Reading the inclusion $L_j \subseteq L_k$ through the isomorphisms gives, for $j \leq k$, an injective endomorphism $\theta_{j,k} = \varphi_k^{-1} \circ \varphi_j$ of $\FL(3)$, which is onto if and only if $j = k$.  The argument has four parts, which we now sketch.  Figure~\ref{fig:map} shows which of these results each proof uses.  Nonstandard terms used here are defined in the sequel.

\begin{enumerate}[label=(\arabic*)]
\item \textit{Theorem~\ref{thm:B}.  A lattice embedding never lowers the complexity of an element, the depth of its canonical form.}

Since $\varphi_0 = \varphi_k \circ \theta_{0,k}$, the complexity of $\theta_{0,k}(g)$ is at most that of $\varphi_0(g)$ for each generator $g$, so the maps $\theta_{0,k}$ take finitely many values, and two of them coincide, $\theta_{0,j} = \theta_{0,k}$ with $j < k$; then $\theta_{j,k}$ is a proper self-embedding of $\FL(3)$ that fixes the copy $\theta_{0,j}(\FL(3))$ of $\FL(3)$ pointwise.  By a diagonal argument the whole chain is, on a subsequence, a coherent chain (Definition~\ref{def:coherent}) of proper self-embeddings $\Theta_0, \Theta_1, \Theta_2, \dots$ of $\FL(3)$ with $\Theta_k \circ \Theta_j = \Theta_j$ for $j < k$ (the coherent form of Theorem~\formofC{coherent}); we call the $L_k$ the links of the chain and the $\Theta_k$ those of the coherent chain.  The diagonal argument uses the infinite pigeonhole principle, but the proof needs only a coherent chain of length three, and four links $L_i \subsetneq L_j \subsetneq L_k \subsetneq L_l$ with $\theta_{i,j} = \theta_{i,l}$ and $\theta_{j,k} = \theta_{j,l}$ give one (the four-link form of Theorem~\formofC{four}).

\item \textit{Theorem~\ref{thm:E}.  Every proper self-embedding of $\FL(3)$ whose fixed sublattice contains a free triple fixes one of the generators.}

A positive power of such a self-embedding, or of its dual, is inflationary, lying above the identity (Section~\ref{sec:inflationary}, Proposition~\ref{prop:inflationary-power}), and the engine of the proof is a count: the elements of $\FL(3)$ whose canonical trees have no generator among the meetands of any meet are fourteen (Proposition~\ref{prop:fourteen}, a proof, not a computation), every fixed point of an inflationary self-embedding that moves every generator is one of the eleven that are not generators (Proposition~\ref{prop:eleven}), and a fixed free triple supplies fixed points of arbitrarily large complexity.

\item \textit{Lemma~\ref{lem:5}.  Every fixed point of an inflationary self-embedding that fixes $z$ and moves $x$ and $y$ lies above $c = z \wedge (x \vee y)$.}

\item \textit{Theorem~\ref{thm:G}.  There is no coherent chain of length three.}

In a coherent chain of length three, $\Theta_1$ fixes a generator, say $z$, and a power of $\Theta_2$ or of its dual is such a self-embedding fixing $\Theta_1(\FL(3))$.  In the first case $\Theta_1(y) \geq c$ by (3), while $\Theta_1(x) \leq x \vee y$ by the partition lemma at $z$ (Section~\ref{sec:whitman}), so $\Theta_1(x) \wedge \Theta_1(z) \leq (x \vee y) \wedge z = c \leq \Theta_1(y)$, which the free triple $\Theta_1(x), \Theta_1(y), \Theta_1(z)$ forbids; the dual case is the same argument for the dual of $\Theta_1$.
\end{enumerate}

Parts (2) to (4) are finite.  Part (1) is not, as stated, but the proof needs only the four-link form of Theorem~\formofC{four}, and finite searches with an induction on the number of values of $\theta_{0,\cdot}$ replace the infinite pigeonhole principle (Remark~\ref{rem:constructive}).  So Theorem~\ref{thm:H} itself has a complete constructive formalization in the Agda proof assistant, for free lattices over generators with decidable equality, which classically is every free lattice; Section~\ref{sec:formal} says exactly what is formally verified, and at which commit of the release of the Agda repository that accompanies the paper.

We note that the configuration part (4) rules out at length three does occur at length two: there is a proper self-embedding of $\FL(3)$ fixing a free triple pointwise (Theorem~\ref{thm:D}), so no shorter argument of this kind is available.

\subsection{Free groups}
\label{sec:free-groups}

Free lattices do not, in general, behave like free groups, but the analogue of Dean's question for free groups also has a negative answer.  Takahasi \cite{Takahasi1950} and, independently, Higman \cite{Higman1951} proved that an ascending chain of subgroups of bounded rank in a free group stabilizes, as Takahasi's later note recalls \cite{Takahasi1951}.  So the analogue, with copies of the free group of rank three in place of copies of $\FL(3)$, has a short answer, and that answer does not need the links to be free, since every subgroup of a free group is free.

The argument of Takahasi's note bounds the lengths of the free generators of a subgroup $U$ containing a finitely generated subgroup $H$ by the lengths of the generators of $H$, provided no proper \emph{free factor} of $U$ (a subgroup generated by part of a free basis of $U$) contains $H$ \cite[Theorem~2]{Takahasi1951}.  In a free lattice the analogous bound fails, although the corresponding proviso is automatic there: two of the free generators of a copy of $\FL(3)$ generate a finite sublattice, which contains no copy of $\FL(3)$.  In the one-step example of Section~\ref{sec:freese-example}, each generator of the larger copy of $\FL(3)$ has complexity~$4$ and each generator of the smaller one complexity~$3$.

Freeness of the links is essential.  Without it the answer reverses: $\FL(3)$ contains an infinite strictly ascending chain of $3$-generated sublattices, none of them isomorphic to $\FL(3)$ (Theorem~\ref{thm:ladder}), so the lattice analogue of Takahasi's chain condition fails.  Theorem~\ref{thm:H} is the statement that it holds \emph{when the links are free}, that is, copies of $\FL(3)$.  Its proof rests on the monotonicity that does hold, Theorem~\ref{thm:B}, which bounds not the generators of the links but the maps between them.

The analogy runs one step further.  The fixed subgroup of an injective endomorphism of a free group of rank $n$ has rank at most $n$ \cite{ImrichTurner1989, BestvinaHandel1992}.  Whether the fixed sublattice of a proper self-embedding of $\FL(3)$ is always finitely generated is open (Section~\ref{sec:questions}); by Theorem~\ref{thm:D} it can contain a copy of $\FL(3)$.

\subsection{What is new}
\label{sec:new}

The main contributions of the paper, beyond the answer to Dean's question, are the following:
\begin{itemize}
\item Theorem~\ref{thm:A}, that an injective endomorphism of a free lattice which maps a finitely generated sublattice onto a superset of it maps it onto itself, closes the self-similar case: no chain is generated by one injective endomorphism of the ambient lattice, nor by its powers (Corollaries~\ref{cor:A1} to~\ref{cor:A3}).
\item Theorem~\ref{thm:B}, that a lattice embedding never lowers complexity, is Lemma~9.13 of FJN read level by level.
\item Theorem~\ref{thm:C} sharpens the finiteness argument of the notes, in which only the pattern by which the generators of a later link cut the term trees of the bottom link's generators, one of finitely many, recurs along the chain, to the recurrence of the whole term that expresses the bottom link in a later one, and in its coherent form it makes Dean's problem equivalent to the existence of an infinite coherent chain (Corollary~\ref{cor:dean-iff-coherent}).
\item Theorem~\ref{thm:D} answers the fixed-sublattice question (Question~\ref{q:fixed}) positively: a proper self-embedding of $\FL(3)$ can fix a copy of $\FL(3)$ pointwise, so the route to a negative answer through that question is closed, and coherent chains of length two exist.
\item Theorem~\ref{thm:E} and Corollary~\ref{cor:E1} describe every proper self-embedding with a fixed free triple: it fixes exactly one generator, and, up to duality and a power, it moves the other two up inside the ideal below their join.
\item Lemma~\ref{lem:5} and Theorem~\ref{thm:G} prove that every such self-embedding is rigid, in that no proper self-embedding fixes its image pointwise, so that no coherent chain of length three exists.%
\item Theorem~\ref{thm:F} is a first form of that rigidity, for a class of self-embeddings.%
\item Theorem~\ref{thm:ladder} shows that the freeness of the links cannot be dropped.
\end{itemize}

The one-step example is verified (Section~\ref{sec:freese-example}), and the almost-proof of the notes is formalized up to its gap, which is stated as a decidable predicate (Section~\ref{sec:almost-proof}).

The formal development (Section~\ref{sec:formal}) adds, besides the check of Theorem~\ref{thm:H} and the other finite results above, in the forms Section~\ref{sec:formal} states, a formalization of Whitman's theory of free lattices, including his rules as an inductive relation on terms, proved sound and complete for the order of every lattice, so that the terms modulo the relation form $\FL(X)$, with a decidable order and the universal property; the canonical form, decided, and the refinement theorem (Theorems~\ref{thm:canonical} and~\ref{thm:refine}, in the forms Section~\ref{sec:formal} states); free generation, decided by Corollary~\ref{cor:1.13}; Lemma~\ref{lem:9.13}; the embedding of $\FL(\omega)$ in $\FL(3)$; and the fixed and periodic points of self-embeddings.

Whitman's recursion has been run inside a proof assistant before: James and Hinze implemented it in Coq as a reflective tactic, proved sound for every lattice and terminating~\cite{JamesHinze2011}, and Guilloud and Pit-Claudel proved a sequent calculus for orthologic sound and complete for ortholattices in Coq, with a proof search of quadratic cost~\cite{GuilloudPitClaudel2025}.  We know of no machine-checked proof that Whitman's rules are complete, and none of the canonical form, of the structure of the sublattices of a free lattice, or of the fixed points of its self-embeddings, in Lean's Mathlib, Isabelle's Archive of Formal Proofs, the Rocq (formerly Coq) and Agda libraries on GitHub, and the Mizar Mathematical Library, searched by name and by full text in October 2026.  As far as we know, our formal development is the first of each, and Theorem~\ref{thm:H} is the first significant theorem about free lattices beyond the word problem to be checked by machine.

\subsection{Contents}

Section~\ref{sec:preliminaries} collects what we use from \cite{FJN} and proves Theorem~\ref{thm:B}.  Section~\ref{sec:freese-example} verifies the one-step example and reads its step from inside $\FL(3)$.  Section~\ref{sec:chains} reformulates a chain as a sequence of self-embeddings and proves that the variant without freeness has a negative answer, and Section~\ref{sec:theorem-a} proves that no step toward a chain is self-similar (Theorem~\ref{thm:A}), which closes the approach a computer search would have followed.  Section~\ref{sec:coherent} proves Theorem~\ref{thm:C} in its finite and coherent forms, states the fixed-sublattice question, and answers it positively (Theorem~\ref{thm:D}).  Section~\ref{sec:structure} proves Theorem~\ref{thm:E} and Section~\ref{sec:rigidity} proves Lemma~\ref{lem:5}, Theorem~\ref{thm:G} and Theorem~\ref{thm:H}, the last both classically and constructively.  Section~\ref{sec:almost-proof} formalizes the almost-proof of the notes up to its gap and decides the gap on the one-step example, as the record of the path that was superseded.  Section~\ref{sec:formal} describes the formal development: what it checks, section by section; how it and the proofs were produced; the name of each statement in the Agda source; what the checking costs; and which computer searches were run and what they found.  Section~\ref{sec:questions} lists what remains.

\section{Preliminaries}
\label{sec:preliminaries}

Everything in this section is in Chapter~I of FJN \cite{FJN}, except Lemma~9.13 of its Chapter~IX, Theorem~\ref{thm:B}, which is ours, and Lemma~\ref{lem:finite}, an elementary count that FJN does not state.  We state each result we use with its number in the book.  \emph{Dually} refers to the statement obtained by exchanging $\vee$ and $\wedge$ and reversing the order; the duals of an element of $\FL(3)$ and of a self-embedding are defined in Section~\ref{sec:inflationary}.

\subsection{Terms, rank and complexity}
\label{sec:terms}

A \emph{lattice term} over $X$ is a generator $x \in X$ or a formal join $t_1 \vee \cdots \vee t_n$ or formal meet $t_1 \wedge \cdots \wedge t_n$ of terms, $n > 1$.  Its \emph{rank} (FJN, Chapter~I, Section~2, where it is also called the length) is $1$ for a generator and one more than the sum of the ranks of $t_1, \dots, t_n$ for a join or a meet; so $x \vee y \vee z$ has rank $4$ and $x \vee (y \vee z)$ rank $5$.  Its \emph{complexity} is $0$ for a generator and one more than the largest complexity among $t_1, \dots, t_n$ for a join or a meet: the depth of the term tree.  The subterms of a term are defined as usual.  A term $t$ over $X$ names an element $t^{L}$ of any lattice $L$ with a map $X \to L$; in $\FL(X)$ every element is named by a term, and we write $s \leq t$ for $s^{\FL(X)} \leq t^{\FL(X)}$.

\subsection{Whitman's solution to the word problem}
\label{sec:whitman}

\begin{theorem}[FJN, Theorem~1.8 and Corollary~1.9]
\label{thm:W}
The free lattice $\FL(X)$ satisfies Whitman's condition
\begin{quote}
$\W$\quad if $v = v_1 \wedge \cdots \wedge v_r \leq u_1 \vee \cdots \vee u_s = u$, then $v_i \leq u$ for some $i$ or $v \leq u_j$ for some $j$.
\end{quote}
The condition does not mention the generators, so every sublattice of a free lattice satisfies $\W$.
\end{theorem}

\begin{lemma}[FJN, Lemma~1.4 and Corollary~1.5]
\label{lem:prime}
Let $X$ generate the lattice $L$ and suppose $\bigwedge S \leq \bigvee T$ implies $S \cap T \neq \emptyset$ for all finite nonempty $S, T \subseteq X$ (which holds in $\FL(X)$, FJN, Lemma~1.2).  Then every $x \in X$ is join prime and meet prime, and if $X$ is the disjoint union of $Y$ and $Z$, then $L$ is the disjoint union of the filter generated by $Y$ and the ideal generated by $Z$.
\end{lemma}

The last clause is the \emph{partition lemma} in the form in which we use it: for a finite set $X$ of at least two generators and $z \in X$ (the partition lemma \emph{at} $z$), every element $w$ of $\FL(X)$ satisfies exactly one of $z \leq w$ and $w \leq \bigvee (X \setminus \{z\})$, and dually exactly one of $w \leq z$ and $w \geq \bigwedge (X \setminus \{z\})$.  In $\FL(3)$: every $w$ lies above $z$ or below $x \vee y$, and not both.  We write $[u, v]$ for the interval between elements $u \leq v$, and $(u]$ and $[u)$ for the principal ideal and the principal filter of $u$; so $\FL(3)$ is the disjoint union of $[z)$ and $(x \vee y]$.

\begin{theorem}[FJN, Theorem~1.11; Whitman]
\label{thm:whitman}
For terms $s$ and $t$ over $X$, whether $s \leq t$ holds in $\FL(X)$ is decided by the following rules, each of which answers the question or reduces it to questions about proper subterms.
\begin{enumerate}[label=(\arabic*)]
\item If $s = x_i$ and $t = x_j$ are generators, then $s \leq t$ if and only if $x_i = x_j$.
\item If $s = s_1 \vee \cdots \vee s_k$, then $s \leq t$ if and only if $s_i \leq t$ for all $i$.
\item If $t = t_1 \wedge \cdots \wedge t_k$, then $s \leq t$ if and only if $s \leq t_i$ for all $i$.
\item If $s = x_i$ and $t = t_1 \vee \cdots \vee t_k$, then $s \leq t$ if and only if $x_i \leq t_j$ for some $j$.
\item If $s = s_1 \wedge \cdots \wedge s_k$ and $t = x_i$, then $s \leq t$ if and only if $s_j \leq x_i$ for some $j$.
\item If $s = s_1 \wedge \cdots \wedge s_k$ and $t = t_1 \vee \cdots \vee t_m$, then $s \leq t$ if and only if $s_i \leq t$ for some $i$ or $s \leq t_j$ for some $j$.
\end{enumerate}
\end{theorem}

Rules (4) and (5) are the primality of the generators, and rule (6) is $\W$.  We call one run of this procedure a \emph{Whitman decision}; every computation in this paper that is reported as checked by the proof assistant is a finite number of Whitman decisions on explicit terms.

\subsection{Free generating sets}
\label{sec:free-generation}

\begin{corollary}[FJN, Corollary~1.6]
\label{cor:1.6}
Let $X$ generate $L$ as in Lemma~\ref{lem:prime}.  If $Y$ generates $L$ then $X \subseteq Y$, and every automorphism of $L$ is induced by a permutation of $X$.
\end{corollary}

The proof of the first statement is the form in which we use it: a generator is join and meet irreducible, so if $x \in X$ lies in the sublattice generated by $Y$ then $x \in Y$.  In particular, for a homomorphism $\sigma \colon \FL(I) \to \FL(X)$ determined by the elements $s_i = \sigma(i)$, a generator $x$ of $\FL(X)$ lies in the image of $\sigma$ only if $x = s_i$ for some $i$.

\begin{corollary}[FJN, Corollaries~1.12 and~1.13]
\label{cor:1.13}
A subset $S$ of a free lattice $\FL(X)$ generates a sublattice isomorphic to a free lattice, freely generated by $S$, if and only if for all $s \in S$ and all finite nonempty $Y \subseteq S$,
\[
s \leq \textstyle\bigvee Y \text{ implies } s \in Y \tag{J}\label{eq:join-condition}
\]
and
\[
\textstyle\bigwedge Y \leq s \text{ implies } s \in Y. \tag{M}\label{eq:meet-condition}
\]
\end{corollary}
We call \eqref{eq:join-condition} the \emph{join condition} and its dual \eqref{eq:meet-condition} the \emph{meet condition}.

We say a family $(s_i)_{i \in I}$ of elements of $\FL(X)$ \emph{freely generates}, or is \emph{free}, when the homomorphism $\FL(I) \to \FL(X)$ sending $i$ to $s_i$ is injective; by Corollary~\ref{cor:1.13} this holds exactly when the $s_i$ are pairwise distinct and the set $\{s_i\}$ satisfies the two conditions.  A \emph{free triple} is a triple $(a, b, c)$ that freely generates, and $\gen{a, b, c}$ is then a copy of $\FL(3)$.  For a triple the two conditions say: $a \not\leq b \vee c$, $b \not\leq a \vee c$, $c \not\leq a \vee b$, $b \wedge c \not\leq a$, $a \wedge c \not\leq b$, $a \wedge b \not\leq c$; the one-element instances ($a \not\leq b$ and so on) follow, and there is no empty instance, since $Y$ is nonempty: lattice terms have no empty joins or meets.  Each of the six is one Whitman decision.  We identify a triple $\sigma = (a, b, c)$ with the homomorphism $\FL(3) \to \FL(X)$ it induces, $x \mapsto a$, $y \mapsto b$, $z \mapsto c$, so that $\sigma(x) = a$, and write $\gen{\sigma}$ for $\gen{a, b, c}$.

An \emph{embedding} is an injective homomorphism of lattices; a \emph{self-embedding} of a lattice is an injective endomorphism, and it is \emph{proper} when it is not onto.  A self-embedding of $\FL(3)$ is determined by the images of $x$, $y$, $z$, which form a free triple, and it is proper exactly when some generator is none of the three images, by Corollary~\ref{cor:1.6}.  We write a self-embedding as a substitution, $\theta \colon x \mapsto a, y \mapsto b, z \mapsto c$, and $\theta(w)$ for the element it assigns to $w$.  An embedding reflects the order: $\theta(s) \leq \theta(t)$ implies $s \leq t$, because $\theta(s \wedge t) = \theta(s)$ forces $s \wedge t = s$.

\subsection{Canonical form}
\label{sec:canonical}

\begin{theorem}[FJN, Theorems~1.17 and~1.18]
\label{thm:canonical}
Each element $w$ of $\FL(X)$ is named by a term of least rank, unique up to commutativity, its \emph{canonical form}.  A term $t = t_1 \vee \cdots \vee t_n$ with $n > 1$ is in canonical form if and only if
\begin{enumerate}[label=(\arabic*)]
\item each $t_i$ is a generator or formally a meet,
\item each $t_i$ is in canonical form,
\item $t_i \not\leq t_j$ for all $i \neq j$ (the $t_i$ form an antichain),
\item if $t_i = \bigwedge_j t_{ij}$ then $t_{ij} \not\leq t$ for all $j$;
\end{enumerate}
dually for a meet; a generator is always in canonical form.
\end{theorem}

\begin{figure}[ht]
\centering
\begin{tikzpicture}[x=8mm, y=9mm, every node/.append style={inner sep=1.5pt}]
  \node (r) at (0,2.4) {$\vee$};
  \node (rx) at (-1.6,1.2) {$x$};
  \node[text=black!50] (rm) at (0,1.2) {$\wedge$};
  \node[text=black!50] (rmx) at (-0.5,0) {$x$};
  \node[text=black!50] (rmy) at (0.5,0) {$y$};
  \node (rn) at (1.6,1.2) {$\wedge$};
  \node (rny) at (1.1,0) {$y$};
  \node (rnj) at (2.1,0) {$\vee$};
  \node (rnjx) at (1.6,-1.2) {$x$};
  \node (rnjz) at (2.6,-1.2) {$z$};
  \draw (rx) -- (r) -- (rn) -- (rny);
  \draw (rn) -- (rnj) -- (rnjx);
  \draw (rnj) -- (rnjz);
  \draw[black!50, densely dashed] (r) -- (rm) -- (rmx);
  \draw[black!50, densely dashed] (rm) -- (rmy);
  \node[text=black!50, font=\footnotesize] at (0,-0.75) {$x \wedge y \leq x$};
  \draw[->] (3.3,1.2) -- node[above=2pt, font=\footnotesize] {canonical form} (5.6,1.2);
  \begin{scope}[xshift=64mm, every label/.append style={font=\footnotesize, text=black!60}]
    \node[label=right:{$3$}] (c) at (0,2.4) {$\vee$};
    \node (cx) at (-1,1.2) {$x$};
    \node[label=right:{$2$}] (cn) at (1,1.2) {$\wedge$};
    \node (cny) at (0.5,0) {$y$};
    \node[label=right:{$1$}] (cnj) at (1.5,0) {$\vee$};
    \node (cnjx) at (1,-1.2) {$x$};
    \node (cnjz) at (2,-1.2) {$z$};
    \draw (cx) -- (c) -- (cn) -- (cny);
    \draw (cn) -- (cnj) -- (cnjx);
    \draw (cnj) -- (cnjz);
  \end{scope}
\end{tikzpicture}
\caption{Left, the term $x \vee (x \wedge y) \vee (y \wedge (x \vee z))$, whose joinand $x \wedge y$ lies below its joinand $x$ and is deleted; right, the canonical tree of the same element, $x \vee (y \wedge (x \vee z))$, with the complexity of each node that is not a generator.  Its canonical joinands are $x$ and $y \wedge (x \vee z)$.}
\label{fig:canonical}
\end{figure}
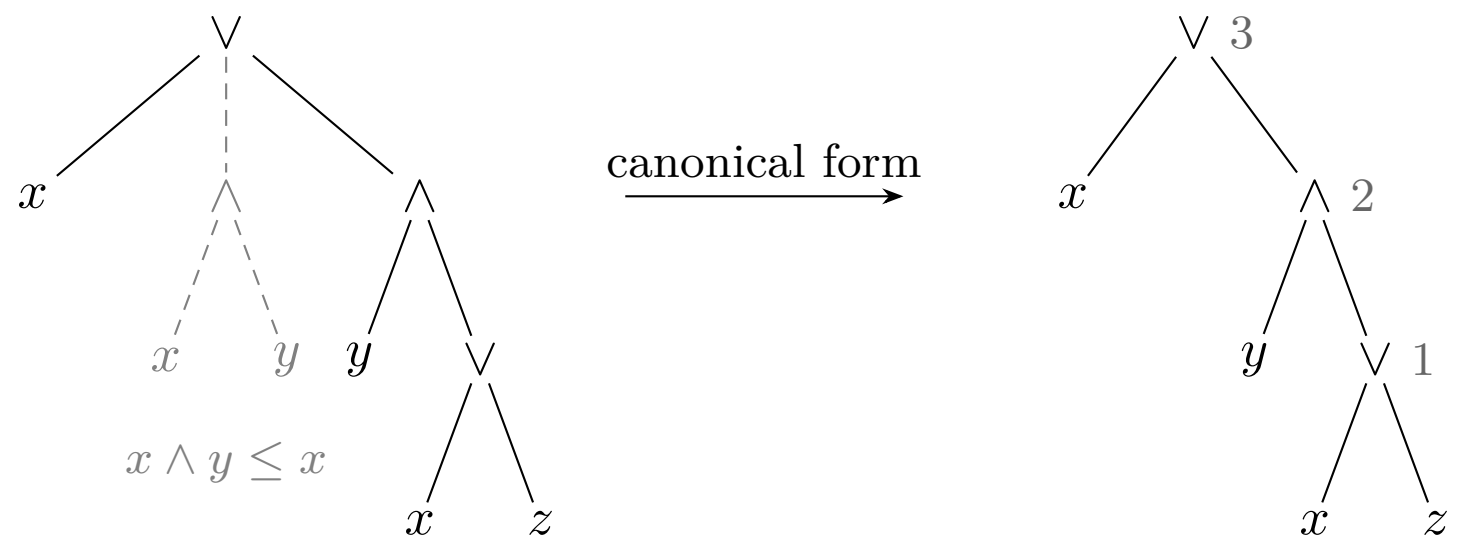

If $w = w_1 \vee \cdots \vee w_n$ with $n > 1$ is the canonical form of $w$, the elements $w_1, \dots, w_n$ are the \emph{canonical joinands} of $w$ and we write $w = w_1 \vee \cdots \vee w_n$ \emph{canonically}; if $w$ is join irreducible, its only canonical joinand is $w$ itself.  Canonical meetands are dual.  The \emph{canonical tree} of $w$ is the term tree of its canonical form (Figure~\ref{fig:canonical}); its nodes are the \emph{subelements} of $w$, each with its own canonical form, by condition (2), and the \emph{children} of a join node are its canonical joinands, those of a meet node its canonical meetands.  The \emph{complexity} $\cx(w)$ of an element is the complexity of its canonical form; so $\cx(w) = 0$ exactly for the generators, and $\cx(w) = 1 + \max_i \cx(w_i)$ when $w = w_1 \vee \cdots \vee w_n$ canonically with $n > 1$, and dually.

\begin{furthercut}
\begin{theorem}[FJN, Theorem~1.19]
\label{thm:refine}
Let $w = w_1 \vee \cdots \vee w_n$ canonically in $\FL(X)$.  If also $w = u_1 \vee \cdots \vee u_m$, then for each $i$ there is a $j$ with $w_i \leq u_j$.  Dually for meets.
\end{theorem}
\end{furthercut}

We use two consequences of the canonical form without further comment.  A canonical joinand $w_i$ of $w$ with $n > 1$ lies strictly below $w$: if $w_i = w$ then $w_j \leq w_i$ for $j \neq i$, against the antichain condition.  And an element of complexity $k + 1$ has a canonical joinand or meetand of complexity $k$, the child of largest complexity in its canonical tree.

\subsection{Embeddings and canonical joinands}
\label{sec:9.13}

\begin{lemma}[FJN, Lemma~9.13]
\label{lem:9.13}
Let $\sigma \colon \FL(X) \to \FL(Y)$ be a lattice embedding and $w = w_1 \vee \cdots \vee w_n$ canonically, $n > 1$.  If $w_1 \notin X$, then $\sigma(w_1)$ is a canonical joinand of $\sigma(w)$.  Dually for meets.
\end{lemma}

This is Lemma~1 of the notes \cite{Freese1997} and the main engine of the argument.  The requirement $w_1 \notin X$ cannot be dropped: the self-embedding $\theta_1$ of Theorem~\ref{thm:D} satisfies $\theta_1(x \vee y) = x \vee y$, whose canonical joinands are $x$ and $y$, and $\theta_1(x) = x \vee (y \wedge z)$ is not among them.

\subsection{Complexity under embeddings}
\label{sec:complexity}

\begin{thmB}
\label{thm:B}
A lattice embedding $\varphi \colon \FL(X) \to \FL(Y)$ never lowers complexity: $\cx(\varphi(w)) \geq \cx(w)$ for every $w \in \FL(X)$.
\end{thmB}

\begin{proof}
By induction on the canonical form of $w$.  A generator has complexity $0$.  Let $w = w_1 \vee \cdots \vee w_n$ canonically with $n > 1$; the case of a meet is dual.  First, $\varphi(w)$ is not a generator: a generator of $\FL(Y)$ in the image of $\varphi$ is the image of a generator of $\FL(X)$, by Corollary~\ref{cor:1.6} applied to the generating set $\varphi(X)$ of the image, so by injectivity $\varphi(w)$ is not the image of a generator.  So $\cx(\varphi(w)) \geq 1$, which settles the case in which every $w_i$ is a generator and $\cx(w) = 1$.  Otherwise choose a joinand $w_i$ that is not a generator and has $\cx(w_i) = \cx(w) - 1$.  By Lemma~\ref{lem:9.13}, $\varphi(w_i)$ is a canonical joinand of $\varphi(w)$, and a proper one: $\varphi(w_i) = \varphi(w)$ would give $w_i = w$ by injectivity, while $w_i < w$.  So $\varphi(w) = u_1 \vee \cdots \vee u_m$ canonically with $m > 1$ and $\varphi(w_i) = u_j$ for some $j$, whence
\[
\cx(\varphi(w)) = 1 + \max_l \cx(u_l) \geq 1 + \cx(\varphi(w_i)) \geq 1 + \cx(w_i) = \cx(w)
\]
by the induction hypothesis.
\end{proof}

\begin{lemma}[Finiteness]
\label{lem:finite}
Let $X$ be finite.  For every $n$, only finitely many elements of $\FL(X)$ have complexity at most $n$.  Hence an infinite sequence of elements of bounded complexity takes some value twice; more generally, an infinite sequence of $m$-tuples of elements of bounded complexity takes some tuple twice.
\end{lemma}

\begin{proof}
Let $E_n$ be the set of elements of complexity at most $n$.  By induction on $n$, $E_n$ is a finite set for all $n$.  Indeed, $E_0 = X$, and an element of complexity at most $n + 1$ is a generator, the join of its canonical joinands, or the meet of its canonical meetands, and these form a nonempty subset of $E_n$; so $|E_{n+1}| \leq |X| + 2^{|E_n| + 1}$.  The rest is the pigeonhole principle: a sequence of elements of complexity at most $n$ takes its values in the finite set $E_n$, and a sequence of $m$-tuples of them in the finite set $E_n^m$, so an infinite sequence of either takes some value twice.
\end{proof}

The complexity of an element of $\FL(X)$ is intrinsic, not relative to a generating set: a sublattice $\FL(X_0) \subseteq \FL(X)$ for $X_0 \subseteq X$ has the same canonical forms.  So a bound that Theorem~\ref{thm:B} gives in $\FL(X)$, for an injective endomorphism $\varphi$ of $\FL(X)$ and its powers, holds in $\FL(X_0)$, where $X_0$ is the finite set of generators that occur in finitely many given elements and Lemma~\ref{lem:finite} applies; the one inequality this needs, that an element of $\FL(X_0)$ has no larger complexity there than in $\FL(X)$, is Theorem~\ref{thm:B} for the inclusion.

\section{The one-step example}
\label{sec:freese-example}

This section reports the one-step example of \cite{Freese1997}, which shows that the new generators of a larger copy of $\FL(3)$ need not be subterms of the old ones, and its verification by the proof assistant.

\subsection{Where the new generators can sit}
\label{sec:where}

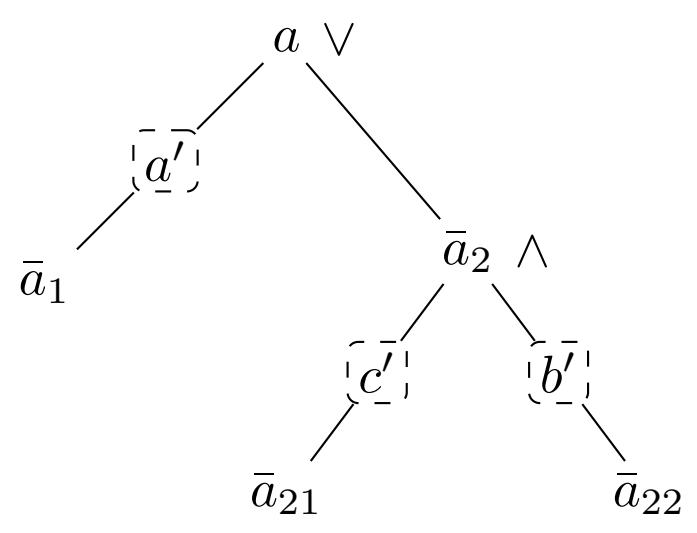
\begin{figure}[ht]
\centering
\begin{tikzpicture}[x=10mm, y=10mm]
  \node[term, label=right:{$\vee$}] (a) at (0,3) {$a$};
  \node[prime] (apr) at (-0.8,2.2) {$a\pr$};
  \node[term] (a1) at (-1.6,1.4) {$\ba{a}_1$};
  \node[term, label=right:{$\wedge$}] (a2) at (1.2,1.6) {$\ba{a}_2$};
  \node[prime] (cpr) at (0.6,0.8) {$c\pr$};
  \node[prime] (bpr) at (1.8,0.8) {$b\pr$};
  \node[term] (a21) at (0,0) {$\ba{a}_{21}$};
  \node[term] (a22) at (2.4,0) {$\ba{a}_{22}$};
  \draw[le] (apr) -- (a);
  \draw[le] (a1) -- (apr);
  \draw[le] (a2) -- (a);
  \draw[le] (a2) -- (cpr);
  \draw[le] (cpr) -- (a21);
  \draw[le] (a2) -- (bpr);
  \draw[le] (bpr) -- (a22);
  \begin{scope}[xshift=55mm]
    \node[inner sep=1.5pt, label=right:{$\vee$}] (ra) at (0,3) {$a$};
    \node[leaf] (rapr) at (-0.8,2.2) {$a\pr$};
    \node[inner sep=1.5pt, label=right:{$\wedge$}] (ra2) at (1.2,1.6) {$\ba{a}_2$};
    \node[leaf] (rcpr) at (0.6,0.8) {$c\pr$};
    \node[leaf] (rbpr) at (1.8,0.8) {$b\pr$};
    \draw (rapr) -- (ra) -- (ra2) -- (rcpr);
    \draw (ra2) -- (rbpr);
  \end{scope}
\end{tikzpicture}
\caption{Figure~1 of \cite{Freese1997}, with our bars.  Left, the tree of the term $\ba{a}_1 \vee (\ba{a}_{21} \wedge \ba{a}_{22})$ naming $a$, with each primed generator between a node and one of its children: $\ba{a}_1 \leq a\pr \leq a$, $\ba{a}_2 \leq c\pr \leq \ba{a}_{21}$ and $\ba{a}_2 \leq b\pr \leq \ba{a}_{22}$.  Right, the term $a\pr \vee (b\pr \wedge c\pr)$ naming $a$, the same shape cut off at the primes.}
\label{fig:truncation}
\end{figure}

\begin{expertread}
Suppose $\gen{a, b, c} \subseteq \gen{a\pr, b\pr, c\pr}$ with both triples free.  Then $a$ is named by a term in $a\pr$, $b\pr$, $c\pr$, and Lemma~\ref{lem:9.13}, applied to the embedding $\FL(3) \to \FL(X)$ sending the generators to $a\pr$, $b\pr$, $c\pr$, limits where $a\pr$, $b\pr$, $c\pr$ can occur in the canonical tree of $a$.  In the words of the notes \cite{Freese1997}: ``This lemma implies that if we look at the term tree of $a$ each path from $a$ to a leaf (i.e., a generator) there is a subterm $u$ and an immediate subterm $v$ of it such either $a\pr$, $b\pr$, or $c\pr$ lies between $u$ and $v$.  Another way of looking at this is that the term in $a\pr$, $b\pr$, and $c\pr$ giving $a$ is just a `truncation' of the term in the free generators giving $a$.''\footnote{The notes are quoted with their slips of spelling silently corrected; the words we supply are in brackets.}  Section~\ref{sec:truncation-theorem} states this as a theorem.

The notes continue: ``The simplest thing one might try is to write each of $a$, $b$ and $c$ as a join and try to place the $a\pr$, $b\pr$ and $c\pr$ between these joinands and the element.  But it is easy to argue that this is not possible.  The next easiest possibility is illustrated in Figure~1.''  In that figure, redrawn as Figure~\ref{fig:truncation}, $a = \ba{a}_1 \vee \ba{a}_2$ with $\ba{a}_2 = \ba{a}_{21} \wedge \ba{a}_{22}$, and the primed generators sit between a node and one of its children: $\ba{a}_1 \leq a\pr \leq a$, $\ba{a}_2 \leq c\pr \leq \ba{a}_{21}$, $\ba{a}_2 \leq b\pr \leq \ba{a}_{22}$, and cyclically for the trees of $b$ and $c$.  (The bars are ours, in the figure too: the nodes are elements that Section~\ref{sec:relations} constructs in a free lattice whose generators carry the unbarred names.)  So $a = a\pr \vee (b\pr \wedge c\pr)$, the same tree shape cut off at the primes.  Reading the three trees together, each primed generator is bounded on both sides: for instance $\ba{a}_1 \vee \ba{b}_2 \vee \ba{c}_2 \leq a\pr \leq a \wedge \ba{b}_{21} \wedge \ba{c}_{22}$.
\end{expertread}

\subsection{The nine generators and the twelve relations}
\label{sec:relations}

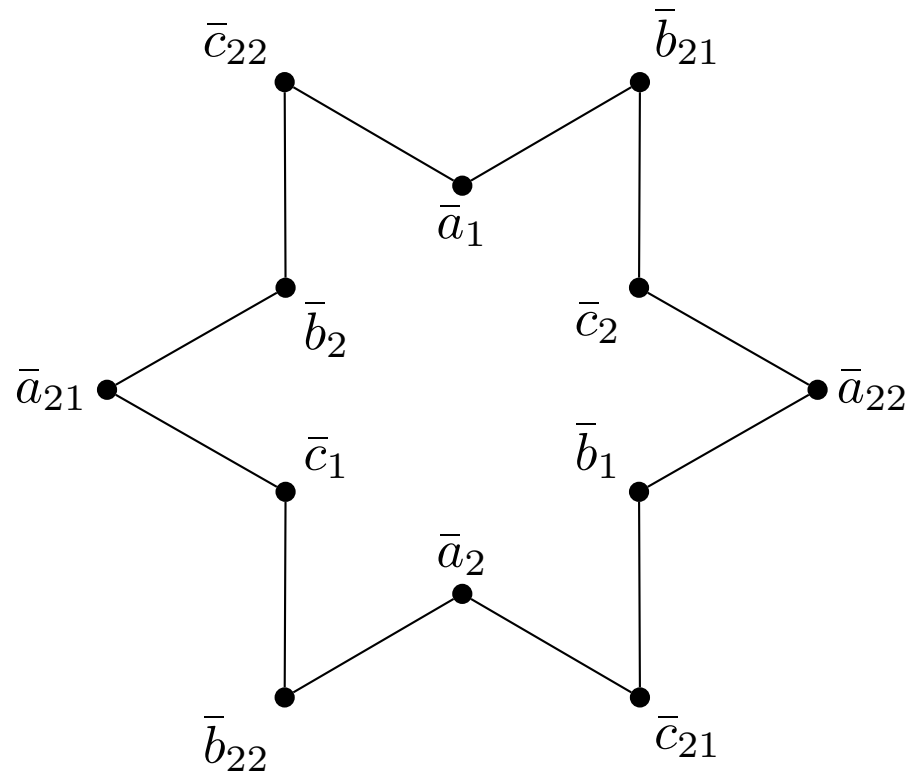
\begin{figure}[ht]
\centering
\begin{tikzpicture}[x=10mm, y=10mm]
  \node[elt, label={270:{$\ba{a}_1$}}] (a1) at (0.000,1.350) {};
  \node[elt, label={60:{$\ba{b}_{21}$}}] (b21) at (1.175,2.035) {};
  \node[elt, label={210:{$\ba{c}_2$}}] (c2) at (1.169,0.675) {};
  \node[elt, label={0:{$\ba{a}_{22}$}}] (a22) at (2.350,0.000) {};
  \node[elt, label={150:{$\ba{b}_1$}}] (b1) at (1.169,-0.675) {};
  \node[elt, label={300:{$\ba{c}_{21}$}}] (c21) at (1.175,-2.035) {};
  \node[elt, label={90:{$\ba{a}_2$}}] (a2) at (0.000,-1.350) {};
  \node[elt, label={240:{$\ba{b}_{22}$}}] (b22) at (-1.175,-2.035) {};
  \node[elt, label={30:{$\ba{c}_1$}}] (c1) at (-1.169,-0.675) {};
  \node[elt, label={180:{$\ba{a}_{21}$}}] (a21) at (-2.350,-0.000) {};
  \node[elt, label={330:{$\ba{b}_2$}}] (b2) at (-1.169,0.675) {};
  \node[elt, label={120:{$\ba{c}_{22}$}}] (c22) at (-1.175,2.035) {};
  \begin{scope}[on background layer]
    \draw[le] (a1) -- (b21);
    \draw[le] (c2) -- (b21);
    \draw[le] (c2) -- (a22);
    \draw[le] (b1) -- (a22);
    \draw[le] (b1) -- (c21);
    \draw[le] (a2) -- (c21);
    \draw[le] (a2) -- (b22);
    \draw[le] (c1) -- (b22);
    \draw[le] (c1) -- (a21);
    \draw[le] (b2) -- (a21);
    \draw[le] (b2) -- (c22);
    \draw[le] (a1) -- (c22);
  \end{scope}
\end{tikzpicture}
\caption{The twelve relations (1) to (4) of \cite{Freese1997} among the nine barred elements and the three meets $\ba{a}_2 = \ba{a}_{21} \wedge \ba{a}_{22}$, $\ba{b}_2$ and $\ba{c}_2$: each element of the inner ring lies below its two neighbors on the outer ring.  The relations close up into one cycle, and the permutation $a \mapsto b \mapsto c \mapsto a$ of the letters turns it by a third of a turn.}
\label{fig:relations}
\end{figure}

\begin{expertread}
The notes work in the free lattice on nine generators
\[
a_1, a_{21}, a_{22}, \quad b_1, b_{21}, b_{22}, \quad c_1, c_{21}, c_{22},
\]
and look for elements $\ba{a}_1, \ba{a}_{21}, \ba{a}_{22}, \dots, \ba{c}_{22}$ of it, one for each generator, satisfying the following twelve relations, which the placement of Figure~\ref{fig:truncation} forces (the bars are ours; the notes write $a_1$ for both the generator and its image):
\begin{align}
\ba{a}_{21} &\geq \ba{c}_1, & \ba{b}_{21} &\geq \ba{a}_1, & \ba{c}_{21} &\geq \ba{b}_1, \tag{1}\\
\ba{a}_{22} &\geq \ba{b}_1, & \ba{b}_{22} &\geq \ba{c}_1, & \ba{c}_{22} &\geq \ba{a}_1, \tag{2}\\
\ba{a}_{21} &\geq \ba{b}_{21} \wedge \ba{b}_{22}, & \ba{b}_{21} &\geq \ba{c}_{21} \wedge \ba{c}_{22}, & \ba{c}_{21} &\geq \ba{a}_{21} \wedge \ba{a}_{22}, \tag{3}\\
\ba{a}_{22} &\geq \ba{c}_{21} \wedge \ba{c}_{22}, & \ba{b}_{22} &\geq \ba{a}_{21} \wedge \ba{a}_{22}, & \ba{c}_{22} &\geq \ba{b}_{21} \wedge \ba{b}_{22}.  \tag{4}
\end{align}
Figure~\ref{fig:relations} draws them.  A Lisp program of the second author's finds solutions, and the notes display one:
{\footnotesize
\begin{align*}
\ba{a}_1 ={}& a_1 \vee [(a_1 \vee b_{21} \vee c_{22}) \wedge (c_1 \vee a_{21} \vee b_{22}) \wedge (a_1 \vee b_1 \vee a_{22} \vee b_{21} \vee c_{21})\\
&\qquad \wedge (b_1 \vee c_1 \vee a_{22} \vee b_{22} \vee c_{21})\\
&\qquad \wedge (a_1 \vee b_1 \vee c_1 \vee a_{21} \vee b_{21} \vee c_{21}) \wedge (a_1 \vee b_1 \vee c_1 \vee a_{22} \vee b_{22} \vee c_{22})]\\
&\vee [(a_1 \vee b_{21} \vee c_{22}) \wedge (b_1 \vee a_{22} \vee c_{21}) \wedge (a_1 \vee c_1 \vee a_{21} \vee b_{22} \vee c_{22})\\
&\qquad \wedge (b_1 \vee c_1 \vee a_{21} \vee b_{22} \vee c_{21}) \wedge (a_1 \vee b_1 \vee c_1 \vee a_{21} \vee b_{21} \vee c_{21})\\
&\qquad \wedge (a_1 \vee b_1 \vee c_1 \vee a_{22} \vee b_{22} \vee c_{22})],\\
\ba{a}_{21} ={}& a_{21} \vee [(a_{21} \vee b_{22} \vee c_1) \wedge (a_1 \vee a_{21} \vee b_{21} \vee c_1 \vee c_{22}) \wedge (a_{22} \vee b_1 \vee b_{22} \vee c_1 \vee c_{21})\\
&\qquad \wedge (a_1 \vee a_{21} \vee b_1 \vee b_{21} \vee c_1 \vee c_{21}) \wedge (a_1 \vee a_{22} \vee b_1 \vee b_{21} \vee c_1 \vee c_{21})\\
&\qquad \wedge (a_1 \vee a_{22} \vee b_1 \vee b_{21} \vee c_1 \vee c_{22}) \wedge (a_1 \vee a_{22} \vee b_1 \vee b_{22} \vee c_1 \vee c_{22})],\\
\ba{a}_{22} ={}& a_{22} \vee [(a_{22} \vee b_1 \vee c_{21}) \wedge (a_1 \vee a_{22} \vee b_1 \vee b_{21} \vee c_{22}) \wedge (a_{21} \vee b_1 \vee b_{22} \vee c_1 \vee c_{21})\\
&\qquad \wedge (a_1 \vee a_{21} \vee b_1 \vee b_{21} \vee c_1 \vee c_{21}) \wedge (a_1 \vee a_{21} \vee b_1 \vee b_{21} \vee c_1 \vee c_{22})\\
&\qquad \wedge (a_1 \vee a_{21} \vee b_1 \vee b_{22} \vee c_1 \vee c_{22}) \wedge (a_1 \vee a_{22} \vee b_1 \vee b_{22} \vee c_1 \vee c_{22})],
\end{align*}}%
with $\ba{b}_1, \ba{b}_{21}, \ba{b}_{22}$ and $\ba{c}_1, \ba{c}_{21}, \ba{c}_{22}$ obtained by applying the permutation $a \mapsto b \mapsto c \mapsto a$ of the letters once and twice.  Then
\[
a = \ba{a}_1 \vee (\ba{a}_{21} \wedge \ba{a}_{22}), \qquad a\pr = a \wedge \ba{b}_{21} \wedge \ba{c}_{22} = a \wedge (a_1 \vee b_{21} \vee c_{22}),
\]
with $b, c, b\pr, c\pr$ again by the permutation.  The notes claim that both $\{a, b, c\}$ and $\{a\pr, b\pr, c\pr\}$ generate sublattices isomorphic to $\FL(3)$, the latter properly containing the former, and that ``we can choose any $a\pr \in a \wedge b_{21} \wedge c_{22} / a_1 \vee b_2 \vee c_2$ and independently choose similar $b\pr$ and $c\pr$.  Thus $a\pr$, $b\pr$, and $c\pr$ can be chosen to be arbitrarily complex'' (here $u/v$ is the interval $[v, u]$, and the $b_2$ and $c_2$ of the notes are $\ba{b}_2 = \ba{b}_{21} \wedge \ba{b}_{22}$ and $\ba{c}_2 = \ba{c}_{21} \wedge \ba{c}_{22}$).  A footnote records a piece of luck: Figure~\ref{fig:truncation} needs $a \geq a\pr \geq \ba{b}_2$, which the definition of $a$ does not guarantee, and ``the first solution to (1)--(4) does not satisfy this but the second one does.''
\end{expertread}

\subsection{The verification}
\label{sec:verification}

The proof assistant checks the following, by Whitman decisions on explicit terms in the free lattice on nine generators (items (i) to (iv) in the module \agda{DeansProblem.\allowbreak Examples.\allowbreak FreeseStep} of the development described in Section~\ref{sec:formal}, item (v) in \agda{DeansProblem.\allowbreak Examples.\allowbreak FreeseBox}).

\begin{enumerate}[label=(\roman*)]
\item The printed $\ba{a}_1, \ba{a}_{21}, \ba{a}_{22}$ and their images under the permutation satisfy all twelve relations; four decisions, the other eight being the images of those under the permutation, which is an automorphism of the free lattice.
\item The triples $(a, b, c)$ and $(a\pr, b\pr, c\pr)$ each freely generate (Corollary~\ref{cor:1.13}, six decisions each).
\item $a = a\pr \vee (b\pr \wedge c\pr)$, and symmetrically; hence $\gen{a, b, c} \subseteq \gen{a\pr, b\pr, c\pr}$.
\item $a\pr \notin \gen{a, b, c}$; hence the containment is proper.
\item The footnote's inequalities $\ba{b}_2 \leq a$ and $\ba{c}_2 \leq a$, decided on the canonical form of $a$\ (Section~\ref{sec:verdict}).
\end{enumerate}

Item (iv) needs no search.  Write $\rho$ for the self-embedding of $\FL(3)$ with $\rho(x) = x \vee (y \wedge z)$, $\rho(y) = y \vee (z \wedge x)$, $\rho(z) = z \vee (x \wedge y)$, and $\iota$ for the embedding $\FL(3) \to \FL(9)$ sending $x, y, z$ to $a\pr, b\pr, c\pr$.  By (iii), $a = \iota(\rho(x))$ and so on, so every element of $\gen{a, b, c}$ has the form $\iota(\rho(t))$ for a term $t$.  If $a\pr = \iota(x)$ were such an element, injectivity of $\iota$ would give $x = \rho(t) \in \gen{\rho(x), \rho(y), \rho(z)}$, so by Corollary~\ref{cor:1.6} the generator $x$ would be one of $\rho(x), \rho(y), \rho(z)$, which three Whitman decisions in $\FL(3)$ refute.

\begin{fullonly}
The decisions run through the table of Section~\ref{sec:what-checked}, the form of Whitman's procedure that decides each pair of subterms once, and on smaller terms equal to the defined ones: the canonical form $\hat{a}$ of $a$, computed once by the canonicalization algorithm and then written out and checked equal to $a$ by one decision, and the simplification $a \wedge (a_1 \vee b_{21} \vee c_{22})$ of $a\pr$ given in the notes, with $\hat{a}$ in place of $a$, checked equal to $a\pr$; the cases of $b$ and $c$ are the images under the permutation.  Free generation and the equations of (iii) transfer along equality.  The smaller terms were once needed and are now a choice.  Through Whitman's recursion, which is exponential in the worst case, the decisions on the defining terms are costly: a module deciding free generation of $(a\pr, b\pr, c\pr)$ exhausts the heap of $6$ GiB allowed to the build, and one deciding the equation $a = a\pr \vee (b\pr \wedge c\pr)$ of (iii) alone checks in $11.8$ seconds at a peak of $5.7$ GiB.  Through the table, the module deciding that equation on the defining terms checks in $2.5$ seconds, and one deciding free generation of $(a, b, c)$ or of $(a\pr, b\pr, c\pr)$ in $1.9$ or $3.3$ seconds; about $1.6$ seconds of each is the loading of the modules it imports.  These modules are not part of the release; they were measured at its earlier commit \texttt{227ffdc} (Section~\ref{sec:cost}).  The smaller terms stay because they are cheaper still, and because the module of item (v) builds on $\hat{a}$.  In the check of the release, the module of items (i) to (iv) takes $1.4$ seconds (Section~\ref{sec:cost}).
\end{fullonly}

\begin{furthercut}
\begin{expertread}
For the record, the canonical form $\hat{a}$ of $a$ is the following.\footnote{The notes display the canonical form of $a$ with two meetands of its third joinand printed differently; the form here is the one computed from the printed $\ba{a}$'s and checked equal to $a$, and the claims of the notes, which concern $a$ as defined, do not depend on the display.}
{\footnotesize
\begin{align*}
\hat{a} ={}& a_1 \vee [(a_1 \vee b_{21} \vee c_{22}) \wedge (c_1 \vee a_{21} \vee b_{22}) \wedge (a_1 \vee b_1 \vee a_{22} \vee b_{21} \vee c_{21})\\
&\qquad \wedge (b_1 \vee c_1 \vee a_{22} \vee b_{22} \vee c_{21})\\
&\qquad \wedge (a_1 \vee b_1 \vee c_1 \vee a_{21} \vee b_{21} \vee c_{21}) \wedge (a_1 \vee b_1 \vee c_1 \vee a_{22} \vee b_{22} \vee c_{22})]\\
&\vee [(a_1 \vee b_{21} \vee c_{22}) \wedge (b_1 \vee a_{22} \vee c_{21}) \wedge (a_1 \vee c_1 \vee a_{21} \vee b_{22} \vee c_{22})\\
&\qquad \wedge (b_1 \vee c_1 \vee a_{21} \vee b_{22} \vee c_{21}) \wedge (a_1 \vee b_1 \vee c_1 \vee a_{21} \vee b_{21} \vee c_{21})\\
&\qquad \wedge (a_1 \vee b_1 \vee c_1 \vee a_{22} \vee b_{22} \vee c_{22})]\\
&\vee [(a_{21} \vee b_{22} \vee c_1) \wedge (a_1 \vee a_{21} \vee b_{21} \vee c_1 \vee c_{22}) \wedge (a_1 \vee a_{21} \vee b_1 \vee b_{21} \vee c_1 \vee c_{21})\\
&\qquad \wedge (a_{22} \vee b_1 \vee c_{21})\\
&\qquad \wedge (a_1 \vee a_{22} \vee b_1 \vee b_{21} \vee c_{22}) \wedge (a_1 \vee a_{22} \vee b_1 \vee b_{22} \vee c_1 \vee c_{22})].
\end{align*}}%
Its first three joinands, $a_1$ and two meets, are the joinands of $\ba{a}_1$; the third meet is new.  That the term is canonical, that is, that the four conditions of Theorem~\ref{thm:canonical} hold at each of its nodes, is itself decided, by five hundred and forty Whitman decisions, in the module \agda{DeansProblem.\allowbreak Examples.\allowbreak FreeseBox}.
\end{expertread}
\end{furthercut}

\subsection{The step seen from inside}
\label{sec:inside}

Through the isomorphism $\iota \colon \FL(3) \to \gen{a\pr, b\pr, c\pr}$, the sublattice $\gen{a, b, c}$ is the image of the self-embedding $\rho$ above.  So the one-step example, read inside the upper copy of $\FL(3)$, is the substitution
\[
\rho \colon x \mapsto x \vee (y \wedge z), \qquad y \mapsto y \vee (z \wedge x), \qquad z \mapsto z \vee (x \wedge y),
\]
whose three values freely generate (six Whitman decisions on terms with at most six generator occurrences) and miss the generator $x$ (three decisions), so that $\rho(\FL(3))$ is a proper copy of $\FL(3)$ in $\FL(3)$ (the module \agda{DeansProblem.\allowbreak Examples.\allowbreak SelfEmbedding}).  Iterating $\rho$ gives the descending chain of Section~\ref{sec:shape}, $\rho^{k+1}(\FL(3)) \subsetneq \rho^{k}(\FL(3))$ for every $k$.  What the nine-generator setting adds is the ambient lattice, in which the new generators are not subterms of the old: $a\pr$ lies between $\ba{a}_1$, the join of three of the four canonical joinands of $a$, and $a$ itself (by relations (1) and (2), $\ba{a}_1 \leq \ba{b}_{21}$ and $\ba{a}_1 \leq \ba{c}_{22}$), it is not in $\gen{a, b, c}$, and the notes add that it may be chosen as complex as one likes.  That is the phenomenon any chain must repeat, and Section~\ref{sec:theorem-a} shows that it must repeat it with ever new embeddings.

\section{Chains as sequences of self-embeddings}
\label{sec:chains}

\begin{definition}
\label{def:chain}
A \emph{chain} (for Dean's problem) in a free lattice $\FL(X)$ is a sequence $L_0 \subsetneq L_1 \subsetneq L_2 \subsetneq \cdots$ of sublattices of $\FL(X)$, each isomorphic to $\FL(3)$; the $L_k$ are its \emph{links}, and $L_0$ is its \emph{bottom}.  Unqualified, a chain is one in this sense; a chain of other sublattices is named by its links, as in a chain of $3$-generated sublattices.  We fix isomorphisms $\varphi_k \colon \FL(3) \to L_k$ and write $\theta_k = \varphi_{k+1}^{-1} \circ \varphi_k$ for the embedding of $\FL(3)$ into itself that reflects the inclusion $L_k \subseteq L_{k+1}$.
\end{definition}

\begin{figure}[ht]
\centering
\begin{furthercut}
\begin{tikzpicture}[x=10mm, y=10mm, every node/.append style={inner sep=2pt}]
  \foreach \k in {0,1,2} {
    \node (F\k) at (2.6*\k, 1.8) {$\FL(3)$};
    \node (L\k) at (2.6*\k, 0) {$L_{\k}$};
    \draw[->] (F\k) -- node[left] {$\varphi_{\k}$} (L\k);
  }
  \node (Fd) at (7.8,1.8) {$\cdots$};
  \node (Ld) at (7.8,0) {$\cdots$};
  \node[anchor=west] (U) at (8.35,0) {$\subseteq \FL(X)$};
  \foreach \j/\k in {0/1, 1/2} {
    \draw[->] (F\j) -- node[above] {$\theta_{\j}$} (F\k);
    \node at ($(L\j)!0.5!(L\k)$) {$\subsetneq$};
  }
  \draw[->] (F2) -- node[above] {$\theta_2$} (Fd);
  \node at ($(L2)!0.5!(Ld)$) {$\subsetneq$};
  \foreach \j/\k in {0/1, 1/2, 2/d} {
    \draw[->, densely dashed] (L\j.south) to[bend right=40] node[below] {$\Phi$} (L\k.south);
  }
\end{tikzpicture}
\caption{Definition~\ref{def:chain} as a diagram: the isomorphisms $\varphi_k$ of $\FL(3)$ onto the links and the proper self-embeddings $\theta_k = \varphi_{k+1}^{-1} \circ \varphi_k$, each square commuting, which makes the $\varphi_k$ compatible with the $\theta_k$; the union of the links is the direct limit of the top row, its maps $\iota_k$ the $\varphi_k$ followed by the inclusions, and{} Proposition~\ref{prop:reformulation} reads the chain off the top row.  The dashed arrows are the self-similar case, one injective endomorphism $\Phi$ of $\FL(X)$ carrying each link onto the next, which Corollary~\ref{cor:A3} rules out.}
\label{fig:colimit}
\end{furthercut}
\end{figure}
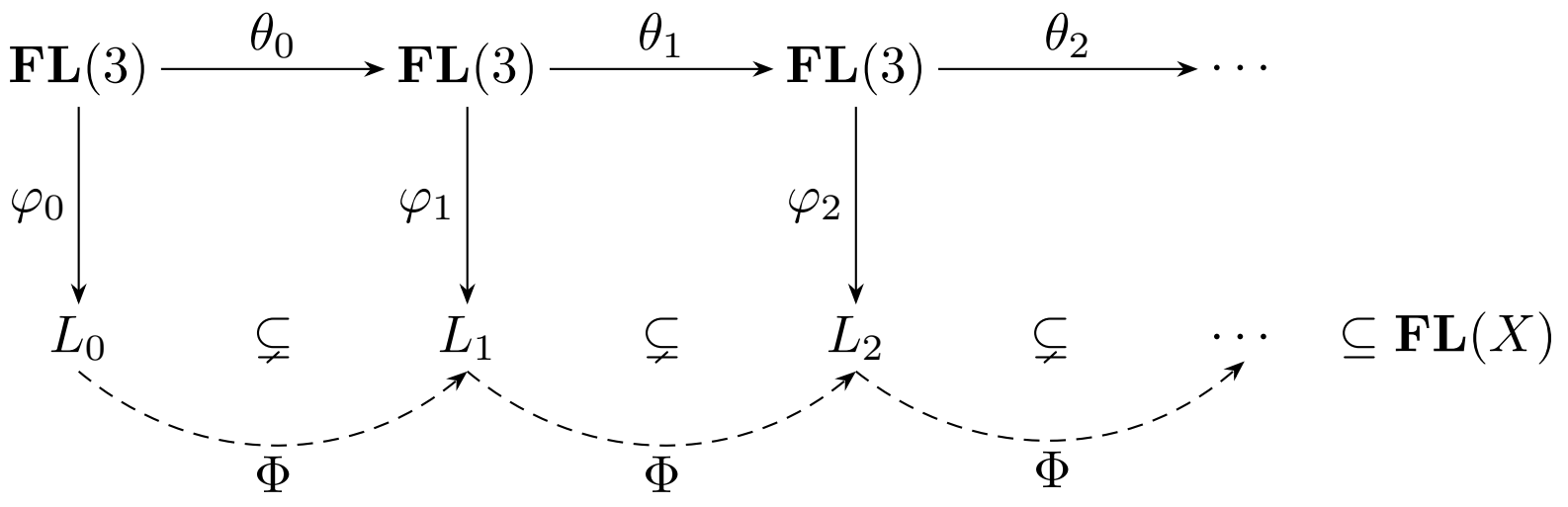

A sequence of homomorphisms $\psi_k$ from $\FL(3)$ into a lattice is \emph{compatible} with the sequence $\theta_k$ if $\psi_{k+1} \circ \theta_k = \psi_k$ for all $k$.  The sequence $\varphi_k$ in the definition above (as maps into $\FL(X)$) is compatible with the $\theta_k$.

The \emph{direct limit} of a sequence of homomorphisms $\theta_k \colon \FL(3) \to \FL(3)$ is a lattice $L$ with homomorphisms $\iota_k \colon \FL(3) \to L$ compatible with the $\theta_k$ such that every element of $L$ is some $\iota_k(u)$, and $\iota_j(u) = \iota_k(v)$ with $j \leq k$ only if $(\theta_{m-1} \circ \cdots \circ \theta_j)(u) = (\theta_{m-1} \circ \cdots \circ \theta_k)(v)$ for some $m \geq k$.
In words, $L$ is made of one copy of $\FL(3)$ for each $k$, glued so that $u$ in the $k$-th copy is identified with $\theta_k(u)$ in the next, and nothing is identified that these identifications do not force.  Figure~\ref{fig:colimit} draws the case of a chain, where the union $L$ of the links is the direct limit, its maps $\iota_k$ being the $\varphi_k$ followed by the inclusions $L_k \subseteq L$ (Proposition~\ref{prop:reformulation}(a)).

Each $\theta_k$ is injective, as a composition of injective maps, and not onto, so a chain gives a sequence of proper self-embeddings of $\FL(3)$.  Conversely, if a free lattice contains a copy of the direct limit of a sequence of proper self-embeddings of $\FL(3)$, then the sequence yields a chain in it.

\begin{furthercut}
\begin{proposition}
\label{prop:reformulation}
Let $k$ range over $\mathbb{N}$.
\begin{enumerate}[label=(\alph*)]
\item If $L_0 \subsetneq L_1 \subsetneq \cdots$ is a chain in $\FL(X)$, then the union $L = \bigcup_k L_k$ is the direct limit of the sequence of maps $\theta_k \colon \FL(3) \to \FL(3)$ and it embeds in $\FL(X)$.
\item Conversely, if $\theta_k$ are proper self-embeddings of $\FL(3)$ whose direct limit $L$ has an embedding $e$ into a free lattice $\FL(X)$, then the images $L_k = e(\iota_k(\FL(3)))$ of the canonical maps $\iota_k \colon \FL(3) \to L$ form a chain in $\FL(X)$.
\end{enumerate}
\end{proposition}

\begin{proof}
(a) The union $L$ is a sublattice of $\FL(X)$, since any two of its elements lie in some $L_k$, which holds their join and meet.  The maps $\varphi_k$ into $L$ are compatible with the $\theta_k$, and two elements with one image in $L$ are identified at a finite stage by injectivity of the $\varphi_k$ (for $j \leq k$, $\varphi_j(u) = \varphi_k(v)$ gives $(\theta_{k-1} \circ \cdots \circ \theta_j)(u) = v$, the case $m = k$); thus, $L$ is the direct limit.

(b) The canonical maps $\iota_k$ are injective because each $\theta_k$ is injective, and $\iota_k = \iota_{k+1} \circ \theta_k$; so each $L_k$ is isomorphic to $\FL(3)$ and $L_k \subseteq L_{k+1}$.  The containment is proper: if $L_{k+1} \subseteq L_k$, then each $e(\iota_{k+1}(w))$ would be some $e(\iota_k(u)) = e(\iota_{k+1}(\theta_k(u)))$, so $w = \theta_k(u)$, since $e \circ \iota_{k+1}$ is injective, and $\theta_k$ would be onto.
\end{proof}
\end{furthercut}

So Dean's problem is a question about sequences of self-embeddings of $\FL(3)$, with the ambient free lattice entering only through the limit (Figure~\ref{fig:colimit}).  We formally prove Proposition~\ref{prop:reformulation} in Agda too, in the module \agda{DeansProblem.Chain.Limit}, which imports the \agda{agda-lattices} library's direct limits.  However, our Agda proof of Theorem~\ref{thm:H} does not pass through the limit; it records a chain inside $\FL(X)$ directly, as we now describe.

A \emph{step} consists of the following data:
\begin{enumerate}[label=(\roman*)]
\item two free triples $\sigma = (a, b, c)$ and $\tau = (a\pr, b\pr, c\pr)$ in $\FL(X)$;
\item terms $p_x, p_y, p_z$ in $x, y, z$ that express $\sigma$ in terms of $\tau$: $a = p_x(a\pr, b\pr, c\pr)$, $b = p_y(a\pr, b\pr, c\pr)$ and $c = p_z(a\pr, b\pr, c\pr)$;
\item a generator of $\FL(3)$ that is none of $p_x, p_y, p_z$.
\end{enumerate}
The terms define a self-embedding $\theta \colon x \mapsto p_x, y \mapsto p_y, z \mapsto p_z$ of $\FL(3)$ (injective because $\sigma$ is free).  By (ii), $\gen{\sigma} \subseteq \gen{\tau}$.  By (iii) and Corollary~\ref{cor:1.6}, $\theta$ is not onto.  Since $\tau$ is free, it follows that $\gen{\sigma} \subsetneq \gen{\tau}$.  If $\sigma$ and $\tau$ are the triples of $L_k$ and $L_{k+1}$, respectively (that is, the images of $x, y, z$ under $\varphi_k$ and $\varphi_{k+1}$), then $\theta$ is the $\theta_k$ of Definition~\ref{def:chain}.

A chain is recorded as a sequence of free triples indexed by $\mathbb{N}$, the \mbox{$k$-th} generating the link $L_k$, with such terms and such a generator for each consecutive pair; since the two records are equivalent (as shown below), such a sequence is called a chain too, and its triples are its links.  Whether a given triple of terms misses a generator is one Whitman decision per term, which is how the one-step example is proved proper in Section~\ref{sec:verification}.

We formalize Definition~\ref{def:chain} in Agda as the \agda{CopyChain X} record type in the \agda{DeansProblem.\allowbreak Chain.\allowbreak Copies} module of the companion release \agda{agda-deans-problem}.  This record type defines a chain as a sequence of embeddings $\varphi_k$ of $\FL(3)$ in $\FL(X)$ via the subalgebra relation of the \agda{agda-algebras} library \cite{agda-algebras}, with the extra constraint that the images ascend strictly.  A sublattice isomorphic to $\FL(3)$ is the image of such an embedding, and the image of one is such a sublattice, so the record is the definition read through the inclusions.  Strictness is the negation of the reverse containment, which, given the containment, is the negation of equality; it does not ask for an element of $L_{k+1}$ outside $L_k$.  The two records are equivalent, with no hypothesis on $X$ (\agda{copyChain→chain} and \agda{chain→copyChain}).

From a chain of sublattices, the $k$-th triple is $\varphi_k(x), \varphi_k(y), \varphi_k(z)$, free by Corollary~\ref{cor:1.13} since $\varphi_k$ is injective, and the containment $L_k \subseteq L_{k+1}$ at the generators gives the terms $p_x, p_y, p_z$.  If each generator $g$ of $\FL(3)$ were equal to one of the terms, say $p_{i(g)}$, then every element $\varphi_{k+1}(t)$ of $L_{k+1}$ would be $\varphi_k(t\pr)$, where $t\pr$ is $t$ with each generator $g$ replaced by $i(g)$, and $L_{k+1} \subseteq L_k$; so some generator is none of the terms, and Whitman decisions in $\FL(3)$ find it.  Conversely, the homomorphisms that a chain of free triples induces are injective by Corollary~\ref{cor:1.13}, their images ascend by the terms, and strictly, since the image of the missed generator lies in $L_{k+1}$ and, by Corollary~\ref{cor:1.6}, not in $L_k$.

\begin{furthercut}
\subsection{The self-similar case}
\label{sec:self-similar}

The simplest conceivable chain would be generated by a single map: one injective endomorphism $\Phi$ of the ambient lattice $\FL(X)$ carrying each link onto the next, so that a single step $\gen{a, b, c} \subsetneq \gen{a\pr, b\pr, c\pr}$ with $\Phi(a) = a\pr$, $\Phi(b) = b\pr$, $\Phi(c) = c\pr$ would iterate to the whole chain.  We call a step \emph{self-similar} when an injective endomorphism of $\FL(X)$ sends its lower triple to its upper one coordinatewise.

\begin{proposition}
\label{prop:self-similar}
Let $\gen{\sigma} \subsetneq \gen{\tau}$ be a step in $\FL(X)$ and $\Phi$ an injective endomorphism of $\FL(X)$ with $\Phi(\sigma) = \tau$ coordinatewise.  Then the triples $\Phi^k(\sigma)$, $k \geq 0$, form a chain.
\end{proposition}

\begin{proof}
Each $\Phi^k$ is injective, and an injective homomorphism carries a free triple to a free triple: the conditions of Corollary~\ref{cor:1.13} are inequalities and non-inequalities among the elements, which an embedding preserves and reflects.  The terms $p$ that express $\sigma$ in $\tau$ express $\Phi^k(\sigma)$ in $\Phi^k(\tau) = \Phi^{k+1}(\sigma)$, since $\Phi^k$ is a homomorphism, and the generator missed by the terms is the same at every step.  So each consecutive pair is a step.
\end{proof}

For the one-step example this would be an injective endomorphism $\Phi$ of $\FL(9)$ solving the three equations $\Phi(a) = a\pr$, $\Phi(b) = b\pr$, $\Phi(c) = c\pr$ in the nine unknown elements $\Phi(a_1), \dots, \Phi(c_{22})$, and the initial plan was to search for one by computer.  Theorem~\ref{thm:A} shows that no injective endomorphism of any free lattice solves any such system, so that the self-similar case never occurs.
\end{furthercut}

\begin{furthercut}
\expandafter\subsection{The variant without freeness}
\label{sec:takahasi}

Dropping the requirement that the sublattices be free leaves a question of its own: does every ascending chain of $3$-generated sublattices of a free lattice stabilize?  A positive answer would have answered Dean's question too, since a strictly ascending chain of copies of $\FL(3)$ is a strictly ascending chain of $3$-generated sublattices.  It is the lattice analogue of Takahasi's theorem for free groups, that an ascending chain of subgroups of bounded rank in a free group stabilizes, which Takahasi \cite{Takahasi1950} and, independently, Higman \cite{Higman1951} proved, as Takahasi's later note \cite{Takahasi1951} recalls; that note proves a finiteness theorem of the same kind by bounding the length of the generators of a larger subgroup by those of a smaller one.  Nothing of the kind holds for lattices; sublattices of free lattices need not be free, and the one-step example of Section~\ref{sec:freese-example} shows that the generators of a larger copy of $\FL(3)$ can be more complex than those of a smaller one.  For two generators the question is trivial, since two elements generate a sublattice of at most four elements.  For three generators the answer is no (Theorem~\ref{thm:ladder}).

Given elements $a$, $b$ and $w$ of a lattice, we call the pair of sequences
\[
A_0 = a \wedge w, \quad B_0 = b \wedge w, \quad A_{k+1} = w \wedge (a \vee B_k), \quad B_{k+1} = w \wedge (b \vee A_k)
\]
the \emph{ladder} over $a$, $b$ and the bound $w$, and its terms \emph{rungs}; we write $A_k(w)$ and $B_k(w)$ when the bound matters.  Each rung is the value at $(a, b, w)$ of a lattice term in three variables, the same for every lattice.  By induction on $k$, $A_k$ and $B_k$ are the values at $(p, q, r) = (w, a, b)$ of $p \wedge q_k$ and $p \wedge r_k$, where $q_0 = q$, $r_0 = r$, $q_{k+1} = q \vee (p \wedge r_k)$ and $r_{k+1} = r \vee (p \wedge q_k)$ are the lattice words by which Kearnes and Kiss write meet semidistributivity as an identity of meet continuous lattices \cite[Corollary~2.25]{KearnesKiss}.  Figure~\ref{fig:ladder} draws the ladder of Theorem~\ref{thm:ladder}.

\begin{figure}[ht]
\centering
\begin{tikzpicture}[x=13mm, y=8mm]
  \node[elt, label=above:{$d$}] (d) at (0,5.7) {};
  \node[elt, label=left:{$A_3$}] (A3) at (-1,4) {};
  \node[elt, label=right:{$B_3$}] (B3) at (1,4) {};
  \node[elt, label=left:{$A_2$}] (A2) at (-1,3) {};
  \node[elt, label=right:{$B_2$}] (B2) at (1,3) {};
  \node[elt, label=left:{$A_1$}] (A1) at (-1,2) {};
  \node[elt, label=right:{$B_1$}] (B1) at (1,2) {};
  \node[elt, label=below left:{$A_0$}] (A0) at (-1,1) {};
  \node[elt, label=below right:{$B_0$}] (B0) at (1,1) {};
  \node[elt, label=left:{$x$}] (x) at (-2.4,1.9) {};
  \node[elt, label=right:{$y$}] (y) at (2.4,1.9) {};
  \node[elt, label=below:{$x \wedge y$}] (xmy) at (0,0) {};
  \begin{scope}[on background layer]
    \draw[le] (xmy) -- (A0);
    \draw[le] (xmy) -- (B0);
    \draw[le] (A0) -- (x);
    \draw[le] (B0) -- (y);
    \draw[le] (A0) -- (A1);
    \draw[le] (A0) -- (B1);
    \draw[le] (B0) -- (B1);
    \draw[le] (B0) -- (A1);
    \draw[le] (A1) -- (A2);
    \draw[le] (A1) -- (B2);
    \draw[le] (B1) -- (B2);
    \draw[le] (B1) -- (A2);
    \draw[le] (A2) -- (A3);
    \draw[le] (A2) -- (B3);
    \draw[le] (B2) -- (B3);
    \draw[le] (B2) -- (A3);
    \draw[le, densely dotted] (A3) -- (d);
    \draw[le, densely dotted] (B3) -- (d);
  \end{scope}
\end{tikzpicture}
\caption{The ladder over $x$, $y$ and $d = (x \vee z) \wedge (y \vee z)$ in $\FL(3)$: its first four rungs on each side, with $x \wedge y$, $x$, $y$ and $d$, ordered as in $\FL(3)$.  The diagram shows the order on these twelve elements alone, not the sublattice they generate; joins and meets such as $A_0 \vee B_0$ and $A_1 \wedge B_1$ are not drawn, and each line joins an element to one above it among the twelve.  The dotted lines stand for the rungs above $A_3$ and $B_3$.  Each side ascends and each rung lies below the next rung of the other side (Lemma~\ref{lem:ladder}); neither side ever stops (Theorem~\ref{thm:ladder}).}
\label{fig:ladder}
\end{figure}
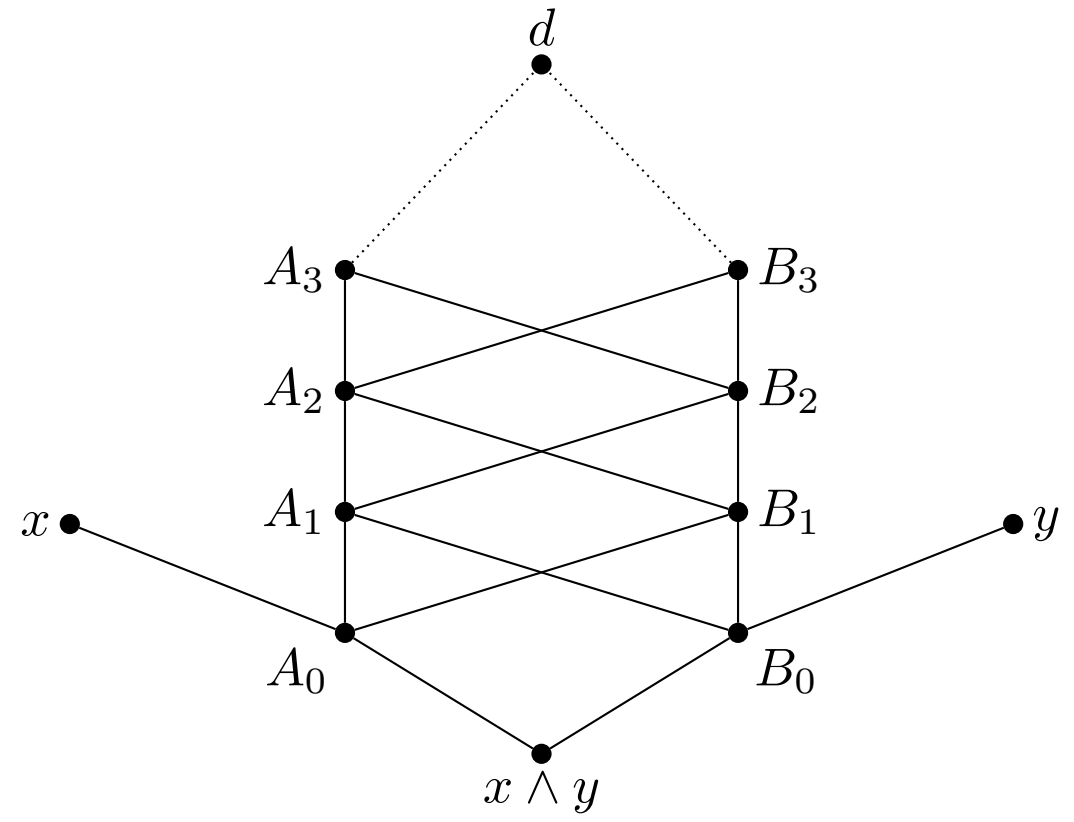

\begin{lemma}
\label{lem:ladder}
In every lattice, every rung lies below $w$; each side ascends, $A_k \leq A_{k+1}$ and $B_k \leq B_{k+1}$; and each rung lies below the next rung of the other side, $B_k \leq A_{k+1}$ and $A_k \leq B_{k+1}$.  If $A_k(w) \vee B_k(w) \leq v \leq w$, then $A_j(v) = A_j(w)$ and $B_j(v) = B_j(w)$ for every $j \leq k$.
\end{lemma}

\begin{proof}
Every rung is a meet with $w$.  $A_0 = a \wedge w \leq w \wedge (a \vee B_0) = A_1$, and if $B_{k-1} \leq B_k$ then $A_k = w \wedge (a \vee B_{k-1}) \leq w \wedge (a \vee B_k) = A_{k+1}$; symmetrically for $B$, so both sides ascend by induction.  And $B_k \leq w \wedge (a \vee B_k) = A_{k+1}$, symmetrically $A_k \leq B_{k+1}$.  For the last claim, induct on $j$.  First $A_0(v) = a \wedge v \leq a \wedge w = A_0(w)$, and $A_0(w) \leq a \wedge v$ because $A_0(w) \leq A_k(w) \leq v$; likewise $B_0(v) = B_0(w)$.  If the rungs agree at $j$ and $j + 1 \leq k$, then $A_{j+1}(v) = v \wedge (a \vee B_j(w))$ lies below $w \wedge (a \vee B_j(w)) = A_{j+1}(w)$, and above it, since $A_{j+1}(w) \leq A_k(w) \leq v$ and $A_{j+1}(w) \leq a \vee B_j(w)$; likewise for $B_{j+1}$.
\end{proof}

\begin{expertread}
\begin{theorem}
\label{thm:ladder}
In $\FL(3)$ let $d = (x \vee z) \wedge (y \vee z)$, let $A_k$ and $B_k$ be the rungs of the ladder over $x$, $y$ and $d$, and put
\[
L_{2k} = \gen{x, y, A_k \vee B_k}, \qquad L_{2k+1} = \gen{x, y, A_{k+1}}.
\]
Then $L_0 \subsetneq L_1 \subsetneq L_2 \subsetneq \cdots$ is an infinite strictly ascending chain of $3$-generated sublattices of $\FL(3)$, none of them isomorphic to $\FL(3)$.  So every free lattice on at least three generators contains an infinite strictly ascending chain of $3$-generated sublattices, and the analogue of Takahasi's theorem fails for lattices.
\end{theorem}
\end{expertread}

\begin{proof}
Write $w_{2k} = A_k \vee B_k$ and $w_{2k+1} = A_{k+1}$, so that $L_n = \gen{x, y, w_n}$ for every $n$; each $w_n$ is the value at $(x, y, d)$ of a term $t_n$, the ladder written out, and all lie below $d$.

\emph{Ascent.}  By Lemma~\ref{lem:ladder}, $w_{2k+1} = A_{k+1}$ lies between $A_k \vee B_k$ and $d$, so the rungs up to $k$ over $(x, y, w_{2k+1})$ are those over $(x, y, d)$, and $w_{2k} = t_{2k}(x, y, w_{2k+1}) \in L_{2k+1}$.  Likewise $w_{2k+2} = A_{k+1} \vee B_{k+1}$ lies between itself and $d$, so $w_{2k+1} = A_{k+1}(w_{2k+2}) \in L_{2k+2}$, and so does $B_{k+1} = B_{k+1}(w_{2k+2})$.  Hence $L_n \subseteq L_{n+1}$ for every $n$.

\emph{Seven facts.}  Whitman's rules (Theorem~\ref{thm:whitman}) decide that $d \not\leq x \vee y$, $x \not\leq y \vee z$, $y \not\leq x \vee z$, $x \wedge d \not\leq y$, $y \wedge d \not\leq x$, $x \not\leq y$ and $y \not\leq x$.  So neither $x$ nor $y$ lies below $d$.  By induction every rung lies below $x \vee y$, since $A_{k+1} \leq x \vee B_k$ and $B_{k+1} \leq y \vee A_k$.  No rung $A_k$ lies below $y$, because $A_0 = x \wedge d \leq A_k$, and no rung $A_{k+1}$ lies below $x$, because $B_0 = y \wedge d \leq B_k \leq A_{k+1}$; symmetrically no $B_k$ lies below $x$ and no $B_{k+1}$ below $y$.

\emph{The sides ascend strictly.}  $A_1 \leq A_0 \leq x$ is false, and symmetrically so is $B_1 \leq B_0 \leq y$.  Suppose $A_{k+2} \leq A_{k+1}$; then $A_{k+2} = d \wedge (x \vee B_{k+1}) \leq x \vee B_k$, and $\W$ leaves the following four cases:
\begin{itemize}
\item $d \leq x \vee B_k \leq x \vee y$, which is false;
\item $x \vee B_{k+1} \leq x \vee B_k$;
\item $A_{k+2} \leq x$, which is false;
\item $A_{k+2} \leq B_k$, which gives $B_{k+1} \leq A_{k+2} \leq B_k$.
\end{itemize}
In either remaining case, $B_{k+1} \leq x \vee B_k$, and $\W$ applied to $B_{k+1} = d \wedge (y \vee A_k)$ leaves only $B_{k+1} \leq B_k$.  Indeed, $d \leq x \vee B_k$ is false, $y \vee A_k \leq x \vee B_k$ would put the join prime generator $y$ below $x$ or below $B_k \leq d$, and $B_{k+1} \leq x$ is false.
So $A_{k+2} \leq A_{k+1}$ implies $B_{k+1} \leq B_k$, and symmetrically $B_{k+2} \leq B_{k+1}$ implies $A_{k+1} \leq A_k$.  Since $A_1 \not\leq A_0$ and $B_1 \not\leq B_0$, induction on $k$ gives $A_{k+1} \not\leq A_k$ and $B_{k+1} \not\leq B_k$ for every $k$.

We prove $A_{k+1} \not\leq B_k$ similarly.  For $k = 0$, $A_1 \leq B_0$ would give $A_1 \leq y$, which is false.  If $A_{k+2} \leq B_{k+1} \leq y \vee A_k$, then $\W$ leaves only $A_{k+2} \leq A_k$, since $d \leq y \vee A_k \leq x \vee y$ is false, $x \vee B_{k+1} \leq y \vee A_k$ would put the join prime generator $x$ below $y$ or below $A_k \leq d$, and $A_{k+2} \leq y$ is false.  But $A_{k+2} \leq A_k$ gives $A_{k+1} \leq A_{k+2} \leq A_k$, contradicting the fact $A_{k+1} \not\leq A_k$ that we just established.

Also, $A_{k+1} \not\leq A_k \vee B_k$, since by $\W$ that would need $d \leq A_k \vee B_k \leq x \vee y$, or $x \leq A_k \vee B_k$ (so $x \leq A_k$ or $x \leq B_k$, both below $d$), or $A_{k+1} \leq A_k$, or $A_{k+1} \leq B_k$.
And $B_{k+1} \not\leq A_{k+1}$, since that would give $B_{k+1} \leq x \vee B_k$ and so $B_{k+1} \leq B_k$, as above.

\emph{The chain ascends strictly.} We first prove $L_{2k} \subsetneq L_{2k+1}$.  Let $S$ be the set of elements $s \leq x \vee y$ such that $A_{k+1} \leq s$ implies $x \vee B_k \leq s$.  It contains $x$, $y$ and $A_k \vee B_k$, none of which lies above $A_{k+1}$ (the last because $A_{k+1} \not\leq A_k \vee B_k$, proved above), and it is closed under meets.  It is closed under joins: if $A_{k+1} = d \wedge (x \vee B_k) \leq s \vee t$ with $s, t \in S$, then by $\W$ either $d \leq s \vee t \leq x \vee y$, which is false, or $x \vee B_k \leq s \vee t$, or $A_{k+1}$ lies below $s$ or below $t$, which then lies above $x \vee B_k$.  So $L_{2k} \subseteq S$, and if $A_{k+1}$ were in $L_{2k}$ we would have $x \leq x \vee B_k \leq A_{k+1} \leq d$, which is false.
Next we show $L_{2k+1} \subsetneq L_{2k+2}$.  The same argument with the generators and the sides exchanged, for the set of elements $s \leq x \vee y$ such that $B_{k+1} \leq s$ implies $y \vee A_k \leq s$, which contains $x$, $y$ and $A_{k+1}$ (the last because $B_{k+1} \not\leq A_{k+1}$, proved above), keeps $B_{k+1}$ out of $L_{2k+1}$; and $B_{k+1} \in L_{2k+2}$ by the ascent.

\begin{expertread}
\emph{No link is a copy of $\FL(3)$.}  $x \wedge y \leq x \wedge d = A_0 \leq w_n$, so the triple $(x, y, w_n)$ fails the meet condition~\eqref{eq:meet-condition} of Corollary~\ref{cor:1.13} and does not freely generate.  If $L_n$ were isomorphic to $\FL(3)$, the images of the free generators would generate $L_n$ as in Lemma~\ref{lem:prime}, so by Corollary~\ref{cor:1.6} they would lie in the generating set $\{x, y, w_n\}$, hence be $x$, $y$ and $w_n$, which do not freely generate.
\end{expertread}

Finally, $\FL(3)$ is a sublattice of every $\FL(X)$ with $|X| \geq 3$, generated by three of the generators, and an embedding carries a strictly ascending chain of $3$-generated sublattices to another.
\end{proof}

The herringbone of \cite[Figure~8.2]{KearnesKiss} is a configuration of the same zigzag kind, which Kearnes and Kiss build from the duals of these words \cite[proof of Theorem~8.3]{KearnesKiss}.  And FJN, Figure~5.3, draws the lattice $\FL(2+2)$ freely generated by two two-element chains $a_0 < a_1$ and $b_0 < b_1$, which is infinite and embeds in a free lattice (FJN, Chapter~V, Section~2).  The homomorphism from $\FL(2+2)$ to $\FL(3)$ that extends $a_0 \mapsto x$, $a_1 \mapsto x \vee z$, $b_0 \mapsto y$ and $b_1 \mapsto y \vee z$ sends $a_1 \wedge b_1$ to $d$, so it maps the ladder over $a_0$, $b_0$ and $a_1 \wedge b_1$ onto the ladder of Theorem~\ref{thm:ladder}, and $\gen{a_0, b_0, a_1 \wedge b_1}$ onto $\gen{x, y, d}$, which contains every link.  So $\gen{a_0, b_0, a_1 \wedge b_1}$, like $\gen{x, y, d}$, is an infinite $3$-generated sublattice of a free lattice.

The companion release \agda{agda-deans-problem} checks the theorem as stated (\agda{Theorem-4·5} in \agda{DeansProblem.\allowbreak Examples.\allowbreak ChainCondition}, which writes $c$ for $d$).  There \agda{ladder} is the chain, \agda{no-link-is-a-copy} proves that no link is a copy of $\FL(3)$, and \agda{ladder-into} proves that the last clause of Theorem~\ref{thm:ladder} holds in every free lattice with three distinct generators, with no assumption of decidable equality.  The chain condition itself, ascending, stabilizing and strictly ascending chains of $n$-generated sublattices, is stated in \agda{DeansProblem.\allowbreak Chain.\allowbreak Condition}, where the record of a chain of copies of $\FL(3)$ is proved a strictly ascending chain in this sense (\agda{copyChain→\allowbreak strict}), so that the chain condition excludes such chains (\agda{chainCondition→\allowbreak no-copy-chain}).  Theorem~\ref{thm:H} is the case of the chain condition in which every link is free, and it holds where the general case fails.  What fails without freeness is not the recurrence on which Theorem~\ref{thm:C} rests; by Lemma~\ref{lem:ladder} the term $t_0$ expresses $w_0$ in every later link, so the terms expressing the bottom link in the later ones recur here too.  It is injectivity that fails, which Theorems~\ref{thm:B} to~\ref{thm:G} use throughout.
\end{furthercut}

\begin{fullonly}
\begin{computation}
\label{comp:ladder}
The links appear to be finite.  Closing each generating triple under joins and meets in Python, with a memoized Whitman procedure, finishes for the first $24$ links: $L_{2k}$ has $10 + 6k$ elements and $L_{2k+1}$ has $12 + 6k$, for $k \leq 11$, and each containment is proper by the closures alone.  We do not prove that every link is finite.  If they are, the finite $3$-generated sublattices of a free lattice have unbounded size, and even an ascending chain of finite $3$-generated sublattices need not stabilize.  The chain was found by a search that decides membership in a finitely generated sublattice of a free lattice exactly: a homomorphism $h$ from a finitely generated lattice into a free lattice is lower bounded (FJN, Corollary~2.11, Theorem~2.6 and Lemma~2.7), so $w$ lies in its image exactly when $h(\beta(w)) = w$, where $\beta(w)$ is the least element whose image lies above $w$, which FJN, Corollary~2.8, computes by a recursion through the minimal nontrivial join covers of $w$, and the search by a recursion on a term for $w$ derived from it.  The search ran inside the sublattice $\gen{x, y, d}$, which contains every $L_n$ and is therefore infinite, but which grows slowly under closure: after the first five rounds it has $9$, $16$, $21$, $24$ and $27$ elements, where $\FL(3)$ has $9$, $23$, $49$, $315$ and $20\,611$, and each of the next seven rounds adds three elements.
\end{computation}
\end{fullonly}

\section{No step is self-similar}
\label{sec:theorem-a}

\begin{thmA}
\label{thm:A}
Let $\Phi$ be an injective endomorphism of a free lattice $\FL(X)$ and $L$ a finitely generated sublattice with $L \subseteq \Phi(L)$.  Then $\Phi(L) = L$.
\end{thmA}

\begin{proof}
Let $g_1, \dots, g_n$ generate $L$.  Since $g_i \in \Phi(L)$, there is $h_i \in L$ with $\Phi(h_i) = g_i$, and $h_i = t_i(g_1, \dots, g_n)$ for a term $t_i$.  Let $\theta$ be the endomorphism of $\FL(n)$ with $\theta(x_i) = t_i$ and $\tau \colon \FL(n) \to \FL(X)$ the homomorphism with $\tau(x_i) = g_i$.  Then $\Phi(\tau(\theta(x_i))) = \Phi(t_i(g_1, \dots, g_n)) = \Phi(h_i) = g_i = \tau(x_i)$ for each $i$, so $\Phi \circ \tau \circ \theta = \tau$, and by induction $\Phi^k \circ \tau \circ \theta^k = \tau$ for every $k \geq 0$.

Fix $i$ and put $u_k = \tau(\theta^k(x_i)) \in L$.  Then $\Phi^k(u_k) = \tau(x_i) = g_i$, and $\Phi^k$ is an injective endomorphism of $\FL(X)$, so Theorem~\ref{thm:B} gives $\cx(u_k) \leq \cx(g_i)$ for every $k$.  The elements $u_k$ lie in $L \subseteq \FL(X_0)$, where $X_0$ is the finite set of generators occurring in $g_1, \dots, g_n$, so by Lemma~\ref{lem:finite} there are $a < b$ with $u_a = u_b$.  Then $g_i = \Phi^a(u_a) = \Phi^b(u_b) = \Phi^a(\Phi^{b-a}(u_b))$, and injectivity of $\Phi^a$ gives $u_a = \Phi^{b-a}(u_b) = \Phi^{b-a}(u_a)$.  Hence
\[
\Phi^{b-a}(g_i) = \Phi^{b-a}(\Phi^a(u_a)) = \Phi^a(\Phi^{b-a}(u_a)) = \Phi^a(u_a) = g_i.
\]
So each generator $g_i$ of $L$ is fixed by a positive power $\Phi^{d_i}$ of $\Phi$.  Let $D$ be a common multiple of $d_1, \dots, d_n$; then $\Phi^D$ fixes every $g_i$, hence every element of $L$, and $\Phi^D(L) = L$.  Applying $\Phi$ repeatedly to $L \subseteq \Phi(L)$ gives $L \subseteq \Phi(L) \subseteq \Phi^2(L) \subseteq \cdots \subseteq \Phi^D(L) = L$, so $\Phi(L) = L$.
\end{proof}

\begin{fullonly}
The machine-checked form (\agda{self-similar-onto} in the module \agda{DeansProblem.\allowbreak SelfSimilar}) is the case in which $L = \tau(\FL(n))$ is free, $\tau$ injective: for an injective $\Phi$ and terms $\theta$ with $\tau = \Phi \circ \tau \circ \theta$, every generator $x_j$ of $\FL(n)$ is one of the $\theta(x_i)$, so $\theta$ is onto and $\Phi(L) = \Phi(\tau(\theta(\FL(n)))) = \tau(\FL(n)) = L$.  Its proof is the one above read inside $\FL(n)$: the orbit $\theta^k(x_j)$ has complexity at most $\cx(\tau(x_j))$, by Theorem~\ref{thm:B} for the embedding $\Phi^k \circ \tau$, so it repeats, $\theta^a(x_j) = \theta^b(x_j)$ with $a < b$; $\theta$ is injective, since $\Phi \circ \tau \circ \theta$ is; so $x_j = \theta^{b-a}(x_j) = \theta(\theta^{b-a-1}(x_j))$ lies in the image of $\theta$, and Corollary~\ref{cor:1.6} makes it one of the $\theta(x_i)$.  Freeness of $L$ spares the common period.  The general form is checked too (\agda{generated-onto} in \agda{DeansProblem.\allowbreak Chain.\allowbreak Condition}), for generators with decidable equality, by a variant of the proof above that needs no common period either: since $\Phi \circ \tau \circ \theta = \tau$, $\Phi(u_{k+1}) = u_k$ for every $k$, so a repeat $u_a = u_b$ with $a < b$ gives $g_i = \Phi^a(u_a) = \Phi^a(u_b) = u_{b-a}$ and then $\Phi(g_i) = u_{b-a-1} \in L$; so $\Phi(L) \subseteq L$, and injectivity of $\Phi$ enters only through Theorem~\ref{thm:B}.
\end{fullonly}

\begin{corAone}[no self-similar step]
\label{cor:A1}
No step $\gen{\sigma} \subsetneq \gen{\tau}$ in a free lattice admits an injective endomorphism $\Phi$ with $\Phi(\sigma) = \tau$.
\end{corAone}

\begin{proof}
With $L = \gen{\sigma}$ we would have $L \subseteq \gen{\tau} = \Phi(L)$, so $\Phi(L) = L$ by Theorem~\ref{thm:A}, that is, $\gen{\tau} = \gen{\sigma}$, against the properness of the step.  (The machine-checked form, \agda{no-self-similar-step}, reads: with $\theta$ the step's terms, $\tau = \Phi \circ \tau \circ \theta$, and the theorem makes the generator the step's terms miss one of those terms.)
\end{proof}

\begin{furthercut}
\begin{corAtwo}[the one-step example]
\label{cor:A2}
The equations $\Phi(a) = a\pr$, $\Phi(b) = b\pr$, $\Phi(c) = c\pr$ of Section~\ref{sec:self-similar} have no injective solution $\Phi$ in the endomorphisms of $\FL(9)$.
\end{corAtwo}

\begin{cutcandidate}
The theorem settles every candidate at once, at no decision cost on the terms of the example.  Three candidates fail for their own reasons, which are worth seeing: the cyclic permutation $\pi$ of the letters sends $a$ to $b$, and $b \not\leq a$ by the join condition~\eqref{eq:join-condition} of the free triple $(a, b, c)$ while $a\pr \leq a$; its square sends $a$ to $c$, likewise; and the natural guess, each generator to its barred element, is not even injective, since relation (1) of Section~\ref{sec:relations} puts $\ba{c}_1$ below $\ba{a}_{21}$ while $c_1 \not\leq a_{21}$ (the module \agda{DeansProblem.\allowbreak Examples.\allowbreak FreeseEquations}).
\end{cutcandidate}

\begin{corAthree}[links]
\label{cor:A3}
In a chain $L_0 \subsetneq L_1 \subsetneq \cdots$ of finitely generated sublattices of a free lattice, no injective endomorphism $\Phi$ has $\Phi(L_j) = L_k$ with $j < k$.
\end{corAthree}

\begin{proof}
$L_j \subseteq L_k = \Phi(L_j)$ forces $\Phi(L_j) = L_j$, by Theorem~\ref{thm:A}, so $L_k = L_j$, contradicting $L_j \subsetneq L_k$.
\end{proof}

So no chain is generated by one embedding, even with a period.  Freeness of the links is not needed, so the same holds for the chains of $3$-generated sublattices of Section~\ref{sec:takahasi}, which do exist (Theorem~\ref{thm:ladder}); the checked form (\agda{no-onto-later}) says that in a strictly ascending chain no injective $\Phi$ maps $L_j$ onto a sublattice containing a later $L_k$.  The maps $\theta_k$ of Definition~\ref{def:chain} are never the powers of one map either, but that is a different claim, since $\Phi$ acts on $\FL(X)$ and $\theta_k$ on $\FL(3)$, and it follows from the bound on the transitions, Lemma~\ref{lem:bounded-transitions}, not from Corollary~\ref{cor:A3}.  If $\theta_k = \psi^{m_k + 1}$ for every $k$, the transition $\theta_{0,k} = \theta_{k-1} \circ \cdots \circ \theta_0$ of Section~\ref{sec:transitions} is $\psi^{M_k}$ with $M_k = \sum_{e < k} (m_e + 1)$, so $M_0 < M_1 < \cdots$, and by the lemma $\theta_{0,j} = \theta_{0,k}$ for some $j < k$.  Then $\psi^{M_j} \circ \psi^{M_k - M_j} = \psi^{M_j}$, and canceling the injective $\psi^{M_j} = \theta_{0,j}$ makes $\theta_{j,k} = \psi^{M_k - M_j}$ the identity, against the properness of $\theta_{j,k}$ for $j < k$.

\begin{corAfour}[periodic points]
\label{cor:A4}
Let $\theta$ be an injective endomorphism of $\FL(n)$.  The intersection $\bigcap_k \theta^k(\FL(n))$ of the descending chain of images is exactly the set of elements fixed by some positive power of $\theta$.
\end{corAfour}

\begin{proof}
Let $m$ lie in every image, $m = \theta^k(m_k)$.  By Theorem~\ref{thm:B}, $\cx(m_k) \leq \cx(m)$ for every $k$, so $m_a = m_b$ for some $a < b$ by Lemma~\ref{lem:finite}.  As in the proof of Theorem~\ref{thm:A}, $\theta^a(m_a) = \theta^b(m_b) = \theta^a(\theta^{b-a}(m_a))$ gives $m_a = \theta^{b-a}(m_a)$, and then $\theta^{b-a}(m) = \theta^{b-a}(\theta^a(m_a)) = \theta^a(\theta^{b-a}(m_a)) = \theta^a(m_a) = m$.  Conversely, if $\theta^d(m) = m$ with $d \geq 1$ then $m = \theta^{kd}(m) \in \theta^k(\FL(n))$ for every $k$.
\end{proof}

Both directions are machine-checked (\agda{image→periodic} and \agda{periodic→\allowbreak image} in \agda{DeansProblem.\allowbreak Fixed}).  We write $\Per(\theta)$ for this set, the \emph{periodic points} of $\theta$, and $\Fix(\theta) = \{w : \theta(w) = w\}$ for the \emph{fixed sublattice} of $\theta$; clearly $\Fix(\theta) \subseteq \Per(\theta)$.

\begin{remark}
\label{rem:what-survives}
Theorem~\ref{thm:A} needs neither the freeness of $L$ nor three generators, and Theorem~\ref{thm:B} is its only ingredient beyond Lemma~\ref{lem:finite}.  What survives of the self-similar idea for an arbitrary chain is a necessary condition: the terms expressing the bottom triple in the $k$-th have bounded complexity, so one proper self-embedding of $\FL(3)$ expresses the bottom triple, by the same three terms, in infinitely many later triples.  That condition is the start of the next section.
\end{remark}
\end{furthercut}

\section{Coherent chains and the fixed-sublattice question}
\label{sec:coherent}

Theorem~\ref{thm:B} gives more than Theorem~\ref{thm:A}: it gives Theorem~\ref{thm:C}, which has four forms, named in their headings and called the one-step, four-link, length-three and coherent forms.  Throughout this section fix a chain $L_0 \subsetneq L_1 \subsetneq \cdots$ in $\FL(X)$ with isomorphisms $\varphi_k \colon \FL(3) \to L_k$, as in Definition~\ref{def:chain}.

\subsection{The transitions of a chain}
\label{sec:transitions}

For $j \leq k$ let $\theta_{j,k} = \varphi_k^{-1} \circ \varphi_j$ be the self-embedding of $\FL(3)$ that reads the inclusion $L_j \subseteq L_k$; so $\theta_{k,k+1}$ is the $\theta_k$ of Definition~\ref{def:chain}.  Three properties are immediate.
\begin{enumerate}[label=(\alph*)]
\item $\varphi_k \circ \theta_{j,k} = \varphi_j$.
\item $\theta_{j,j}$ is the identity, and $\theta_{j,l} = \theta_{k,l} \circ \theta_{j,k}$ for $j \leq k \leq l$ (the cocycle identity).
\item $\theta_{j,k}$ is injective, and it is onto exactly when $j = k$: by (a) its image under $\varphi_k$ is $L_j$, which is all of $L_k$ only for $j = k$.
\end{enumerate}
We call $\theta_{j,k}$ the \emph{transition} from link $j$ to link $k$.  In the formal development a chain is given by its triples and the terms expressing each in the next (Section~\ref{sec:chains}), and $\theta_{j,k}$ is the composite of the terms of the steps $j, \dots, k-1$; its properness for $j < k$ comes from the generator the last step misses, through Corollary~\ref{cor:1.6}, since $\theta_{j,k}$ is $\theta_{j,k-1}$ followed by the terms of step $k - 1$.

\begin{lemma}
\label{lem:bounded-transitions}
For every $j$, every $k \geq j$ and every generator $g$ of $\FL(3)$, $\cx(\theta_{j,k}(g)) \leq \cx(\varphi_j(g))$.  Hence, for fixed $j$, the transitions $\theta_{j,k}$ with $k \geq j$ take only finitely many values.
\end{lemma}

\begin{proof}
By (a), $\varphi_j(g) = \varphi_k(\theta_{j,k}(g))$, and $\varphi_k$ is an embedding $\FL(3) \to \FL(X)$, so Theorem~\ref{thm:B} gives the bound.  A self-embedding is determined by the triple of images of the generators, and by Lemma~\ref{lem:finite} only finitely many triples have coordinates of complexity at most $\max_g \cx(\varphi_j(g))$.
\end{proof}

Write $B(j) = \max_g \cx(\varphi_j(g))$, and $N(j)$ for the number of triples of elements of $\FL(3)$ of complexity at most $B(j)$; the transitions $\theta_{j,k}$, $k \geq j$, take at most $N(j)$ values.

\subsection{Theorem~C: the finite forms}
\label{sec:theorem-c-finite}

\begin{thmC}[one step]
\label{thm:C}
There are $j < k$ with $\theta_{0,j} = \theta_{0,k}$.  For such $j$ and $k$, the proper self-embedding $\theta_{j,k}$ of $\FL(3)$ fixes pointwise the sublattice $M = \theta_{0,j}(\FL(3))$, which is a proper sublattice of $\FL(3)$ isomorphic to $\FL(3)$.
\end{thmC}

\begin{proof}
By Lemma~\ref{lem:bounded-transitions} with $j = 0$, among $\theta_{0,0}, \dots, \theta_{0,N(0)}$ two coincide: $\theta_{0,j} = \theta_{0,k}$ with $j < k$.  By the cocycle identity, $\theta_{j,k} \circ \theta_{0,j} = \theta_{0,k} = \theta_{0,j}$, which says that $\theta_{j,k}$ is the identity on $M = \theta_{0,j}(\FL(3))$.  $M$ is isomorphic to $\FL(3)$ because $\theta_{0,j}$ is injective; $M \neq \FL(3)$ because $j \geq 1$ (if $j = 0$ then $\theta_{0,k} = \theta_{0,0}$ would be onto, against (c)); and $\theta_{j,k}$ is not onto, by (c), since $j < k$.
\end{proof}

The proof is constructive, a finite pigeonhole among the first $N(0) + 1$ links, and the formal development follows it (\agda{chain→fixed-sublattice} in \agda{DeansProblem.\allowbreak Coherent}): the conclusion is a record of a proper self-embedding $\theta$, a generator it misses, and a free triple $\Theta$ with $\theta(\Theta(g)) = \Theta(g)$ for each generator $g$, the configuration we now name.

\begin{definition}
\label{def:coherent}
A \emph{coherent chain of length $n$} is a sequence $\Theta_0, \dots, \Theta_{n-1}$ of proper self-embeddings of $\FL(3)$ with $\Theta_k \circ \Theta_j = \Theta_j$ for all $j < k$; an \emph{infinite coherent chain} is the same with the index set $\mathbb{N}$.  The $\Theta_j$ are its \emph{links}, counted from $\Theta_0$.  The condition says that each later embedding fixes every element of the image of each earlier one.
\end{definition}

\begin{lemma}
\label{lem:coherent-images}
Let $\Theta_0, \dots, \Theta_{n-1}$ be a coherent chain and $M_j = \Theta_j(\FL(3))$.  Then $M_0 \subsetneq M_1 \subsetneq \cdots \subsetneq M_{n-1} \subsetneq \FL(3)$, and each $M_j$ is a copy of $\FL(3)$.  A coherent chain of length two is exactly the configuration of the one-step form of Theorem~\ref{thm:C}: a proper self-embedding $\Theta_1$ fixing pointwise the copy $M_0 = \Theta_0(\FL(3))$ of $\FL(3)$.
\end{lemma}

\begin{proof}
For $j < k$, $M_j = \Theta_j(\FL(3)) = \Theta_k(\Theta_j(\FL(3))) \subseteq \Theta_k(\FL(3)) = M_k$.  If $M_j = M_k$ then $\Theta_k(\Theta_j(\FL(3))) = \Theta_k(\FL(3))$, and injectivity of $\Theta_k$ gives $\Theta_j(\FL(3)) = \FL(3)$, against the properness of $\Theta_j$; likewise $M_{n-1} \neq \FL(3)$.  The last sentence is the definition read at $n = 2$.
\end{proof}

So the strictness of the ascent costs nothing in a coherent chain: it is the properness of the lower embeddings (Figure~\ref{fig:coherent}).  Theorem~\ref{thm:D} below exhibits a coherent chain of length two.

\begin{figure}[ht]
\centering
\begin{tikzpicture}[x=10mm, y=10mm]
  \foreach \k in {0,1,2} {
    \begin{scope}[xshift=40*\k mm]
      \ifnum\k>0
        \fill[region] (0.25,0.2) rectangle ({1.25+0.75*(\k-1)}, {0.85+0.5*(\k-1)});
      \fi
      \draw[link] (0,0) rectangle (3.5,2.5);
      \node[anchor=north west] at (0,2.5) {$\FL(3)$};
      \foreach \j in {0,1,2} {
        \pgfmathsetmacro{\w}{\j==\k ? 1.2 : 0.4}
        \draw[link, line width=\w pt] (0.25,0.2) rectangle (1.25+0.75*\j, 0.85+0.5*\j);
        \node[anchor=north east, inner sep=1.5pt] at (1.25+0.75*\j, 0.85+0.5*\j) {$M_{\j}$};
      }
      \node at (1.75,-0.35) {$\Theta_{\k}$};
    \end{scope}
  }
\end{tikzpicture}
\caption{Definition~\ref{def:coherent} at length three: $\Theta_k$ maps $\FL(3)$ onto $M_k$, drawn bold, and fixes every element of each earlier $M_j$, shaded.  Theorem~\ref{thm:D} gives a coherent chain of length two, and Theorem~\ref{thm:G} shows that none of length three exists.}
\label{fig:coherent}
\end{figure}
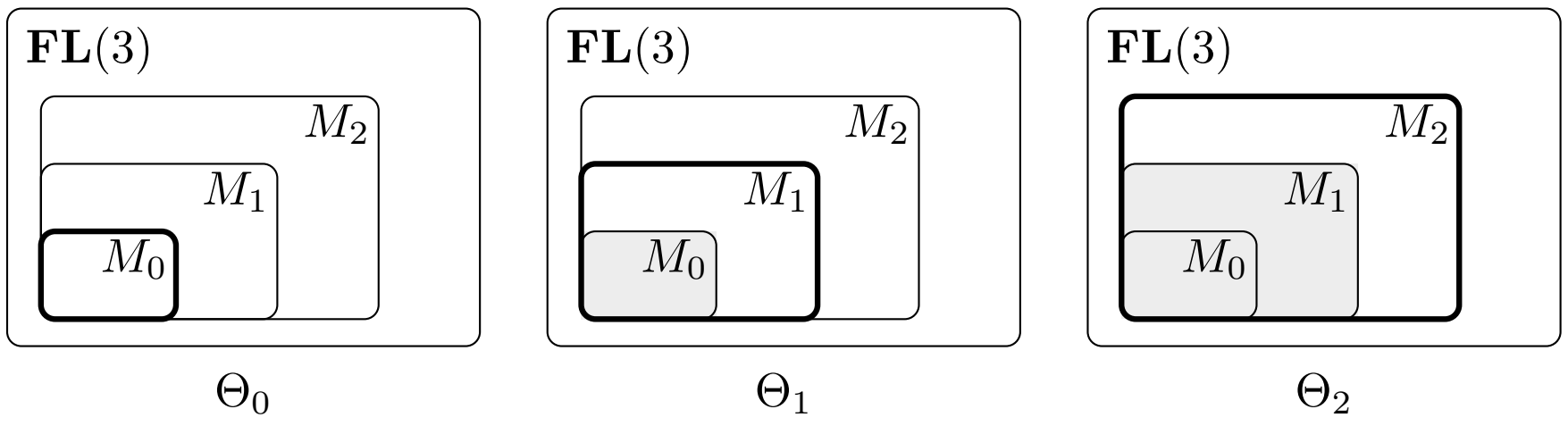  The next two statements extract longer ones from a chain.

\begin{expertread}
\begin{thmCfour}
\label{thm:C-four}
Four links $i < j < k < l$ of a chain with $\theta_{i,j} = \theta_{i,l}$ and $\theta_{j,k} = \theta_{j,l}$ give a coherent chain of length three, namely $\Theta_0 = \theta_{i,l}$, $\Theta_1 = \theta_{j,l}$, $\Theta_2 = \theta_{k,l}$.
\end{thmCfour}

\begin{proof}
Each $\theta_{\cdot,l}$ is a proper self-embedding by (c).  The cocycle identity turns the two hypotheses into the three coherence equations:
\begin{align*}
\Theta_1 \circ \Theta_0 &= \theta_{j,l} \circ \theta_{i,l} = \theta_{j,l} \circ \theta_{i,j} = \theta_{i,l} = \Theta_0,\\
\Theta_2 \circ \Theta_1 &= \theta_{k,l} \circ \theta_{j,l} = \theta_{k,l} \circ \theta_{j,k} = \theta_{j,l} = \Theta_1,
\end{align*}
and, since $\theta_{i,k} = \theta_{j,k} \circ \theta_{i,j} = \theta_{j,l} \circ \theta_{i,l} = \theta_{i,l}$ by the two hypotheses and the first equation,
\[
\Theta_2 \circ \Theta_0 = \theta_{k,l} \circ \theta_{i,l} = \theta_{k,l} \circ \theta_{i,k} = \theta_{i,l} = \Theta_0.  \qedhere
\]
\end{proof}
\end{expertread}

The third equation $\theta_{i,k} = \theta_{i,l}$, which one might expect to need as a hypothesis, follows from the other two.

\begin{furthercut}
\begin{thmCthree}
\label{thm:C-three}
Suppose some link $j$ of a chain has a transition from the bottom that recurs beyond every bound: for every $m$ there is $l > m$ with $\theta_{0,l} = \theta_{0,j}$.  Then the chain yields a coherent chain of length three.
\end{thmCthree}

\begin{proof}
First, $j \geq 1$: for $l > 0$ the transition $\theta_{0,l}$ is proper, so it is not $\theta_{0,0}$, the identity.  From the hypothesis read off recurrences $j < l_0 < l_1 < \cdots$ with $\theta_{0,l_t} = \theta_{0,j}$ for every $t$.  By Lemma~\ref{lem:bounded-transitions} the transitions $\theta_{j,l_t}$ take at most $N(j)$ values, so $\theta_{j,l_t} = \theta_{j,l_{t\pr}}$ for some $t < t\pr$.  The four links $0 < j < l_t < l_{t\pr}$ satisfy the hypotheses of the four-link form, with $i = 0$, $k = l_t$ and $l = l_{t\pr}$.
\end{proof}
\end{furthercut}

\subsection{Theorem~C: the coherent form}
\label{sec:theorem-c-coherent}

\begin{expertread}
\begin{thmCcoh}
\label{thm:C-coherent}
After passing to a subsequence of the chain, the transition $\theta_{j,k}$ depends only on $j$: there are proper self-embeddings $\Theta_0, \Theta_1, \dots$ with $\theta_{j,k} = \Theta_j$ for all $k > j$, and they form an infinite coherent chain.  Conversely, the images $M_j = \Theta_j(\FL(3))$ of any infinite coherent chain form a chain, inside $\FL(3)$.
\end{thmCcoh}

\begin{proof}
For a fixed $j$ the transitions $\theta_{j,k}$, $k > j$, take finitely many values (Lemma~\ref{lem:bounded-transitions}), so some value is taken on an infinite set of indices $k$.  Choose $k_0 = 0$ and an infinite set $S_0$ of indices $k > k_0$ on which $\theta_{k_0,k}$ is constant; then $k_1 = \min S_0$ and an infinite $S_1 \subseteq S_0$ of indices $k > k_1$ on which $\theta_{k_1,k}$ is constant; and so on.  Since $k_{j} \in S_{j-1} \subseteq \cdots \subseteq S_i$ for $i < j$, the value $\theta_{k_i,k_j}$ is the constant value of $\theta_{k_i,\cdot}$ on $S_i$, independent of $j > i$; call it $\Theta_i$.  The subsequence $L_{k_0} \subsetneq L_{k_1} \subsetneq \cdots$ is a chain whose transitions from link $i$ are all $\Theta_i$, a proper self-embedding by (c).  The cocycle identity for $i < j < l$ reads $\theta_{k_i,k_l} = \theta_{k_j,k_l} \circ \theta_{k_i,k_j}$, that is, $\Theta_i = \Theta_j \circ \Theta_i$: the sequence is coherent.

For the converse, Lemma~\ref{lem:coherent-images} gives $M_0 \subsetneq M_1 \subsetneq \cdots$ with every $M_j$ a copy of $\FL(3)$, a chain in $\FL(3)$.
\end{proof}

\begin{corollary}
\label{cor:dean-iff-coherent}
Dean's problem has a positive answer if and only if an infinite coherent chain exists.
\end{corollary}
\end{expertread}

\begin{furthercut}
\begin{remark}[which form is classical]
\label{rem:classical}
The one-step form is a finite pigeonhole, and the four-link form is finite algebra.  The length-three form needs a link $j$ with $N(j) + 1$ later recurrences of $\theta_{0,j}$, where $N(j)$ depends on the complexity of the generators of $L_j$, which nothing bounds in advance; classically some value of $\theta_{0,\cdot}$ recurs infinitely often and its first occurrence is such a $j$, but that is the infinite pigeonhole principle, not a finite one.  So the formal statement of the length-three form takes the recurrence as a hypothesis.  The coherent form repeats the same choice at every level; the formal development checks its diagonal argument with the infinite pigeonhole principle taken as a hypothesis (\agda{Theorem-C-coherent} in \agda{DeansProblem.\allowbreak Coherent.\allowbreak Infinite}).  Theorem~\ref{thm:H} needs neither: Remark~\ref{rem:constructive} proves it from the four-link form by finite searches, and Section~\ref{sec:which-step} says how that proof is checked.
\end{remark}

\begin{proposition}
\label{prop:unbounded}
No infinite coherent chain has generator images of bounded complexity.
\end{proposition}

\begin{proof}
If $\cx(\Theta_k(g)) \leq n$ for all $k$ and all generators $g$, then by Lemma~\ref{lem:finite} two links coincide, $\Theta_j = \Theta_k$ with $j < k$.  Then $\Theta_j = \Theta_k \circ \Theta_j = \Theta_j \circ \Theta_j$, so $\Theta_j(w) = \Theta_j(\Theta_j(w))$ and injectivity gives $\Theta_j(w) = w$ for every $w$: $\Theta_j$ is the identity, which is not proper.
\end{proof}

So the links of an infinite coherent chain, which is any chain read inside $\FL(3)$, have unbounded complexity, and no compactness argument under a uniform complexity bound is available for passing from coherent chains of every finite length to an infinite one.  Theorem~\ref{thm:G} makes the point moot.
\end{furthercut}

\subsection{The fixed-sublattice question}
\label{sec:question}

\begin{expertread}
By Definition~\ref{def:coherent} and Lemma~\ref{lem:coherent-images}, the link $\Theta_1$ of any coherent chain is a proper self-embedding of $\FL(3)$ whose fixed sublattice $\Fix(\Theta_1)$ contains a copy of $\FL(3)$, and by Theorem~\ref{thm:C} every chain yields such a link.  A negative answer to the following question would therefore have answered Dean's question in the negative at once.
\end{expertread}

\begin{question}
\label{q:fixed}
Does $\FL(3)$ have a proper self-embedding $\theta$ whose fixed sublattice $\Fix(\theta) = \{w : \theta(w) = w\}$ contains a sublattice isomorphic to $\FL(3)$?
\end{question}

The answer is positive, and the argument has to go one link further.

\begin{thmD}
\label{thm:D}
The self-embedding
\[
\theta_1 \colon x \mapsto x \vee (y \wedge z), \qquad y \mapsto y \vee (x \wedge z), \qquad z \mapsto z
\]
of $\FL(3)$ is proper, and it fixes the three elements
\[
m_1 = x \vee c, \qquad m_2 = z \vee ((x \vee z) \wedge q), \qquad m_3 = ((z \vee (p \wedge q)) \wedge (x \vee y)) \vee y,
\]
where $c = z \wedge (x \vee y)$, $p = x \vee c$ and $q = y \vee c$, and these three freely generate.  So $\Theta_0 \colon x, y, z \mapsto m_1, m_2, m_3$ and $\Theta_1 = \theta_1$ form a coherent chain of length two, and $\gen{m_1, m_2, m_3} \subsetneq \theta_1(\FL(3)) \subsetneq \FL(3)$.
\end{thmD}

\begin{proof}
Every claim is a Whitman decision in $\FL(3)$ on terms with at most thirty generator occurrences, and the module \agda{DeansProblem.\allowbreak Examples.\allowbreak FixedTriple} decides them all (nine decisions besides the properness of $\theta_1$).  We write out the checks, since they are short and the theorem bears on the shape of the whole argument.%

\begin{fullonly}
\emph{$\theta_1$ is a proper self-embedding.}  Write $r_1 = x \vee (y \wedge z)$, $r_2 = y \vee (x \wedge z)$, $r_3 = z$.  The three join conditions~\eqref{eq:join-condition} of Corollary~\ref{cor:1.13} hold because a generator is join prime and lies below no other generator and below no meet of two generators: $r_1 \leq r_2 \vee r_3$ would put $x$ below one of $y$, $x \wedge z$, $z$; $r_2 \leq r_1 \vee r_3$ would put $y$ below one of $x$, $y \wedge z$, $z$; and $r_3 \leq r_1 \vee r_2$ would put $z$ below one of $x$, $y \wedge z$, $y$, $x \wedge z$.  For the meet conditions~\eqref{eq:meet-condition}, $r_2 \wedge r_3 \leq r_1$, that is, $(y \vee (x \wedge z)) \wedge z \leq x \vee (y \wedge z)$, leaves by $\W$ four cases, each false: $y \vee (x \wedge z) \leq x \vee (y \wedge z)$ puts $y$ below $x$ or below $y \wedge z$; $z \leq x \vee (y \wedge z)$ puts $z$ below $x$ or $y$; $(y \vee (x \wedge z)) \wedge z \leq x$ puts $y \vee (x \wedge z)$ or $z$ below $x$; and $(y \vee (x \wedge z)) \wedge z \leq y \wedge z$ puts $y \vee (x \wedge z)$ or $z$ below $y$, the first needing $x \wedge z \leq y$.  The condition $r_1 \wedge r_3 \leq r_2$ is the same with $x$ and $y$ exchanged, and $r_1 \wedge r_2 \leq r_3$ would put $r_1$ or $r_2$ below $z$, hence $x$ or $y$ below $z$.  And $x$ is none of the three values: $x = x \vee (y \wedge z)$ would need $y \wedge z \leq x$, $x = y \vee (x \wedge z)$ would need $x \leq y$ or $x \leq z$, and $x \neq z$.
\end{fullonly}

\emph{The three elements are fixed.}  $\theta_1$ fixes $z$ and $x \vee y$, since $\theta_1(x \vee y) = x \vee (y \wedge z) \vee y \vee (x \wedge z) = x \vee y$; hence it fixes $c$, and $x \vee z$ too.  It fixes $p = x \vee c$ because $\theta_1(p) = x \vee (y \wedge z) \vee c = x \vee c$, as $y \wedge z \leq c$; likewise $q$.  So $m_1 = p$ and $m_2$, built from $z$, $x \vee z$ and $q$, are fixed.  For $m_3$, $\theta_1(m_3) = m_3 \vee (x \wedge z)$, and $x \wedge z \leq (z \vee (p \wedge q)) \wedge (x \vee y)$, so $\theta_1(m_3) = m_3$; equivalently, $m_3 = ((z \vee (p \wedge q)) \wedge (x \vee y)) \vee q$, a term in fixed elements, since $c \leq (z \vee (p \wedge q)) \wedge (x \vee y)$.

\begin{fullonly}
\emph{The three elements freely generate.}  For the join conditions~\eqref{eq:join-condition} it suffices to find a generator joinand of each $m_i$ that is below neither of the other two $m$'s nor any of their joinands: $x \leq m_2 \vee m_3$ fails because $x$ is join prime and below none of $z$, $(x \vee z) \wedge q$ (as $x \not\leq q$, which needs $x \leq y$ or $x \leq c$), $(z \vee (p \wedge q)) \wedge (x \vee y)$ (as $x \leq z \vee (p \wedge q)$ needs $x \leq q$) and $y$; $z \leq m_1 \vee m_3$ fails because $z$ is below none of $x$, $c$, $y$ and $(z \vee (p \wedge q)) \wedge (x \vee y)$ (as $z \not\leq x \vee y$); and $y \leq m_1 \vee m_2$ fails because $y$ is below none of $x$, $c$, $z$ and $(x \vee z) \wedge q$.  For the meet conditions~\eqref{eq:meet-condition}, each is an instance of $\W$ whose four cases fail: $m_2 \wedge m_3 \leq m_1 = x \vee c$ needs $m_2 \leq x \vee c$ (false, as $z \not\leq x \vee c$), $m_3 \leq x \vee c$ (false, as $y \not\leq x \vee c$), $m_2 \wedge m_3 \leq x$ (false, since $m_2 \not\leq x$ and $m_3 \not\leq x$) or $m_2 \wedge m_3 \leq c$, which needs $m_2 \wedge m_3 \leq z$, false since neither $m_2$ nor $m_3$ is below $z$.  $m_1 \wedge m_3 \leq m_2$ needs $m_1 \leq m_2$ (false, as $x \not\leq m_2$), $m_3 \leq m_2$ (false, as $y \not\leq m_2$), $m_1 \wedge m_3 \leq z$ (false) or $m_1 \wedge m_3 \leq (x \vee z) \wedge q$, which needs $m_1 \wedge m_3 \leq q = y \vee c$; that in turn needs $m_1 \leq q$ (false, $x \not\leq q$), $m_3 \leq q$ (false: its joinand $(z \vee (p \wedge q)) \wedge (x \vee y)$ is not below $y \vee c$, since neither meetand is, and it is below neither $y$ nor $c$, as it is not below $z$), $m_1 \wedge m_3 \leq y$ (false) or $m_1 \wedge m_3 \leq c$ (false, as it needs $m_1 \wedge m_3 \leq z$).  Finally $m_1 \wedge m_2 \leq m_3$ needs $m_1 \leq m_3$ (false, $x \not\leq m_3$), $m_2 \leq m_3$ (false, $z \not\leq m_3$), $m_1 \wedge m_2 \leq y$ (false) or $m_1 \wedge m_2 \leq (z \vee (p \wedge q)) \wedge (x \vee y)$, which needs $m_1 \wedge m_2 \leq z \vee (p \wedge q)$; that needs $m_1 \leq z \vee (p \wedge q)$ (false, $x \not\leq q$), $m_2 \leq z \vee (p \wedge q)$ (false: its joinand $(x \vee z) \wedge q$ is not below $z \vee (p \wedge q)$, because $x \vee z$ and $q$ are not, it is not below $z$, and it is not below $p \wedge q$ since it is not below $p = x \vee c$), $m_1 \wedge m_2 \leq z$ (false) or $m_1 \wedge m_2 \leq p \wedge q$ (false, as $m_1 \wedge m_2 \leq q$ needs $m_1 \leq q$, $m_2 \leq q$, $m_1 \wedge m_2 \leq y$ or $m_1 \wedge m_2 \leq c$, all false).  The antichain conditions follow from the join conditions~\eqref{eq:join-condition}.
\end{proof}
\end{fullonly}%

\begin{figure}[ht]
\centering
\begin{tikzpicture}[x=10mm, y=9mm]
  \node[anchor=north west] at (-3.9,6.15) {$[z)$};
  \node[anchor=south east] at (4.1,-0.45) {$(x \vee y]$};
  \node[fixed, upper, label=above:{$x \vee y \vee z$}] (top) at (0.6,5.45) {};
  \node[fixed, lower, label=left:{$x \vee y$}] (xy) at (-0.9,4.4) {};
  \node[fixed, triple, lower, label=right:{$m_3$}] (m3) at (0.3,3.4) {};
  \node[fixed, triple, upper, label=right:{$m_2$}] (m2) at (3,3.6) {};
  \node[fixed, triple, lower, label=left:{$p = m_1$}] (p) at (-2,2.6) {};
  \node[fixed, lower, label=left:{$q$}] (q) at (0.3,2.2) {};
  \node[moved, lower, label=left:{$\theta_1(x)$}] (t1x) at (-3,1.6) {};
  \node[moved, lower, label=right:{$\theta_1(y)$}] (t1y) at (1.3,1.95) {};
  \node[fixed, upper, label=right:{$z$}] (z) at (3,1.6) {};
  \node[moved, lower, label=left:{$x$}] (x) at (-3,0.5) {};
  \node[moved, lower, label=left:{$y$}] (y) at (1.45,1.05) {};
  \node[fixed, lower, label=below:{$c$}] (c) at (0,0) {};
  \begin{scope}[on background layer]
    \fill[region] (2.35,0.75) -- (4.1,0.75) -- (4.1,6.15) -- (-3.9,6.15)
      -- (-3.9,4.95) -- (2.35,4.95) -- cycle;
    \draw[le] (xy) -- (top);
    \draw[le] (m2) -- (top);
    \draw[le] (m3) -- (xy);
    \draw[le] (p) -- (xy);
    \draw[le] (q) -- (m3);
    \draw[le] (z) -- (m2);
    \draw[le] (t1x) -- (p);
    \draw[le] (t1y) -- (q);
    \draw[le] (x) -- (t1x);
    \draw[le] (y) -- (t1y);
    \draw[le] (c) -- (p);
    \draw[le] (c) -- (q);
    \draw[le] (c) -- (z);
  \end{scope}
  \draw[maps] (x) to[bend right=50] (t1x);
  \draw[maps] (y) to[bend right=50] (t1y);
\end{tikzpicture}
\caption{Twelve elements of $\FL(3)$, ordered as in $\FL(3)$: the shaded block is the filter $[z)$ of the partition lemma and the rest is the ideal $(x \vee y]$.  The self-embedding $\theta_1$ of Theorem~\ref{thm:D} fixes the filled elements and moves the open ones, $x$ and $y$ as the dashed arrows show; the ringed elements are its fixed free triple $m_1 = p$, $m_2$, $m_3$.  Every fixed element drawn lies above $c = z \wedge (x \vee y)$, as Lemma~\ref{lem:5} says it must, and no moved one does.}
\label{fig:theorem-d}
\end{figure}
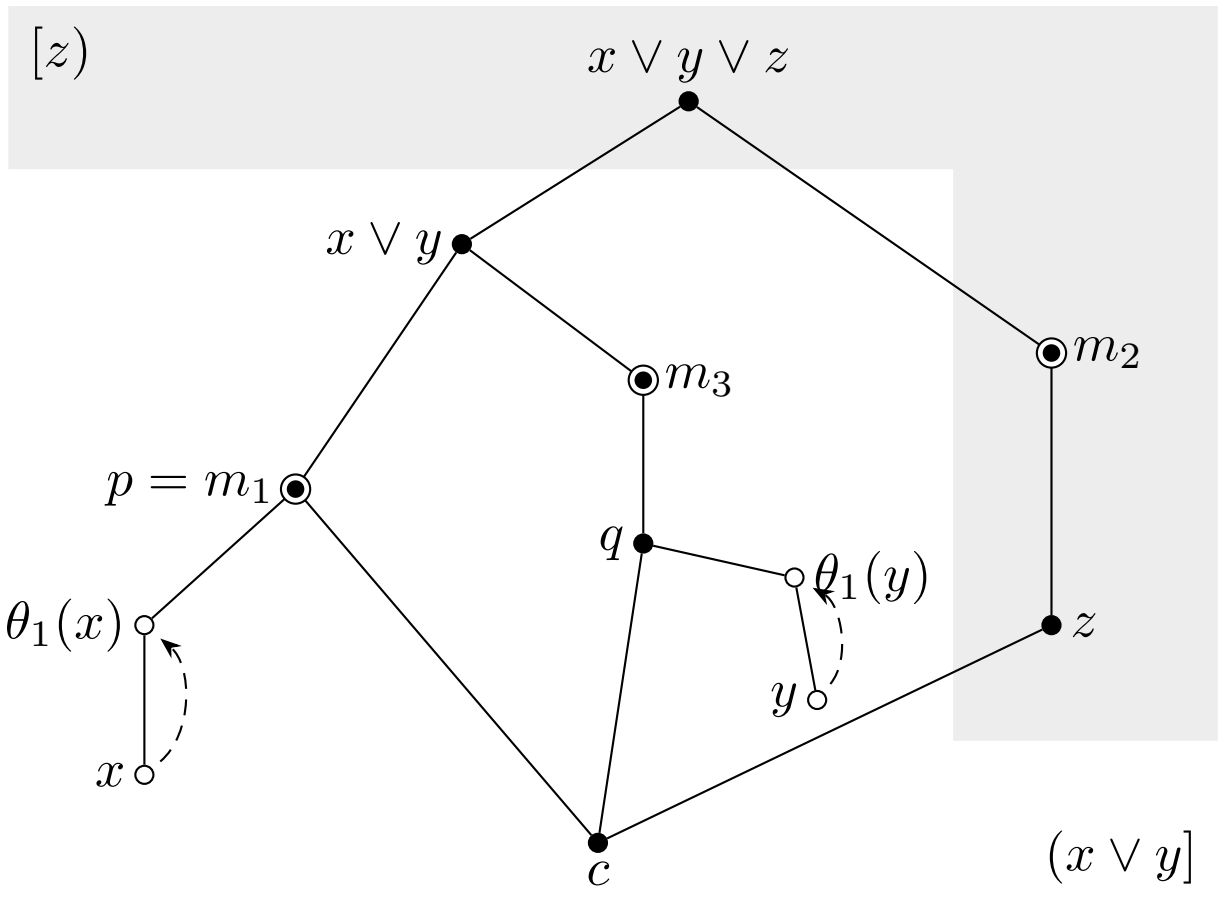

\begin{furthercut}
Where the triple comes from.  Since $\theta_1$ fixes $z$, it preserves the ideal $(x \vee y]$ and the filter $[z)$ of the partition lemma (Figure~\ref{fig:theorem-d}); inside the ideal it fixes $c$ and moves $x$ to $x \vee (y \wedge c)$ and $y$ to $y \vee (x \wedge c)$, so it fixes $p$ and $q$, and $\Fix(\theta_1)$ contains the sublattice $\gen{p, q, z}$, in which $p \wedge z = q \wedge z = c$.  That sublattice is infinite.  A search by computer over the $1291$ elements of $\FL(3)$ of complexity at most~$4$ (there are $3, 8, 14, 38$ and $1228$ of complexity $0$ to $4$) found $18$ fixed elements of $\theta_1$, all in $\gen{p, q, z}$ and among them $p$, $q$ and $z$ (of complexity $3$, $3$ and $0$), so the sublattice they generate is $\gen{p, q, z}$ itself.  The proof assistant certifies the triple; the enumeration and the counts are a computation, not a theorem.
\end{furthercut}

\begin{furthercut}
\begin{remark}
\label{rem:what-d-means}
Theorem~\ref{thm:D} closes the route to a negative answer through Question~\ref{q:fixed}.  It also shows that $\Fix(\theta_1)$ is large in the sense of embeddings: it contains a copy of $\FL(3)$, hence, by Theorem~1.28 of FJN, a copy of $\FL(\omega)$.  That does not decide the analogue of the theorem of Imrich and Turner \cite{ImrichTurner1989} (after Bestvina and Handel \cite{BestvinaHandel1992}) that the fixed subgroup of an injective endomorphism of a free group of rank $n$ has rank at most $n$: $\FL(3)$ itself contains $\FL(\omega)$, and whether $\Fix(\theta_1)$ is generated by three elements, as its known part $\gen{p, q, z}$ is, is open (Section~\ref{sec:questions}).\footnote{See \cite[Theorem~6.1]{Lupini2026} and the footnote{} on that preprint in Section~\ref{sec:problem}.}  A chain of three copies of $\FL(3)$ inside $\FL(3)$ is nothing new in itself (iterate any proper self-embedding); what Theorem~\ref{thm:D} gives is the coherence, the bottom copy fixed pointwise by the embedding that defines the middle one.
\end{remark}

The pair of Theorem~\ref{thm:D} extends to no third link above.

\begin{proposition}
\label{prop:rigid}
An endomorphism $\Phi$ of $\FL(3)$ that fixes $x \vee (y \wedge z)$, $y \vee (x \wedge z)$ and $z$ is the identity.  Hence no proper self-embedding $\Theta_2$ fixes $\theta_1(\FL(3))$ pointwise, and no coherent chain of length three has $\Theta_0 = (m_1, m_2, m_3)$ and $\Theta_1 = \theta_1$.
\end{proposition}

\begin{proof}
Write $a = \Phi(x)$, $b = \Phi(y)$; then $\Phi(z) = z$ and $a \vee (b \wedge z) = x \vee (y \wedge z)$.  Since $x$ is join prime, $x \leq a$ or $x \leq b \wedge z$, and the latter puts $x$ below $z$; so $x \leq a$, and likewise $y \leq b$.  Now $b \wedge z \leq x \vee (y \wedge z)$, and $\W$ leaves four cases: $b \leq x \vee (y \wedge z)$ puts $y \leq b$ below $x$ or below $y \wedge z$, hence below $z$; $z \leq x \vee (y \wedge z)$ puts $z$ below $x$ or below $y$; $b \wedge z \leq x$ puts $b$ or $z$ below $x$, and $b \leq x$ would put $y$ below $x$; and $b \wedge z \leq y \wedge z$ gives $b \wedge z \leq y$, so $b \leq y$ or $z \leq y$.  Only the last survives, so $b \leq y$ and $b = y$; by the symmetry exchanging $x$ and $y$, $a = x$.  The two consequences follow, since a proper self-embedding is not the identity.
\end{proof}

The proof uses Whitman's rules only, and it needs no injectivity of $\Phi$.  Section~\ref{sec:rigidity} proves that nothing proper fixes the image of \emph{any} self-embedding with a fixed free triple.
\end{furthercut}

\subsection{Three facts about fixed sublattices}
\label{sec:facts}

\begin{lemma}
\label{lem:two-generators}
Let $X$ be finite with $|X| \geq 2$, and let $\theta$ be an injective endomorphism of $\FL(X)$ that fixes every generator but $z$.  Then $\theta$ fixes $z$ too.  In particular a self-embedding of $\FL(3)$ that fixes two generators is the identity, $\Fix(\theta)$ contains at most one generator when $\theta$ is proper, and a triple $(x, y, r)$ is free only for $r = z$.
\end{lemma}

\begin{proof}
Let $Y = X \setminus \{z\}$.  By the partition lemma, $\theta(z) \leq \bigvee Y$ or $z \leq \theta(z)$.  The first would read $\theta(z) \leq \bigvee_{y \in Y} \theta(y)$, and the join condition~\eqref{eq:join-condition} of Corollary~\ref{cor:1.13} for the free family $(\theta(g))_{g \in X}$ would put $z \in Y$; so $z \leq \theta(z)$.  Dually $\theta(z) \leq z$.  For the last claim, if $(x, y, r)$ is free then $x, y, z \mapsto x, y, r$ is a self-embedding fixing $x$ and $y$, hence the identity, so $r = z$.
\end{proof}

\begin{lemma}
\label{lem:descent}
Let $\theta$ be a self-embedding of $\FL(n)$ and $w$ a fixed point of $\theta$.
\begin{enumerate}[label=(\alph*)]
\item If $u$ is a canonical joinand or meetand of $w$ that is not a generator, then some positive power of $\theta$ fixes $u$.
\item If $\cx(w) \geq 1$, there are $d \geq 1$ and $v$ with $\cx(v) = 1$ and $\theta^d(v) = v$.  More precisely, if $\cx(w) \geq 2$, there are an element $u$ of complexity $2$ fixed by a positive power $\theta^e$ of $\theta$, and a canonical joinand or meetand $v$ of $u$ of complexity $1$ fixed by a positive power of $\theta^e$.
\end{enumerate}
\end{lemma}

\begin{proof}
(a) Let $u$ be a canonical joinand of $w$, not a generator (the meet case is dual).  By Lemma~\ref{lem:9.13}, $\theta(u)$ is a canonical joinand of $\theta(w) = w$, and $\theta(u)$ is not a generator either: a generator in the image of $\theta$ is the image of a generator (Corollary~\ref{cor:1.6}), and $u$ is not one.  So the orbit $u, \theta(u), \theta^2(u), \dots$ stays in the finite set of canonical joinands of $w$ (the same set for $\theta(w)$ as for $w$, by the uniqueness of the canonical form), it repeats, $\theta^a(u) = \theta^b(u)$ with $a < b$, and injectivity of $\theta^a$ gives $\theta^{b-a}(u) = u$.

(b) Induction on $\cx(w)$.  If $\cx(w) = 1$ take $v = w$ and $d = 1$.  If $\cx(w) = m + 1 \geq 2$, let $u$ be a canonical joinand or meetand of $w$ with $\cx(u) = m \geq 1$, the child of largest complexity; it is not a generator, so by (a) a positive power $\theta^{d_1}$ fixes it, and the induction hypothesis for the self-embedding $\theta^{d_1}$ and its fixed point $u$ gives $v$ of complexity $1$ fixed by a positive power of $\theta^{d_1}$, hence of $\theta$.  Keeping the last parent gives the more precise form: when $\cx(w) = 2$, take $u = w$, $e = 1$, and $v$ its child of complexity $1$, fixed by a positive power of $\theta$ by (a); when $\cx(w) > 2$, descend one level and apply the induction hypothesis.
\end{proof}

In the formal development the lemma is \agda{Lemma-6·11} in \agda{DeansProblem.\allowbreak Fixed}, in $\FL(n)$ as stated: part (a) is \agda{joinand-periodic} and \agda{meetand\allowbreak -periodic}, and part (b) is \agda{descent}, with \agda{descent-parent} for the element $u$ of its more precise form.  So a fixed point of positive complexity descends to a join or meet of generators fixed by a power of $\theta$; $\theta_1$ fixes $x \vee y$, $x \vee z$ and $y \vee z$ outright.

\begin{fullonly}
\begin{computation}
\label{comp:fixed}
The following were computed in Python by a memoized implementation of Whitman's procedure, over the elements of $\FL(3)$ of complexity at most~$4$; the positive facts they rest on (a triple is fixed and free) are the ones decided in Agda, and no proof in this paper cites a count.
\begin{itemize}
\item The self-embedding $\rho \colon x \mapsto x \vee (y \wedge z)$ of Section~\ref{sec:inside}, cyclically, fixes $11$ elements of complexity at most~$4$: the eight-element Boolean lattice generated by $x \vee y$, $y \vee z$, $z \vee x$, with bottom $\bot^{B} = (x \vee y) \wedge (y \vee z) \wedge (z \vee x)$, and the three elements $x \vee \bot^{B}$, $y \vee \bot^{B}$, $z \vee \bot^{B}$; they form a sublattice with no free triple.  Corollary~\ref{cor:freese-eleven} proves that these are all the fixed points of $\rho$.
\item $\theta_1$ fixes $18$ elements of complexity at most~$4$, all in $\gen{p, q, z}$.
\item Of the $48$ self-embeddings of $\FL(3)$ whose generator images have complexity at most~$2$, exactly six have a free triple in their fixed sublattice (within the bounds of the search): those that fix one generator and move the other two by meets with it, or dually by joins, that is, $\theta_1$ up to permutation of the generators and duality (Section~\ref{sec:inflationary}).  $\rho$ and its dual fix the $11$ elements above; the other $40$ fix between $2$ and $8$ elements of complexity at most~$4$ and no free triple.
\item No pair $(a, b)$ of elements of complexity at most~$4$ other than $(x, y)$ satisfies $a \vee (b \wedge z) = \theta_1(x)$ and $b \vee (a \wedge z) = \theta_1(y)$ with $(a, b, z)$ free; Proposition~\ref{prop:rigid} proves this for all $a$, $b$.
\end{itemize}
\end{computation}
\end{fullonly}

\section{Self-embeddings with a fixed free triple}
\label{sec:structure}

Let $\Theta$ be a proper self-embedding of $\FL(3)$ whose fixed sublattice contains a free triple, the middle link of any coherent chain of length three.  This section proves that $\Theta$ fixes a generator (Theorem~\ref{thm:E}), and describes its shape (Corollary~\ref{cor:E1}).

\subsection{Inflationary self-embeddings and duality}
\label{sec:inflationary}

\begin{definition}
\label{def:inflationary}
A self-embedding $\Theta$ of $\FL(3)$ is \emph{inflationary} when $g \leq \Theta(g)$ for each generator $g$, and \emph{strictly inflationary} when moreover $\Theta(g) \neq g$ for each $g$; \emph{deflationary} is dual.
\end{definition}

If $\Theta$ is inflationary then $w \leq \Theta(w)$ for every $w$,
as follows from a simple induction.

The \emph{dual} of an element $w$ is the element $\delta(w)$ named by the term obtained from a term naming $w$ by exchanging $\vee$ and $\wedge$; it is well defined because Whitman's rules are self-dual, and $\delta$ is an order-reversing involution of $\FL(3)$ fixing the generators.  The \emph{dual} of a self-embedding $\Theta$ is $\Theta^{\delta} = \delta \circ \Theta \circ \delta$, the self-embedding with $\Theta^{\delta}(g) = \delta(\Theta(g))$.  It is injective; it misses the same generators as $\Theta$, so it is proper when $\Theta$ is; it fixes $\delta(w)$ exactly when $\Theta$ fixes $w$; its powers are the duals of the powers of $\Theta$; a triple is free exactly when its \emph{dual triple}, the triple of duals, is, since the join condition~\eqref{eq:join-condition} of Corollary~\ref{cor:1.13} for the duals is the meet condition~\eqref{eq:meet-condition} for the originals; and $\Theta^{\delta}$ is inflationary exactly when $\Theta$ is deflationary.  (In the formal development: \agda{dualˢ}, \agda{dual-free}, \agda{dual-misses}, \agda{dual-fixed}, \agda{dual-iter}.)  A statement about $\Theta$ holds \emph{up to duality} when it holds for $\Theta$ or for $\Theta^{\delta}$.

\subsection{Four lemmas}
\label{sec:four-lemmas}

\begin{lemone}[a fixed join of two generators]
\label{lem:1}
Let $\Theta$ be a self-embedding of $\FL(3)$ with $\Theta(x \vee y) = x \vee y$.  Then $\Theta(x), \Theta(y) \leq x \vee y$ and $z \leq \Theta(z)$, and either $x \leq \Theta(x)$ and $y \leq \Theta(y)$, or $x \leq \Theta(y)$ and $y \leq \Theta(x)$; in the second case $\Theta^2$ is inflationary.
\end{lemone}

\begin{proof}
$\Theta(x) \vee \Theta(y) = x \vee y$ gives the first claim.  Since $x \leq \Theta(x) \vee \Theta(y)$ and $x$ is join prime, $x \leq \Theta(x)$ or $x \leq \Theta(y)$, and likewise for $y$.  If both $x$ and $y$ lay below $\Theta(x)$, then $\Theta(y) \leq x \vee y \leq \Theta(x)$, and an embedding reflects the order, so $y \leq x$, which is false; likewise both cannot lie below $\Theta(y)$.  This leaves the two cases.  In the second, $x \leq \Theta(y) \leq \Theta(\Theta(x)) = \Theta^2(x)$ and likewise $y \leq \Theta^2(y)$, and $z \leq \Theta(z) \leq \Theta^2(z)$ once the claim about $z$ is proved.  For $z$: by the partition lemma, $\Theta(z) \leq x \vee y$ or $z \leq \Theta(z)$.  The first would give $\Theta(x \vee y \vee z) = \Theta(x \vee y) \vee \Theta(z) \leq x \vee y = \Theta(x \vee y)$, hence $x \vee y \vee z \leq x \vee y$, implying $z \leq x \vee y$, which is false.
\end{proof}

\begin{lemtwo}[a generator fixed by a power]
\label{lem:2}
Let $\Theta$ be a proper self-embedding of $\FL(3)$ and $n \geq 1$.  If $\Theta^n$ fixes a generator $g$, so does $\Theta$.
\end{lemtwo}

\begin{proof}
Put $g\pr = \Theta^{n-1}(g)$, so that $\Theta(g\pr) = g$.  The generator $g$ lies in the image of $\Theta$, so by Corollary~\ref{cor:1.6} it is the image $\Theta(h)$ of a generator $h$, and injectivity gives $g\pr = h$: $g\pr$ is a generator.  If $g\pr = g$ we are done.  Otherwise $\Theta^n$ fixes $g\pr$ as well, since $\Theta^n(g\pr) = \Theta^{n-1}(\Theta(g\pr)) = \Theta^{n-1}(g) = g\pr$; so $\Theta^n$ fixes two generators and is the identity by Lemma~\ref{lem:two-generators}.  Then every $w$ equals $\Theta(\Theta^{n-1}(w))$, so $\Theta$ is onto, which it is not.
\end{proof}

\begin{expertread}
\begin{lemthree}[fixed points of an inflationary self-embedding]
\label{lem:3}
Let $\Theta$ be an inflationary self-embedding of $\FL(3)$.
\begin{enumerate}[label=(\alph*)]
\item If $w$ is a fixed point of $\Theta$ and $u$ is a canonical joinand or meetand of $w$ that is not a generator, then $\Theta(u) = u$.  Hence $\Theta$ fixes every non-generator node of the canonical tree of each of its fixed points.
\item Suppose $w$ is a fixed point of $\Theta$ and that a generator $g$ is a canonical meetand of $w$.  Then $\Theta(g) = g$.
\begin{cutcandidate}
\item If $\Theta$ is strictly inflationary, no fixed point of $\Theta$ lies below a generator.
\end{enumerate}
\end{cutcandidate}
\end{lemthree}

\begin{proof}
(a) By Lemma~\ref{lem:9.13}, $\Theta(u)$ is a canonical joinand of $\Theta(w) = w$, and $u \leq \Theta(u)$ since $\Theta$ is inflationary.  The canonical joinands of $w$ form an antichain, so two of them that are comparable are equal: $\Theta(u) = u$.  Meetands are dual, and the second sentence follows by induction down the tree.

(b) Let $w = g \wedge u_2 \wedge \cdots \wedge u_n$ canonically.  Then $\Theta(w) = \Theta(g) \wedge \Theta(u_2) \wedge \cdots \wedge \Theta(u_n) = w \leq g$, so by rule (5) of Theorem~\ref{thm:whitman} some meetand's image lies below $g$.  If $\Theta(u_i) \leq g$ for some $i \geq 2$, then $u_i \leq \Theta(u_i) \leq g$, against the antichain condition on the canonical meetands of $w$ (if $u_i$ is itself a generator, $u_i \leq g$ gives $u_i = g$, the same contradiction).  So $\Theta(g) \leq g$, and with $g \leq \Theta(g)$ we get $\Theta(g) = g$.

\begin{cutcandidate}
(c) Suppose $w = \Theta(w) \leq g$.  Assume $w$ is of minimal complexity among the fixed points of $\Theta$ that lie below $g$.
If $w$ is a generator, it is $g$, since distinct generators are incomparable.  But then $\Theta(g) = g$, violating strictness.
If $w$ is a proper meet, we have $g \geq w = w_1 \wedge w_2 \wedge \cdots \wedge w_n$ canonically, and so, by rule (5) of Theorem~\ref{thm:whitman}, $w_i \leq g$ for some $i$.  Using parts (a) and (b), we conclude that $w_i$ is a fixed point of $\Theta$, violating minimality.
For the remaining case, when $w$ is a proper join, none of the joinands can be a generator.  Indeed, if a joinand $u$ is a generator, then $u < w \leq g$, an impossibility.  So (a) implies that the canonical joinands are fixed points of $\Theta$, below $w \leq g$, violating minimality.
\end{proof}
\end{cutcandidate}
\end{expertread}

\begin{cutcandidate}
So a proper inflationary self-embedding of $\FL(3)$ fixes no meet of two generators, and if it fixes $z$ then no fixed point of it has a meet node with $x$ or $y$ among its meetands: the canonical trees of its fixed points hang $x$ and $y$ on joins only.
\end{cutcandidate}

\begin{fullonly}
\begin{expertread}
\begin{lemfour}[inflation transfers through a fixed element]
\label{lem:4}
Let $\Phi$ be a self-embedding of $\FL(3)$ with $\Phi(z) = z$ that fixes an element $a$ with $x \leq a$ and $y \not\leq a$.  Then $x \leq \Phi(x)$.
\end{lemfour}

\begin{proof}
Let $h \colon \FL(3) \to \mathbf{2}$ be the homomorphism onto the two-element lattice with $h(x) = 1$ and $h(y) = h(z) = 0$.  By the partition lemma at $x$, every $w$ satisfies exactly one of $x \leq w$ and $w \leq y \vee z$; the first gives $h(w) = 1$ and the second $h(w) = 0$, so $h(w) = 1$ exactly when $x \leq w$.  For any term $t$ naming $w$, $h(\Phi(w))$ is the value of $t$ in $\mathbf{2}$ under the valuation $v \mapsto h(\Phi(v))$, since $h \circ \Phi$ is a homomorphism.  Now $h(a) = 1$ and $a = \Phi(a)$, so $t_a$, a term naming $a$, takes the value $1$ under $v \mapsto h(\Phi(v))$, a valuation with $h(\Phi(z)) = h(z) = 0$.  If also $h(\Phi(x)) = 0$, the valuation is pointwise below the valuation $v \mapsto [v = y]$, and lattice terms are monotone in their valuations, so $t_a$ takes the value $1$ under $v \mapsto [v = y]$ as well; but that value is $h\pr(a)$ for the homomorphism $h\pr$ with $h\pr(y) = 1$, $h\pr(x) = h\pr(z) = 0$, which is $1$ exactly when $y \leq a$, by the partition lemma at $y$.  Since $y \not\leq a$, we conclude $h(\Phi(x)) = 1$, that is, $x \leq \Phi(x)$.
\end{proof}
\end{expertread}
\end{fullonly}

\subsection{The fourteen elements}
\label{sec:fourteen}

\begin{figure}[ht]
\centering
\begin{tikzpicture}[x=12mm, y=9mm]
  \node[fixed, label=above:{$x \vee y \vee z$}] (top) at (0,5) {};
  \node[fixed, label=left:{$x \vee y$}] (xy) at (-2,4) {};
  \node[fixed, label={[yshift=4pt]right:{$x \vee z$}}] (xz) at (0,4) {};
  \node[fixed, label=right:{$y \vee z$}] (yz) at (2,4) {};
  \node[fixed, label=left:{$A$}] (A) at (-2,3) {};
  \node[fixed, label={[yshift=-5pt]right:{$B$}}] (B) at (0,3) {};
  \node[fixed, label=right:{$C$}] (C) at (2,3) {};
  \node[fixed, label=left:{$x \vee \bot^{B}$}] (xb) at (-2,2) {};
  \node[fixed, label=right:{$y \vee \bot^{B}$}] (yb) at (0,2) {};
  \node[fixed, label=right:{$z \vee \bot^{B}$}] (zb) at (2,2) {};
  \node[fixed, label=below:{$\bot^{B}$}] (bot) at (0,0.4) {};
  \node[moved, label=left:{$x$}] (x) at (-2.35,1.35) {};
  \node[moved, label=left:{$y$}] (y) at (-0.35,1.35) {};
  \node[moved, label=right:{$z$}] (z) at (2.35,1.35) {};
  \begin{scope}[on background layer]
    \draw[le] (xy) -- (top);
    \draw[le] (xz) -- (top);
    \draw[le] (yz) -- (top);
    \draw[le] (A) -- (xy);
    \draw[le] (A) -- (xz);
    \draw[le] (B) -- (xy);
    \draw[le] (B) -- (yz);
    \draw[le] (C) -- (xz);
    \draw[le] (C) -- (yz);
    \draw[le] (xb) -- (A);
    \draw[le] (yb) -- (B);
    \draw[le] (zb) -- (C);
    \draw[le] (bot) -- (xb);
    \draw[le] (bot) -- (yb);
    \draw[le] (bot) -- (zb);
    \draw[le] (x) -- (xb);
    \draw[le] (y) -- (yb);
    \draw[le] (z) -- (zb);
  \end{scope}
\end{tikzpicture}
\caption{The fourteen elements of $T_{14}$, ordered as in $\FL(3)$, each line joining an element to one that covers it among the fourteen.  The eleven filled are $T_{11}$, the fixed points of the self-embedding $\rho$ of Section~\ref{sec:inside} (Corollary~\ref{cor:freese-eleven}), which moves the three generators, drawn open; by Proposition~\ref{prop:eleven} a strictly inflationary self-embedding fixes nothing outside $T_{11}$.}
\label{fig:fourteen}
\end{figure}
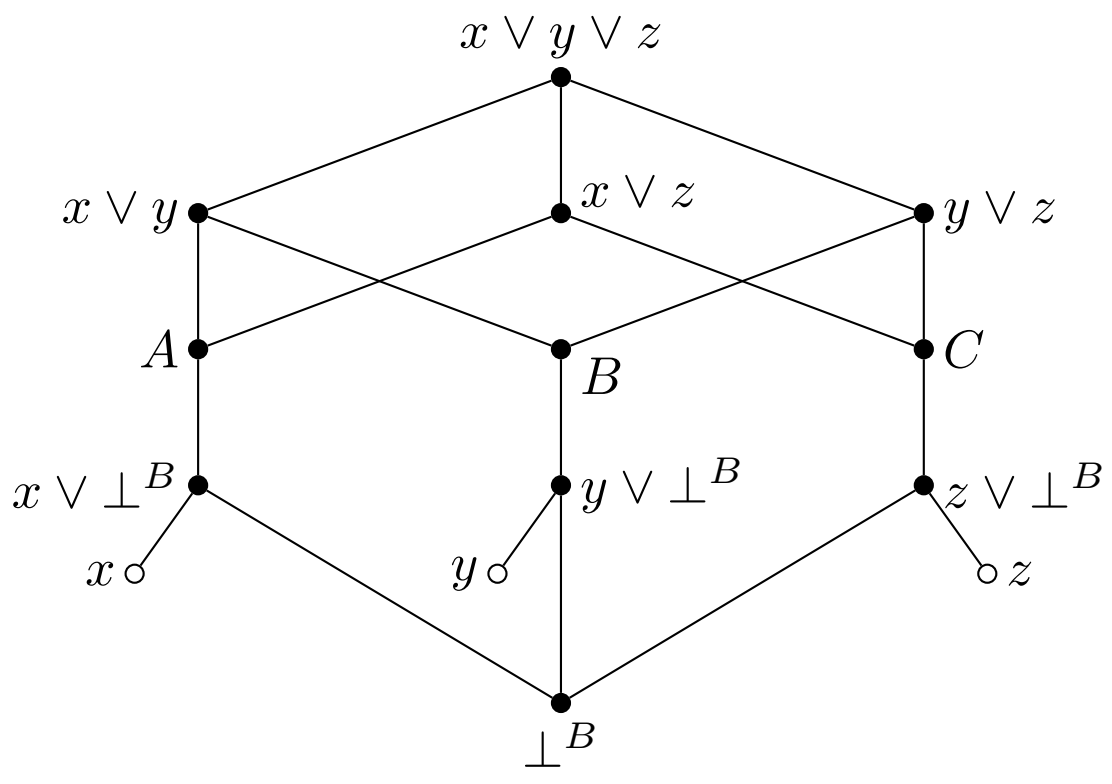

Let $S$ be the set of elements of $\FL(3)$ whose canonical trees have no generator among the meetands of any meet node.  Put
\begin{gather*}
A = (x \vee y) \wedge (x \vee z), \qquad B = (x \vee y) \wedge (y \vee z), \qquad C = (x \vee z) \wedge (y \vee z),\\
\bot^{B} = (x \vee y) \wedge (x \vee z) \wedge (y \vee z),
\end{gather*}
where the superscript of $\bot^{B}$ is for the Boolean lattice of Computation~\ref{comp:fixed}, not for the element $B$, and let $T_{14}$ be the set of the fourteen elements
\[
x, \ y, \ z, \quad x \vee y, \ x \vee z, \ y \vee z, \ x \vee y \vee z, \quad A, \ B, \ C, \ \bot^{B}, \quad x \vee \bot^{B}, \ y \vee \bot^{B}, \ z \vee \bot^{B},
\]
and $T_{11}$ the eleven non-generators among them (Figure~\ref{fig:fourteen}).

\begin{proposition}
\label{prop:fourteen}
$S = T_{14}$.  Every element of $T_{11}$ has complexity at most~$3$.
\end{proposition}

\begin{proof}

The fourteen terms displayed are in canonical form, by the conditions of Theorem~\ref{thm:canonical}, and none has a generator meetand; so $T_{14} \subseteq S$, and the eleven non-generators have complexities $1, 1, 1, 1, 2, 2, 2, 2, 3, 3, 3$.  

For the converse, we first note that $T_{14}$ is closed under joins and $T_{11}$ under meets; in fact $T_{11}$ is an interval sublattice of $\FL(3)$, dual to the connected component of $0$ (FJN, Example~3.44 and Theorem~7.5(1)).  Both closure facts are straightforward free lattice calculations, finitely many Whitman decisions, which the proof assistant runs (\agda{T₁₄-∨} and \agda{T₁₁-∧} in \agda{DeansProblem.\allowbreak FixedGenerator}).

Now let $w \in S$ and argue by induction on the complexity, $\cx(w)$, of $w$.  A generator is in $T_{14}$.  If $w$ is a canonical join, each joinand is a generator or a canonical meet; a meet joinand is a node of the tree of $w$, so it is in $S$, and it has smaller complexity, so it is in $T_{14}$ by induction; hence $w$, a join of elements of $T_{14}$, is in $T_{14}$.  If $w$ is a canonical meet, no meetand is a generator, since $w \in S$; so each meetand is a canonical join in $S$ of smaller complexity, in $T_{14}$ by induction and not a generator, hence in $T_{11}$; and $w$, a meet of elements of $T_{11}$, is in $T_{11}$.
\end{proof}

\begin{proposition}
\label{prop:eleven}
A strictly inflationary self-embedding $\Psi$ of $\FL(3)$ fixes only elements of $T_{11}$; in particular it has at most eleven fixed points, each of complexity at most~$3$.
\end{proposition}

\begin{proof}
Let $w$ be a fixed point.  By Lemma~\ref{lem:3}(a), $\Psi$ fixes every non-generator node of the canonical tree of $w$, and by Lemma~\ref{lem:3}(b) a meet node with a generator meetand $g$ would make $\Psi(g) = g$, which strictness forbids.  So no meet node of the tree of $w$ has a generator meetand: $w \in S = T_{14}$.  And $w$ is not a generator, again by strictness, so $w \in T_{11}$.  The complexity bound is Proposition~\ref{prop:fourteen}.
\end{proof}

\subsection{Theorem~E}
\label{sec:theorem-e}

\begin{expertread}
\begin{proposition}
\label{prop:inflationary-power}
Let $\Theta$ be a proper self-embedding of $\FL(3)$ fixing a free triple pointwise.  Then a positive power of $\Theta$, or of its dual $\Theta^{\delta}$, is inflationary and fixes the triple, or the dual triple, pointwise.
\end{proposition}

\begin{proof}
Let $\Theta$ fix the free triple $(t_1, t_2, t_3)$, and let $\iota \colon \FL(3) \to \FL(3)$ be the embedding sending $x, y, z$ to $t_1, t_2, t_3$.  Then $\Theta$ fixes every element of the image of $\iota$, since it fixes the generators of that sublattice.  Let $w_4 = x \vee (y \wedge (z \vee (x \wedge y)))$, an element of complexity $4$ (the term is canonical), and $u_4 = \iota(w_4)$.  Then $\Theta(u_4) = u_4$ and $\cx(u_4) \geq 4$ by Theorem~\ref{thm:B}.  By Lemma~\ref{lem:descent}(b) there are an element $u$ of complexity $2$ fixed by a positive power $\Theta^e$, and a canonical joinand or meetand $v$ of $u$ of complexity $1$ fixed by a positive power $\Psi$ of $\Theta^e$, hence of $\Theta$.  An element of complexity $1$ is a join or a meet of two or three distinct generators.  It is not $x \vee y \vee z$: that element is above every element of $\FL(3)$, so if it were a canonical joinand of $u$ the other joinands of $u$ would lie below it, against the antichain condition, and if it were a canonical meetand of $u$ the other meetands would lie below it, against the antichain condition again; either way $u$ would not be a canonical join or meet of two or more elements with it among them.  Dually $v$ is not $x \wedge y \wedge z$.  So, after renaming the generators, $v = x \vee y$ or $v = x \wedge y$.

If $v = x \vee y$, Lemma~\ref{lem:1} applied to $\Psi$ says that $\Psi$ or $\Psi^2$ is inflationary; it is a positive power of $\Theta$, and every power of $\Theta$ fixes $t_1, t_2, t_3$.  If $v = x \wedge y$, then $\Theta^{\delta}$ is a proper self-embedding fixing the free triple $(\delta(t_1), \delta(t_2), \delta(t_3))$, and the dual of $\Psi$, a positive power of $\Theta^{\delta}$, fixes $\delta(v) = x \vee y$; the join case for $\Theta^{\delta}$ gives an inflationary positive power of $\Theta^{\delta}$ fixing the dual triple.
\end{proof}

\begin{thmE}
\label{thm:E}
Every proper self-embedding of $\FL(3)$ whose fixed sublattice contains a free triple fixes a generator.
\end{thmE}

\begin{proof}
Let $\Theta$ be such a self-embedding, with the fixed free triple $(t_1, t_2, t_3)$, and let $u_4$ be the fixed element of complexity at least $4$ of the proof of Proposition~\ref{prop:inflationary-power}.  By that proposition, a positive power $\Phi$ of $\Theta$ or of $\Theta^{\delta}$ is inflationary.  Suppose first that $\Phi = \Theta^m$ with $m \geq 1$.  $\Phi$ is proper, since a power of a proper self-embedding is proper (its image lies inside the image of $\Theta$), and it fixes $u_4$.  If $\Phi$ moved every generator it would be strictly inflationary, and by Proposition~\ref{prop:eleven} all its fixed points would have complexity at most $3$; but $\cx(u_4) \geq 4$.  So $\Phi$ fixes a generator $g$, and by Lemma~\ref{lem:2} so does $\Theta$.  If instead $\Phi = (\Theta^{\delta})^m$, the same argument applied to $\Theta^{\delta}$, a proper self-embedding fixing the free triple $(\delta(t_1), \delta(t_2), \delta(t_3))$ and hence, as in the proof of Proposition~\ref{prop:inflationary-power}, a fixed point of complexity at least $4$, gives a generator $g$ fixed by $\Theta^{\delta}$; then $\Theta$ fixes $\delta(g) = g$.
\end{proof}
\end{expertread}

In the formal development (\agda{fixes-generator} in \agda{DeansProblem.\allowbreak FixedGenerator}) the choice between ``fixes a generator'' and ``moves every generator'' for $\Phi$ is made by three Whitman decisions, and the second branch is refuted by Proposition~\ref{prop:eleven} exactly as above.

\begin{furthercut}
\begin{expertread}
\begin{corollary}
\label{cor:freese-eleven}
The self-embedding $\rho \colon x \mapsto x \vee (y \wedge z)$, $y \mapsto y \vee (z \wedge x)$, $z \mapsto z \vee (x \wedge y)$ of Section~\ref{sec:inside} fixes exactly the eleven elements of $T_{11}$.
\end{corollary}

\begin{proof}
$\rho$ is strictly inflationary: $x \leq x \vee (y \wedge z)$, and $x \vee (y \wedge z) \neq x$ since $y \wedge z \not\leq x$; cyclically.  So by Proposition~\ref{prop:eleven} its fixed points lie in $T_{11}$.  Conversely $\rho$ fixes $x \vee y$, since $\rho(x \vee y) = x \vee (y \wedge z) \vee y \vee (z \wedge x) = x \vee y$, and likewise $x \vee z$ and $y \vee z$; hence it fixes the sublattice they generate, which contains $x \vee y \vee z$, $A$, $B$, $C$ and $\bot^{B}$; and it fixes $x \vee \bot^{B}$, since $\rho(x \vee \bot^{B}) = x \vee (y \wedge z) \vee \bot^{B} = x \vee \bot^{B}$ as $y \wedge z \leq \bot^{B}$, and likewise $y \vee \bot^{B}$ and $z \vee \bot^{B}$.
\end{proof}
\end{expertread}

So the fixed sublattice of $\rho$, the one-step example read inside $\FL(3)$, is the eleven-element lattice of Computation~\ref{comp:fixed}, and it contains no free triple; Theorem~\ref{thm:E} rules out a fixed free triple for every proper strictly inflationary self-embedding.
\end{furthercut}

\begin{furthercut}
\begin{corEone}[the shape]
\label{cor:E1}
Let $\Theta$ be a proper self-embedding of $\FL(3)$ with a fixed free triple.  Then $\Theta$ fixes exactly one generator, say $z$; $\Theta(x \vee y) \leq x \vee y$ and $\Theta(x \wedge y) \geq x \wedge y$; and, up to duality, $\Theta$ or $\Theta^2$ is inflationary and fixes $x \vee y$, $x \vee z$ and $y \vee z$, with
\[
x \leq \Theta(x) \leq (x \vee y) \wedge (x \vee z) \qquad\text{and}\qquad y \leq \Theta(y) \leq (x \vee y) \wedge (y \vee z)
\]
for the inflationary one.
\end{corEone}

\begin{proof}
$\Theta$ fixes a generator by Theorem~\ref{thm:E}, and two fixed generators would make it the identity (Lemma~\ref{lem:two-generators}); call the fixed one $z$.  By the partition lemma, $\Theta(x \vee y) \leq x \vee y$ or $z \leq \Theta(x \vee y)$; the latter would read $\Theta(z) \leq \Theta(x \vee y)$ and give $z \leq x \vee y$ by order reflection.  Dually $\Theta(x \wedge y) \geq x \wedge y$.  By Lemma~\ref{lem:descent}, a positive power $\Psi = \Theta^e$ fixes a join or a meet of two generators (not of three, as in the proof of Proposition~\ref{prop:inflationary-power}); say a join, else pass to the dual.  If it is $x \vee y$, Lemma~\ref{lem:1} applies to $\Psi$ and gives $\Psi$ or $\Psi^2$ inflationary, with $x \vee y$ fixed.  If it is $x \vee z$, Lemma~\ref{lem:1} at the join $x \vee z$ gives $x \leq \Psi(x) \leq x \vee z$ and $y \leq \Psi(y)$ (the swapped case, $x \leq \Psi(z) = z$, is impossible), so $\Psi$ is inflationary; then $x \vee y \leq \Psi(x \vee y) \leq x \vee y$, the second inequality because $\Theta(x \vee y) \leq x \vee y$ iterates, and $\Psi$ fixes $x \vee y$.  The case $y \vee z$ is symmetric.  In every case some positive power $\Theta^d$ fixes $x \vee y$, and the inequalities $x \vee y = \Theta^d(x \vee y) \leq \Theta^{d-1}(x \vee y) \leq \cdots \leq \Theta(x \vee y) \leq x \vee y$ (each is $\Theta^{k}$ applied to $\Theta(x \vee y) \leq x \vee y$) collapse, so $\Theta$ itself fixes $x \vee y$, and Lemma~\ref{lem:1} for $\Theta$ makes $\Theta$ or $\Theta^2$ inflationary; call it $\Theta\pr$.  For $\Theta\pr$: $y \leq \Theta\pr(x)$ would give $\Theta\pr(x) \geq x \vee y \geq \Theta\pr(x)$, so $\Theta\pr(x) = x \vee y \geq \Theta\pr(y)$ and $x \geq y$ by order reflection, which is false; so $y \not\leq \Theta\pr(x)$, and the partition lemma at $y$ puts $\Theta\pr(x) \leq x \vee z$.  Hence $\Theta\pr(x \vee z) = \Theta\pr(x) \vee z \leq x \vee z \leq \Theta\pr(x \vee z)$: $x \vee z$ is fixed, and $\Theta\pr(x) \leq (x \vee y) \wedge (x \vee z)$.  Likewise for $y$.
\end{proof}
\end{furthercut}

\begin{fullonly}
Corollary~\ref{cor:E1} is \agda{shape : Shape Θ} in \agda{DeansProblem.\allowbreak Shape}, and its formal proof is shorter than the one above.  Proposition~\ref{prop:inflationary-power} gives an inflationary positive power $\Theta^e$, after replacing $\Theta$ by its dual $\Theta^{\delta}$ if necessary, which fixes $z$ as well and so satisfies $\Theta^{\delta}(x \vee y) \leq x \vee y$ by the same argument.  Then $x \vee y \leq \Theta^e(x \vee y) \leq \Theta(x \vee y) \leq x \vee y$, the middle inequality because $\Theta(x \vee y) \leq x \vee y$ iterates, so $\Theta$ itself fixes $x \vee y$ with no case analysis on the join that Lemma~\ref{lem:descent} supplies; Lemma~\ref{lem:1} and the partition bounds finish as above.  The formal proof of Theorem~\ref{thm:G} does not use this corollary.
\end{fullonly}

\section{Rigidity, and the answer}
\label{sec:rigidity}

\begin{expertread}
Call a self-embedding $\Theta_1$ of $\FL(3)$ \emph{rigid} when no proper self-embedding fixes $\Theta_1(\FL(3))$ pointwise.  Proposition~\ref{prop:rigid} says that $\theta_1$ is rigid, and more: no endomorphism of $\FL(3)$ but the identity fixes $\theta_1(\FL(3))$ pointwise.  This section proves that every proper self-embedding with a fixed free triple is rigid, which is exactly the statement that no coherent chain of length three exists.  In such a chain $\Theta_0, \Theta_1, \Theta_2$, the link $\Theta_1$ is a proper self-embedding that fixes the free triple $(\Theta_0(x), \Theta_0(y), \Theta_0(z))$, and $\Theta_2$ is a proper self-embedding that fixes $\Theta_1(\FL(3))$ pointwise.  Conversely, if a proper self-embedding $\Theta_1$ fixes a free triple and a proper $\Theta_2$ fixes $\Theta_1(\FL(3))$ pointwise, then $\Theta_1$ and $\Theta_2$, with $\Theta_0$ the self-embedding that sends $x$, $y$, $z$ to that triple, form such a chain; $\Theta_0$ is proper, since a self-embedding that fixes a generating triple is the identity.
\end{expertread}

\begin{fullonly}
\subsection{A first form}
\label{sec:theorem-f}

The first form of rigidity we found asks for a meet node in a fixed point of the upper embedding, $\Phi$ below; it decided the examples found by search, and Lemma~\ref{lem:5} superseded it.  We state it because it is short.

\begin{thmF}[rigidity for a class]
\label{thm:F}
Let $\Theta_1 \colon x, y, z \mapsto a, b, z$ be a self-embedding fixing $z$ with $x \leq a$, $y \not\leq a$, $y \leq b$, $x \not\leq b$, and let $\Phi$ be a self-embedding fixing $a$, $b$ and $z$, hence every element of $\gen{a, b, z}$.  If the canonical tree of some fixed point of $\Phi$ has a meet node with $x$ or $y$ among its meetands, then $\Phi$ is the identity.
\end{thmF}

\begin{proof}
By Lemma~\ref{lem:4} applied to $a$ and to $b$, $x \leq \Phi(x)$ and $y \leq \Phi(y)$, and $\Phi(z) = z$; so $\Phi$ is inflationary.  By Lemma~\ref{lem:3}(a) it fixes every non-generator node of the canonical tree of its fixed point, in particular the meet node, and by Lemma~\ref{lem:3}(b) it then fixes the generator meetand of that node.  With $z$ that is two fixed generators, so $\Phi$ is the identity by Lemma~\ref{lem:two-generators}.
\end{proof}

For $\theta_1$ the node is $y \wedge z$, the meet joinand of $\theta_1(x)$, and Theorem~\ref{thm:F} recovers Proposition~\ref{prop:rigid}.  The search of Computation~\ref{comp:searches} found seven members of the class within complexity $4$, all with such a node; the seventh,
\[
\theta_2 \colon x \mapsto (x \vee (y \wedge (x \vee z))) \wedge (x \vee c), \qquad y \mapsto (y \vee (x \wedge (y \vee z))) \wedge (y \vee c), \qquad z \mapsto z,
\]
where $c = z \wedge (x \vee y)$ as in Theorem~\ref{thm:D}, has it as $y \wedge (x \vee z)$, and it fixes the free triple of Theorem~\ref{thm:D}, so it is a second coherent chain of length two, and it is rigid (the module \agda{DeansProblem.\allowbreak Examples.\allowbreak Rigidity}).  Whether every member of the class has a meet node with a moved generator in some fixed point of the upper embedding was open for a while; the answer is that no node is needed, because the upper embedding pins every element it fixes above one element.
\end{fullonly}

\subsection{Lemma~5 and Theorem~G}
\label{sec:lemma-5}

Recall that $c = z \wedge (x \vee y)$.

\begin{expertread}
\begin{lemfive}
\label{lem:5}
Let $\Psi$ be an inflationary self-embedding of $\FL(3)$ that fixes $z$ and moves $x$ and $y$.  Then every fixed point of $\Psi$ lies above $c$.
\end{lemfive}

\begin{proof}
By induction on the canonical form of a fixed point $w$.  If $w$ is a generator, it is $z$, since $x$ and $y$ are moved, and $c \leq z$.  Let $w = w_1 \vee \cdots \vee w_n$ canonically, $n \geq 2$.  If some $w_i$ is a meet, it is a non-generator node of the tree of $w$, so it is fixed by Lemma~\ref{lem:3}(a) and lies above $c$ by induction, hence so does $w$.  Otherwise every $w_i$ is a generator, and two of them are distinct; $c$ lies below every join of two distinct generators ($c \leq x \vee y$ by definition, and $c \leq z$ gives $c \leq x \vee z$ and $c \leq y \vee z$), hence $c \leq w$.  Let $w = w_1 \wedge \cdots \wedge w_n$ canonically, $n \geq 2$.  A generator meetand is fixed by Lemma~\ref{lem:3}(b), so it is $z$, which lies above $c$; a join meetand is a non-generator node, fixed by Lemma~\ref{lem:3}(a), and above $c$ by induction.  So every meetand lies above $c$, and so does $w$.
\end{proof}

Figure~\ref{fig:theorem-d} shows Lemma~\ref{lem:5} for $\theta_1$, which is inflationary, fixes $z$ and moves $x$ and $y$.

\begin{lemma}[the core]
\label{lem:core}
Let $\Theta_1$ be a self-embedding of $\FL(3)$ fixing $z$, and let $\Psi$ be a proper inflationary self-embedding fixing $\Theta_1(x)$, $\Theta_1(y)$ and $\Theta_1(z)$.  This is impossible.
\end{lemma}

\begin{proof}
Write $a = \Theta_1(x)$ and $b = \Theta_1(y)$; then $\Theta_1(z) = z$, and $(a, b, z)$ is a free triple since $\Theta_1$ is injective.  $\Psi$ fixes $z$.  It moves $x$ and $y$: if it fixed $x$, it would fix two generators and be the identity (Lemma~\ref{lem:two-generators}), which a proper self-embedding is not; likewise $y$ is moved by $\Psi$.  So Lemma~\ref{lem:5} applies to $\Psi$, and its fixed point $b$ lies above $c$.  On the other hand, the partition lemma at $z$ gives $\Theta_1(x \vee y) \leq x \vee y$ or $z \leq \Theta_1(x \vee y)$, and the latter would read $\Theta_1(z) \leq \Theta_1(x \vee y)$ and give $z \leq x \vee y$ by order reflection; so $a \leq \Theta_1(x \vee y) \leq x \vee y$.  Hence
\[
a \wedge z \leq (x \vee y) \wedge z = c \leq b.
\]
But $a \wedge z \not\leq b$ is one of the meet conditions~\eqref{eq:meet-condition} of Corollary~\ref{cor:1.13} for the free triple $(a, b, z)$.
\end{proof}

\begin{thmG}
\label{thm:G}
There is no coherent chain of length three.
\end{thmG}

\begin{proof}
Let $\Theta_0, \Theta_1, \Theta_2$ be one.  $\Theta_1$ fixes the free triple $(\Theta_0(x), \Theta_0(y), \Theta_0(z))$, so by Theorem~\ref{thm:E} it fixes a generator, which after renaming the generators we call $z$.  $\Theta_2$ fixes $\Theta_1(\FL(3))$ pointwise, in particular the free triple $(\Theta_1(x), \Theta_1(y), \Theta_1(z))$; so by Proposition~\ref{prop:inflationary-power} a positive power $\Psi$ of $\Theta_2$ or of $\Theta_2^{\delta}$ is inflationary and fixes that triple or its dual.  In the first case $\Psi$ is a proper inflationary self-embedding fixing $\Theta_1(x)$, $\Theta_1(y)$, $\Theta_1(z)$, and Lemma~\ref{lem:core} applies to $\Theta_1$ and $\Psi$.  In the second case, $\Psi$ fixes $\delta(\Theta_1(x))$, $\delta(\Theta_1(y))$, $\delta(\Theta_1(z))$, which are the images of $x$, $y$, $z$ under the dual self-embedding $\Theta_1^{\delta}$, and $\Theta_1^{\delta}$ fixes $z$ since $\Theta_1$ does; so Lemma~\ref{lem:core} applies to $\Theta_1^{\delta}$ and $\Psi$.  Either way there is a contradiction.
\end{proof}
\end{expertread}

\begin{fullonly}
So the rigidity of Theorem~\ref{thm:F} holds in a strong form, with no condition on the middle link: nothing proper fixes the image of a self-embedding with a fixed free triple pointwise.  The mechanism of Theorem~\ref{thm:D}, a fixed generator with the other two moved inside the ideal below their join, is the only one there is (Theorem~\ref{thm:E}, Corollary~\ref{cor:E1}), and it does not stack: the link above pins every element it fixes above $c$, and a free triple containing $z$ cannot lie there.
\end{fullonly}

\begin{expertread}
\begin{thmH}[Dean's problem]
\label{thm:H}
No free lattice contains an infinite strictly ascending chain of sublattices isomorphic to $\FL(3)$.
\end{thmH}

\begin{proof}
Such a chain yields, by the coherent form of Theorem~\formofC{coherent}, an infinite coherent chain $\Theta_0, \Theta_1, \Theta_2, \dots$, whose first three links form a coherent chain of length three, which Theorem~\ref{thm:G} forbids.
\end{proof}
\end{expertread}

The coherent form of Theorem~\formofC{coherent} chooses its subsequence by the infinite pigeonhole principle (Remark~\ref{rem:classical}).  Theorem~\ref{thm:H} does not need it: Remark~\ref{rem:constructive} proves Theorem~\ref{thm:H} from the four-link form of Theorem~\formofC{four} and Theorem~\ref{thm:G} by finite searches, and that is the proof checked in Agda.

\subsection{What is checked where}
\label{sec:which-step}

\begin{expertread}
Theorem~\ref{thm:G} is a theorem in the formal development (\agda{no-coherent₃} in \agda{DeansProblem.\allowbreak LengthThree}), with every lemma it uses, and so is its consequence that no infinite coherent chain exists (\agda{no-coherentω}).  So is Theorem~\ref{thm:H} as stated, as \agda{no-copy-chain} in \agda{DeansProblem.\allowbreak NoChain}, whose type
\begin{center}
{\small\agda{no-copy-chain : DecidableEquality X → CopyChain X → ⊥}}
\end{center}
says that the free lattice over a type $X$ of generators with decidable equality contains no chain in the sense of Definition~\ref{def:chain}, stated as a sequence of embeddings of $\FL(3)$ whose images ascend strictly (Section~\ref{sec:chains}).  It is \agda{no-chain}, which says the same of the record of a chain by free triples, composed with \agda{copyChain→chain}, which turns a chain of sublattices into that record and assumes nothing about $X$.  So decidable equality is the one hypothesis: it is what the development needs to compute complexity, since a canonical form compares generators, and classically every set has it.  The proof of \agda{no-chain} is that of Remark~\ref{rem:constructive}.  It uses no classical principle, no postulate and no hypothesis on the chain; nor does it use the coherent form of Theorem~\formofC{coherent}, which the development checks only with the infinite pigeonhole principle as a hypothesis (Section~\ref{sec:formal}).  The earlier formal statement, that no chain has a link whose transition from the bottom recurs beyond every bound (\agda{no-chain-recurs}), remains, as a consequence of Theorem~\ref{thm:G} and the length-three form of Theorem~\formofC{three}.  These are the statements of the release's tag \texttt{v0.1.0} (Section~\ref{sec:formal}).

\begin{remark}[the constructive proof]
\label{rem:constructive}
Fix a chain.  Equality in $\FL(3)$ is decidable (Theorem~\ref{thm:whitman}), and by Lemma~\ref{lem:bounded-transitions} the transitions $\theta_{0,k}$ take their values in the finite set $T$ of triples of elements of complexity at most $B(0)$, which has $N(0)$ elements.  Call a strictly increasing map $h \colon \mathbb{N} \to \mathbb{N}$ a \emph{selection}, and say that a sequence $u$ in a set with decidable equality \emph{escapes} if for every selection $h$ and all $r, m \in \mathbb{N}$ there is $s > m$ with $u(h(s)) \neq u(h(r))$.

\emph{The transitions from the bottom escape}.  Given $h$, $r$ and $m$, let $q = \max(r, m)$.  Some $s$ among $q + 1, \dots, q + N(h(r)) + 1$ has $\theta_{0,h(s)} \neq \theta_{0,h(r)}$, and since equality is decidable, a search through these $N(h(r)) + 1$ positions finds it.  Otherwise each of them selects a link $l > h(r)$ with $\theta_{0,l} = \theta_{0,h(r)}$.  Then $h(r) \geq 1$, since $\theta_{0,l}$ is proper and $\theta_{0,0}$ is the identity, and by Lemma~\ref{lem:bounded-transitions} two of these $N(h(r)) + 1$ links, $l_1 < l_2$, have $\theta_{h(r),l_1} = \theta_{h(r),l_2}$.  With $\theta_{0,h(r)} = \theta_{0,l_2}$, the four links $0 < h(r) < l_1 < l_2$ satisfy the hypotheses of the four-link form of Theorem~\formofC{four}, and the coherent chain of length three it gives contradicts Theorem~\ref{thm:G}.

\emph{No sequence in a finite set with decidable equality escapes}.  By induction on the number of elements of the set.  If the set is empty, there is no $u(0)$.  Otherwise, take for $h$ the identity and $r = 0$: the escape of $u$ gives $g(0) > 0$ with $u(g(0)) \neq u(0)$, then $g(1) > g(0)$ with $u(g(1)) \neq u(0)$, and so on, a selection $g$ with $u(g(t)) \neq u(0)$ for every $t$.  So $u \circ g$ takes its values in the set without $u(0)$, which has one element fewer, and it escapes, because a selection $h$ from $u \circ g$ is the selection $g \circ h$ from $u$.  The induction hypothesis rules it out.

The two together prove Theorem~\ref{thm:H}, since the transitions from the bottom of a chain would be a sequence in $T$ that escapes.  The argument uses decidable equality in $\FL(3)$, finite searches and induction on a natural number, and nothing classical.  It is the proof of \agda{no-chain}, in which \agda{recurrences\allowbreak -bounded} and \agda{chain\allowbreak -escapes} are the first part and \agda{finite\allowbreak -colors\allowbreak -exhausted} is the second.  There a triple is coded by the positions of its coordinates in a finite list that contains an element equal to each element of complexity at most $B(j)$, the enumeration in the proof of Lemma~\ref{lem:finite}; so $N(j)$ becomes the cube of the length of the list, an upper bound for the number of triples, and the induction runs over the codes.  Equal codes give equal triples, which is all the argument uses.  Counting gives a second constructive proof, which the development does not check.  By the argument of the first part, the value of $\theta_{0,\cdot}$ first seen at link $j$ recurs at most $N(j)$ times after $j$; so the pair (number of values not yet seen, number of recurrences still allowed for the values seen) decreases in the lexicographic order of $\mathbb{N} \times \mathbb{N}$ at every link, and that order is well founded.  Both proofs use bounds $N(j)$ that depend on the links of the chain, as they must: finite chains of every length exist (Section~\ref{sec:problem}), so no bound independent of the chain limits their length.
\end{remark}
\end{expertread}

\begin{furthercut}
\section{The almost-proof}
\label{sec:almost-proof}

The notes \cite{Freese1997} end with ``an almost proof that an ascending chain cannot exist.''  Before the argument of Sections~\ref{sec:coherent} to~\ref{sec:rigidity} was found, we formalized it up to its gap, stated the gap as a decidable predicate, and decided the predicate on the one-step example.  This section is the record of that path; nothing in it is used by the proof of Theorem~\ref{thm:H}, whose finiteness argument sharpens that of the notes and whose rigidity argument replaces their Lemma~2.

\subsection{The argument}
\label{sec:almost-argument}

\begin{expertread}
In the words of the notes: ``whenever the sublattice generated by $a$, $b$, and $c$ is contained [in] that generated by $a\pr$, $b\pr$, and $c\pr$ and both are free lattices, then the primed generators determine a truncation of the term tree of $a$, $b$, and $c$.  This means that for each path of the tree there is an element and an immediate subelement such that one of $a\pr$, $b\pr$ or $c\pr$ lies in the interval between the element and the subelement.  There are only finitely many possible truncations even counting the labelling.  Thus [if] we had a[n] ascending chain $\mathrm{Sb}(a_k, b_k, c_k)$ of free sublattices we could assume that the truncation pattern of $\{a\pr, b\pr, c\pr\}$ determined by $\{a_k, b_k, c_k\}$ is the same for all $k$.''  The pattern here, made precise as a labeled pattern in Section~\ref{sec:truncation-theorem}, is the one that each link $\{a_k, b_k, c_k\}$ determines in the trees of the \emph{bottom} triple, whose elements stay fixed along the chain; only so are the bounding elements the same for every link.

\begin{quote}
\textbf{Lemma~2 of the notes}.  Suppose $a^- \leq a^+$, $b^- \leq b^+$, and $c^- \leq c^+$ are strongly join and meet irredundant, i.e., $a^- \not\leq b^+ \vee c^+$ and symmetrically and dually.  If $a^- \leq a \leq a^+$ and similarly for $b$ and $c$, then $a$ is the only element of the sublattice generated by $a$, $b$, and $c$ in $a^+/a^-$.

\emph{Proof.}  Suppose $a^- \leq w$ and $w$ is in the sublattice generated by $a$, $b$, and $c$.  Then either $w \geq a$ or $w \leq b \vee c$.  The latter gives $a^- \leq w \leq b \vee c \leq b^+ \vee c^+$, a contradiction.  So $w \geq a$ and similarly $w \leq a$.
\end{quote}

With one pattern recurring along a subchain, the elements of the trees of the bottom triple that bound the primed generators are the same for every link, so every $a_k$ lies in one interval $[a^-, a^+]$, and similarly for $b_k$ and $c_k$.  Under the hypothesis of Lemma~2 of the notes, $a_k$ is the only element of $\gen{a_k, b_k, c_k}$ in that interval; but for $j < k$ the element $a_j$ of $\gen{a_k, b_k, c_k}$ lies there too, so $a_j = a_k$, and the chain is constant.  The gap is the hypothesis: nothing forces the bounding elements of a recurring pattern to be strongly irredundant.
\end{expertread}

\subsection{The truncation theorem}
\label{sec:truncation-theorem}

Lemma~\ref{lem:9.13} applies at every level of the canonical tree, and the result is a theorem of the formal development (\agda{truncation} and \agda{truncation-canonical} in \agda{DeansProblem.\allowbreak Truncation}).  Let $\tau \colon \FL(I) \to \FL(Y)$ be an embedding, $p$ a term over $I$ with canonical form $C$, and $D$ the canonical form of $\tau(p)$.  Say that $C$ \emph{truncates} $D$ (through $\tau$), written $C \sqsubseteq D$, when the following hold recursively.  A generator $i$ truncates $D$ when $\tau(i)$ is the element named by $D$.  A canonical join $C$ truncates a canonical join $D$ when they name the same element through $\tau$, every meet child of $C$ truncates some child of $D$ (its \emph{match}), and every child of $D$ that is not matched lies below $\tau(i)$ for a generator child $i$ of $C$ (a \emph{cut}, with $\tau(i)$ below the node $D$); dually at meet nodes.

\begin{theorem}[the truncation theorem]
\label{thm:truncation}
With $\tau$, $p$, $C$ and $D$ as above, $C \sqsubseteq D$.
\end{theorem}

\begin{proof}
By induction on the two canonical terms, as in the proof of the uniqueness of canonical forms.  At a join node, Lemma~\ref{lem:9.13} (in the form that every non-generator canonical joinand of $p$ is sent to a canonical joinand of $\tau(p)$) supplies the match of each meet child, and Theorem~\ref{thm:refine}, applied to $\tau(p) = \bigvee_i \tau(p_i)$ with $p_i$ the children of $C$, puts each child of $D$ below some $\tau(p_i)$; a child of $D$ below $\tau(p_i)$ for a meet child $p_i$ is the match of $p_i$ (two canonical joinands of the same element that are comparable are equal), and one below $\tau(i)$ for a generator child $i$ is a cut.  The induction hypothesis handles the matched pairs.  Meet nodes are dual.
\end{proof}

This is the sentence of the notes ``the term in $a\pr$, $b\pr$, and $c\pr$ giving $a$ is just a truncation of the term in the free generators giving $a$,'' with the correspondence of nodes and the labels of the cuts recorded.  A \emph{labeled pattern} is the data the relation records: which node of $D$ each node of $C$ names, and at each cut which $\tau(i)$ sits between the node and its child.  For a fixed canonical tree $D$ and a finite $I$ there are finitely many labeled patterns on $D$, whatever $\tau$ and $p$.  With $D$ the tree of an element of the bottom triple, which stays fixed along the chain, and $I$ the three generators, this is the ``finitely many possible truncations even counting the labelling'' of the notes.

\subsection{Lemma~2 of the notes in any lattice, and the gap as a predicate}
\label{sec:gap}

Lemma~2 of the notes holds in any lattice, with the proof given there.  The claim ``either $w \geq a$ or $w \leq b \vee c$'' is the partition lemma for the sublattice generated by $a$, $b$, $c$: by induction on a term in $a, b, c$ naming $w$, a join is above $a$ when either joinand is and below $b \vee c$ when both are, and a meet is below $b \vee c$ when either meetand is and above $a$ when both are.  Read at one index, the lemma needs only the two irredundance conditions at that index.  In the formal development: \agda{above-or-below}, \agda{lemma-2-at} and \agda{lemma-2} in \agda{DeansProblem.\allowbreak Irredundant}, stated for any lattice and any family of three elements, and specialized to $\FL(X)$ in \agda{DeansProblem.\allowbreak AlmostProof}.

A \emph{box} in $\FL(X)$ is a triple of intervals $[l_i, h_i]$, $i \in \{x, y, z\}$.  It is \emph{irredundant at} $i$ when $l_i \not\leq h_j \vee h_k$ and $l_j \wedge l_k \not\leq h_i$, where $\{j, k\}$ are the other two indices, and \emph{strongly irredundant} when it is irredundant at every index; narrowing a box, raising lower ends or lowering upper ends, keeps irredundance.  A triple $\tau$ \emph{lies in} the box when $l_i \leq \tau(i) \leq h_i$ for each $i$.

\begin{theorem}[\agda{box-agree}, \agda{no-step-in-box}]
\label{thm:box-agree}
Let $\sigma\pr$ and $\tau\pr$ be free triples in $\FL(X)$ lying in one box, with $\sigma\pr$ in the sublattice generated by $\tau\pr$.  Then $\sigma\pr(i) = \tau\pr(i)$ at every index $i$ at which the box is irredundant.  Hence no step $\gen{\sigma} \subsetneq \gen{\tau}$ has both of its triples in a box irredundant at the generator the step's terms miss, and at most one link of a chain lies in a strongly irredundant box.
\end{theorem}

\begin{proof}
$\sigma\pr(i)$ is an element of $\gen{\tau\pr}$ in the interval $[l_i, h_i]$, so by Lemma~2 of the notes at the index $i$ it is $\tau\pr(i)$.  For a step, equality at the missed index $f$ would put $\tau(f) = \sigma(f)$ in $\gen{\sigma}$, and then the generator $f$ of $\FL(3)$ would lie in the sublattice generated by the step's terms, hence be one of them by Corollary~\ref{cor:1.6}.  Two links of a chain in a strongly irredundant box would agree at every index, and so generate the same sublattice, against the strictness of the chain.
\end{proof}

\begin{expertread}
The \emph{box of a labeled pattern} is built from its cuts: at each cut the child lies below $\tau(i)$, which lies below the node, so the cut children bound $\tau(i)$ from below and the nodes from above (dually at meet nodes).  Every link with that pattern lies in the box.  So the argument of the notes closes on any recurring labeled pattern whose box is irredundant at a generator missed by the terms of the step between two of its links; \textbf{the gap is the predicate ``strongly irredundant'' on the box of a labeled pattern that recurs}, decidable on explicit boxes (\agda{stronglyIrredundant?}).
\end{expertread}

The gap does not move inside $\FL(3)$.  The recurring datum of the one-step form of Theorem~\ref{thm:C}, the transition $\theta_{0,j}$, fixes the terms of the truncation pattern.  The rest of the labeled pattern is a further datum on the fixed trees of the bottom triple, and by Section~\ref{sec:truncation-theorem} the labeled patterns that a link can have on those trees are finitely many; let $P$ be their number.  The pigeonhole of the one-step form, run on the transition and the labeled pattern together, then finds two links among the first $N(0) \cdot P + 1$ that agree in both.  Two links are all that Theorem~\ref{thm:box-agree} needs, so the gap of the almost-proof is a property of the box of their common pattern.  The ends of the box are elements of the trees of the bottom triple in the ambient lattice, which in general lie outside every link: the generator $a_1$ of the one-step example, which lies below $a\pr$ in $\FL(9)$, lies outside $\gen{a\pr, b\pr, c\pr}$ by Corollary~\ref{cor:1.6}.  Inside a link read as $\FL(3)$ the cut children are the generators themselves, and the box of a generator that is both a joinand and a meetand of the terms shrinks to a point.
\end{furthercut}

\begin{fullonly}
\subsection{The verdict on the one-step example}
\label{sec:verdict}

On the one-step example the predicate holds, so its step is not a witness of the gap (the module \agda{DeansProblem.\allowbreak Examples.\allowbreak FreeseBox}, about two seconds).  The box of the step, $[\ba{a}_1 \vee \ba{b}_2 \vee \ba{c}_2, \ a \wedge \ba{b}_{21} \wedge \ba{c}_{22}]$ and its images under the permutation, lies inside the \emph{generator box}
\[
[a_1, \ a_1 \vee b_{21} \vee c_{22}], \qquad [b_1, \ b_1 \vee c_{21} \vee a_{22}], \qquad [c_1, \ c_1 \vee a_{21} \vee b_{22}],
\]
which is strongly irredundant on its face, since a generator is not below a join of other generators and a meet of two generators is not below a join of generators excluding them (six Whitman decisions on terms with at most six generator occurrences); narrowing keeps the property, so the box of the step is strongly irredundant, and the step lies in it, which needs the ``lucky'' inequalities $\ba{b}_2, \ba{c}_2 \leq a$ of the footnote of the notes (Section~\ref{sec:relations}).  More is true, through the truncation theorem: \emph{every} free triple $\tau$ with $(a, b, c) = (\tau(x) \vee (\tau(y) \wedge \tau(z)), \dots)$, that is, every triple that expresses the base $(a, b, c)$ of the example by the terms of the example, lies in the generator box, whatever the rest of its pattern (\agda{freese-terms-in-box}).  The lower ends come from join primality: $a_1 \leq a = \tau(x) \vee (\tau(y) \wedge \tau(z))$ and $a_1 \not\leq b$ put $a_1 \leq \tau(x)$.  The upper ends come from truncating the canonical forms of $a$, $b$, $c$, written out and certified canonical by the decision procedure for the conditions of Theorem~\ref{thm:canonical}: the match of the meet joinand of the terms is forced, and its meetand $b_{21} \vee c_{22} \vee a_1$ lies above $\tau(x)$.  Hence two such triples, one inside the other's sublattice, are equal (\agda{freese-terms-once}): in a chain over that base, the terms of the example are the transition $\theta_{0,k}$ for at most one $k$.

\subsection{What survives}
\label{sec:survives}

The finiteness argument of the almost-proof survives as Theorem~\ref{thm:C}, sharpened from a truncation pattern to the whole term, because an embedding never lowers complexity (Theorem~\ref{thm:B}); Lemma~2 of the notes and the gap are never used.  What remained for a proof along the lines of the notes was the predicate for every labeled pattern that can recur over a base; where it fails, at some index either $l_i \leq h_j \vee h_k$ or $l_j \wedge l_k \leq h_i$, and the failure lives in the cut children.  The proof found instead does not look at boxes at all: it reads the chain inside $\FL(3)$ as a coherent chain and shows that its link $\Theta_2$ cannot fix the image of $\Theta_1$.
\end{fullonly}

\section{The formal development}
\label{sec:formal}

Theorem~\ref{thm:H}, every part of its proof, and every other result of this paper are theorems of a development in the proof assistant Agda \cite{Agda}, which we describe here: what ``checked'' means, where each statement of the paper lives and in what form, and what the checking costs.  Each result is one Agda definition whose type is its statement (Table~\ref{tab:names}), and each is checked as stated, over generators with decidable equality where a decision needs it, but two: the coherent form of Theorem~\formofC{coherent} and Corollary~\ref{cor:dean-iff-coherent} take the infinite pigeonhole principle as a hypothesis, since Agda's constructive logic does not prove it, and so does the recurrence in infinitely many links of Remarks~\ref{rem:what-survives} and~\ref{rem:classical}.  The checked proof of Theorem~\ref{thm:H} uses none of them (Section~\ref{sec:which-step}).  Two facts that the text states without proof are checked in part: of Section~\ref{sec:shape}~(iii), one direction of the characterization of the free generators as the elements both join and meet prime, and not the conclusion that a copy of $\FL(3)$ determines its free triple; and of the two consequences of the canonical form stated after Theorem~\ref{thm:refine}, the second, in $\FL(n)$.  What stays in prose is the following: the second, counting proof of Remark~\ref{rem:constructive}; the pigeonhole of Section~\ref{sec:gap} on the transition and the labeled pattern together, of which only the finiteness of the patterns is checked; the claim of the notes that the new generators of the one-step example may be chosen as complex as one likes, and the comparison of its canonical form with the display in the notes (Section~\ref{sec:freese-example}); and three remarks that no proof uses, that complexity is intrinsic (Section~\ref{sec:complexity}) and that finite links of the chain of Theorem~\ref{thm:ladder} would have unbounded size (Computation~\ref{comp:ladder}) and be pairwise non-isomorphic (Section~\ref{sec:questions}).  Table~\ref{tab:by-section} counts the claims of each section by how they are checked.  Theorem~\ref{thm:H} is checked as stated, for chains of sublattices as Definition~\ref{def:chain} defines them (Section~\ref{sec:which-step}).  The development is released with this paper as the repository \agda{agda-deans-problem} \cite{agda-deans-problem}, which holds this paper's modules, under the namespace \agda{DeansProblem}, and takes the general lattice theory they use from the library \agda{agda-lattices} \cite{agda-lattices}, pinned at one commit.  The names and the counts below, and the costs of checking the release, are those of the release's tag \texttt{v0.1.0} (commit \texttt{d562c05}, of 19 November 2026), with the library at the commit it pins.

\begin{fullonly}
\subsection{What is checked, and how}
\label{sec:what-checked}

The release consists of this paper's $38$ modules of literate Agda, $16\,673$ lines of source of which $10\,575$ lines are Agda code and the rest prose, and takes $42$ modules from the library, those they import, $16\,821$ lines of which $10\,375$ are code; the library holds $91$ modules at the pinned commit, and the rest of it is not used here.  They are checked by a development version of Agda 2.9.0 with version 2.3 of its standard library, patched to check under that version \cite{agda-stdlib}, on top of the Agda Universal Algebra Library \cite{agda-algebras} for the notions of lattice, homomorphism and equational theory.  Every module carries the pragma
\begin{center}
\agda{\{-\# OPTIONS -{}-without-K -{}-exact-split -{}-safe \#-\}}
\end{center}
and the development has no postulate and no termination pragma: under \agda{-{}-safe} Agda refuses postulates, unsafe primitives and non-terminating definitions, so every theorem is proved from the type theory alone.

The construction is syntactic.  Lattice terms over a type $X$ of generators are the inductive type \agda{LatTerm X} with constructors \agda{ℊ} (a generator), \agda{\_∨̇\_} and \agda{\_∧̇\_}; Whitman's rules of Theorem~\ref{thm:whitman} are an inductive relation \agda{\_≤ʷ\_} on terms, one constructor per rule and per disjunct, so that an inhabitant of \agda{s ≤ʷ t} is a derivation by the procedure; equality in the free lattice, \agda{s ≈ʷ t}, is the conjunction of \agda{s ≤ʷ t} and \agda{t ≤ʷ s}.  The rules are proved sound and complete for the order of every lattice (\agda{≤ʷ→≤}, \agda{≤→≤ʷ}), so the quotient of the terms by \agda{\_≈ʷ\_} is the free lattice \agda{FL X} with its universal property (\agda{FL-extend}), rather than a model read off Day's doubling construction as in FJN.  Canonical forms are a separate inductive type of $n$-ary terms, with the canonicalization algorithm (\agda{canonical}) and the conditions of Theorem~\ref{thm:canonical} (\agda{IsCanonical}, decided by \agda{isCanonical?}).  Theorems~1.17 and~1.18 are checked in full (\agda{FJN-1·17} and \agda{FJN-1·18} in \agda{Lattices.\allowbreak Free.\allowbreak Canonical.\allowbreak LeastRank}): every element has a term of least rank, every term of least rank naming it agrees with that one up to permutation (\agda{\_≈ᶜ\_}), and a term has least rank exactly when it satisfies the conditions.  The rank is that of Section~\ref{sec:terms}, on the $n$-ary terms in which no join has a join as a child and no meet a meet, which hold every term of least rank, since removing the parentheses around a join inside a join lowers the rank by one.  Theorem~1.19 is \agda{canonicalJoinands-refine}.  A substitution \agda{σ : Sub X Y} is a map from $X$ to terms over $Y$; \agda{t [ σ ]} is the term with the generators replaced, \agda{σ ⨟ τ} the composite, and \agda{iter σ k} the $k$-th power.  A family freely generates, \agda{FreelyGenerates σ}, when it satisfies the join condition~\eqref{eq:join-condition} and the meet condition~\eqref{eq:meet-condition} of Corollary~\ref{cor:1.13} on lists of indices, and that is equivalent to the injectivity of the homomorphism it induces (\agda{freelyGenerates⇔injective}).

``Checked'' means one of two things.  A statement proved by structural induction is checked when its proof term has the stated type.  A statement about explicit terms is checked by computation: Whitman's procedure is implemented as a decision \agda{s ≤ʷ? t} of type \agda{Dec (s ≤ʷ t)}, which Agda evaluates while it checks the module, and the result is turned into a proof by \agda{from-yes} (a positive answer) or \agda{from-no} (a negative one); so a theorem such as ``the triple of Theorem~\ref{thm:D} freely generates'' holds because the procedure ran, during type checking, on the six inequalities of Corollary~\ref{cor:1.13}.  Free generation of a family over a finite type is decided the same way (\agda{freelyGenerates?}).  The decision \agda{s ≤ʷ? t} is the recursion of Theorem~\ref{thm:whitman}, exponential in the worst case, since it may decide one pair of subterms many times (FJN, Example~11.32).  The library decides the order a second way, as \agda{s ≤ᵖ? t} in \agda{Lattices.\allowbreak Free.\allowbreak Whitman.\allowbreak Polynomial}, by the table of FJN, Chapter~XI, Section~8: one entry for each pair of a subterm of $s$ and a subterm of $t$, each entry built by the rules of Theorem~\ref{thm:whitman} from at most four others, so that the cost is proportional to the product of the ranks of $s$ and $t$.  Each entry is itself a decision built from the rules, so the table needs no proof of correctness of its own, and its verdicts are proved to be those of the recursion (\agda{≤ᵖ?-agrees}).  The decisions of items (i) to (iv) of Section~\ref{sec:verification} run through the table, on smaller terms checked equal to the defined ones, with the statements transported along equality, and so do a few more on the terms of the example, in \agda{DeansProblem.\allowbreak Examples.\allowbreak FreeseStep.\allowbreak Complexity} and \agda{DeansProblem.\allowbreak Examples.\allowbreak FreeseGap}; every other decision of the release runs through the recursion, those on the canonical form of $a$ of item (v) and of Section~\ref{sec:verdict} among them.

Table~\ref{tab:by-section} counts the claims of each section of this paper, a claim being every numbered or lettered result, every assertion of a numbered remark or computation, and every other claim the text states with a proof, a reason or a source, counted once where it is first so stated, and each claim being checked as stated (F, with decidable equality on the generators where a decision needs it), checked in part (P: in a narrower form, in one direction, or from the infinite pigeonhole principle taken as a hypothesis), left in prose (N), or computed by a search and not checked (E).

\begin{table}[!tbp]
\centering
\caption{What the development checks, section by section, at the release's tag \texttt{v0.1.0}: the claims of each section, and how many are checked as stated (F), checked in part (P), left in prose (N), or computed and not checked (E).}
\label{tab:by-section}
\footnotesize
\setlength{\tabcolsep}{4pt}
\begin{tabular}{@{}>{\raggedright}p{0.19\textwidth}rrrrr>{\raggedright\arraybackslash}p{0.48\textwidth}@{}}
\toprule
Section & Claims & F & P & N & E & What is checked, and what is not \\
\midrule
1 The problem & 5 & 4 & 1 & 0 & 0 & Theorem~\ref{thm:H} as stated; the reduction to $\FL(3)$, the descending chain, and the complexities of the printed step; of (iii), the primality of the generators and the automorphisms, not the rest \\
2 Preliminaries & 16 & 14 & 1 & 1 & 0 & every result, Theorems~\ref{thm:canonical} and~\ref{thm:refine} in full; the two consequences of the canonical form in part; that complexity is intrinsic, not \\
3 The one-step example & 9 & 7 & 0 & 2 & 0 & every claim of Section~\ref{sec:verification} and where $a\pr$ lies; not the claim of the notes that the new generators may be chosen arbitrarily complex, nor the comparison with their display \\
4 Chains & 13 & 11 & 0 & 1 & 1 & every result, Proposition~\ref{prop:reformulation} and Theorem~\ref{thm:ladder} whole, and the claims around them; the sizes of the links computed, and what would follow if all are finite, in prose \\
5 No step is self-similar & 10 & 9 & 1 & 0 & 0 & every result and claim; the recurrence of Remark~\ref{rem:what-survives} from the infinite pigeonhole principle \\
6 Coherent chains & 20 & 15 & 3 & 0 & 2 & every result, the coherent form of Theorem~\formofC{coherent} and Corollary~\ref{cor:dean-iff-coherent} from the infinite pigeonhole principle, as is the recurrence of Remark~\ref{rem:classical}; the counts after Theorem~\ref{thm:D} and of Computation~\ref{comp:fixed} computed \\
7 A fixed free triple & 13 & 13 & 0 & 0 & 0 & every claim \\
8 Rigidity, and the answer & 9 & 8 & 0 & 1 & 0 & every result and claim, and the constructive proof of Remark~\ref{rem:constructive}; not its second proof, by counting \\
9 The almost-proof & 11 & 10 & 0 & 1 & 0 & every claim but the pigeonhole of Section~\ref{sec:gap} on the transition and the labeled pattern together \\
10 The formal development & 1 & 0 & 0 & 0 & 1 & the searches of Computation~\ref{comp:searches}, computed \\
11 Remarks and questions & 1 & 0 & 0 & 1 & 0 & the questions are open; that finite links of Theorem~\ref{thm:ladder} would be pairwise non-isomorphic, not \\
\midrule
Total & 108 & 91 & 6 & 7 & 4 & \\
\bottomrule
\end{tabular}
\end{table}
\end{fullonly}

\subsection{How the development and the proofs were produced}
\label{sec:how-produced}

The Agda development, the computer searches reported in Computations~\ref{comp:ladder}, \ref{comp:fixed} and~\ref{comp:searches}, and the first draft of this paper were produced by large language models working as agents under the direction of the first author, between 6 and 30 October 2026: Claude Fable~5.1 and Claude Opus~5.5 (Anthropic), in Claude Code; GPT-6-Astra (OpenAI), in the Codex CLI; and, for part of the canonical forms, Kimi~K3 (Moonshot AI), in Kimi Code CLI.  The proofs of the new results were worked out in the same sessions.  The problem is the one that the second author's notes of 1997 describe \cite{Freese1997}, and the example of Section~\ref{sec:freese-example} and the argument of Section~\ref{sec:almost-proof} are the second author's.  The first author set the question each agent session took up and suggested approaches to try.  Both authors decided which results to keep.  The agents developed most of the arguments, starting from the second author's notes, found the examples by search, and formalized the proofs.  The first author directed and refined the Agda formalization.  Both authors revised the manuscript, and the second author corrected some arguments and replaced others with better proofs.

The agents wrote the Agda against the type checker through the \agda{agda-mcp} server of the Agda-native AIR project \cite{agda-native-air}, which exposes Agda's type checker, its goal state and a search over the library's corpus to an agent as tools; a field report on the server measures its use in these sessions \cite{agda-mcp-field-report}.

The mathematics was checked in two ways: every statement this section reports as checked, by Agda's type checker under \agda{-{}-safe}; and every proof in this paper, by the authors, who read each against the Agda statement it reports or, for the statements that stay in prose (named above), on its own.  The searches guided the work but are cited by no proof.  The authors take responsibility for every statement and proof in this paper.

\begin{fullonly}
\subsection{Where each statement lives}
\label{sec:table}

The development uses three Agda repositories, each building on the one before.  The Agda Universal Algebra Library, \agda{agda-algebras} \cite{agda-algebras}, gives the general notions of lattice, homomorphism and equational theory.  The library \agda{agda-lattices} \cite{agda-lattices} develops lattice theory on that base: free lattices, Whitman's solution of the word problem and canonical forms, with results of FJN, most of them from its first chapter, under FJN's numbers.  The release, \agda{agda-deans-problem} \cite{agda-deans-problem}, holds this paper's own modules, under the namespace \agda{DeansProblem}, and takes the lattice theory it needs from \agda{agda-lattices}, pinned at one commit.  It is a repository of its own so that every result of this paper can be found, formalized, in one place: its module \agda{DeansProblem.\allowbreak Paper} collects the statement and the proof of each result, in the paper's order and under the paper's numbers.

In Agda a statement is a type, and a proof of it is a definition of that type.  Table~\ref{tab:names} gives, for each result of this paper, the name of its statement and the module that defines it; its proof has the same name with a lower-case initial, as in \agda{theorem-4·5 : Theorem-4·5}.  A module whose name begins with \agda{DeansProblem} is in the release, and one whose name begins with \agda{Lattices} is in \agda{agda-lattices}.  Two of the release's statements restate the library's statements of results of FJN: \agda{Theorem-2·3} is defined as \agda{WhitmansTheorem}, in \agda{Lattices.\allowbreak Free.\allowbreak Universal}, and \agda{Corollary-2·5} has the definition of \agda{FJN-1·13}, in \agda{Lattices.\allowbreak Free.\allowbreak Generation}.  The table ends with the definitions of Sections~\ref{sec:chains}, \ref{sec:coherent} and~\ref{sec:structure}.  Figure~\ref{fig:modules} shows how this paper's modules import one another.

\begin{figure}[ht]
\centering
\begin{tikzpicture}[x=10mm, y=9mm]
  \node[module, minimum width=64mm] (Lib) at (-.2,0) {the library \agda{agda-lattices}};
  \node[ours] (Truncation) at (-4.3,-1.1) {\agda{Truncation}};
  \node[ours] (Irredundant) at (-1.9,-1.1) {\agda{Irredundant}};
  \node[ours] (Chain) at (-0.3,-1.1) {\agda{Chain}};
  \node[ours] (Ladder) at (1.15,-1.1) {\agda{Ladder}};
  \node[ours] (Preliminaries) at (2.94,-1.1) {\agda{Preliminaries}};
  \node[ours] (Partition) at (5,-1.1) {\agda{Partition}};
  \node[ours, align=center] (Truncation-Pattern) at (-4.8,-2.20) {\agda{Truncation.}\\\agda{Pattern}};
  \node[ours] (AlmostProof) at (-2.4,-2.1) {\agda{AlmostProof}};
  \node[ours] (SelfSimilar) at (-1.8,-2.75) {\agda{SelfSimilar}};
  \node[ours] (Chain-Copies) at (0,-2.2) {\agda{Chain.Copies}};
  \node[ours] (Coherent) at (1.8,-2.2) {\agda{Coherent}};
  \node[ours] (Fixed) at (3.58,-2.2) {\agda{Fixed}};
  \node[ours] (Chain-Reduction) at (-1.5,-5.5) {\agda{Chain.Reduction}};
  \node[ours] (Chain-Limit) at (-1.65,-3.35) {\agda{Chain.Limit}};
  \node[ours] (Chain-Powers) at (3.4,-3.3) {\agda{Chain.Powers}};
  \node[ours] (Inflationary) at (5.5,-3.3) {\agda{Inflationary}};
  \node[ours] (Chain-Condition) at (0.9,-3.35) {\agda{Chain.Condition}};
  \node[ours] (FixedGenerator) at (5.2,-4.4) {\agda{FixedGenerator}};
  \node[ours, align=center] (AlmostProof-Links) at (-3.2,-4.40) {\agda{AlmostProof.}\\\agda{Links}};
  \node[ours, align=center] (Chain-TwoGenerated) at (0.3,-4.40) {\agda{Chain.}\\\agda{TwoGenerated}};
  \node[ours] (Coherent-Finite) at (2.7,-4.4) {\agda{Coherent.Finite}};
  \node[ours] (Shape) at (4.75,-5.50) {\agda{Shape}};
  \node[ours] (Pigeonhole) at (-4.9,-5.3) {\agda{Pigeonhole}};
  \node[ours] (Coherent-Infinite) at (-1.3,-6.5) {\agda{Coherent.Infinite}};
  \node[ours] (LengthThree) at (2.8,-5.5) {\agda{LengthThree}};
  \node[ours] (NoChain) at (2,-6.55) {\agda{NoChain}};
  \begin{scope}[on background layer]
    \draw[imports] ([shift={(-2.6,0)}]Lib.center) -- (Truncation);
    \draw[imports] ([shift={(-1,0)}]Lib.center) -- (Irredundant);
    \draw[imports] (Lib) -- (Chain);
    \draw[imports] ([shift={(0.9,0)}]Lib.center) -- (Ladder);
    \draw[imports] ([shift={(1.7,0)}]Lib.center) -- (Preliminaries);
    \draw[imports] ([shift={(2.2,0)}]Lib.center) -- (Partition);
    \draw[imports] ([shift={(-2.3,0)}]Lib.center) -- (Pigeonhole);
    \draw[imports] (Truncation) -- (Truncation-Pattern);
    \draw[imports] (Irredundant) -- (AlmostProof);
    \draw[imports] (Chain) -- (AlmostProof);
    \draw[imports] (Chain) -- (SelfSimilar);
    \draw[imports] (Chain) -- (Chain-Copies);
    \draw[imports] (Chain) -- (Coherent);
    \draw[imports] (Chain) -- (Fixed);
    \draw[imports] (Chain-Copies) -- (Chain-Reduction);
    \draw[imports] ([shift={(-0.1,0)}]Chain-Copies.center) -- ([shift={(0.65,0.15)}]Chain-Limit.center);
    \draw[imports] (Chain-Copies) -- (Chain-Powers);
    \draw[imports] (Coherent) -- (Chain-Powers);
    \draw[imports] (Chain-Copies) -- (Chain-Condition);
    \draw[imports] (Fixed) -- (Inflationary);
    \draw[imports] (Partition) -- (Inflationary);
    \draw[imports] (Inflationary) -- (FixedGenerator);
    \draw[imports] ([shift={(-0.4,-0.1)}]AlmostProof.center) -- (AlmostProof-Links);
    \draw[imports] (Chain-Condition) -- (AlmostProof-Links);
    \draw[imports] (Chain-Condition) -- (Chain-TwoGenerated);
    \draw[imports] (Chain-Condition) -- (Coherent-Finite);
    \draw[imports] (Coherent) -- (Coherent-Finite);
    \draw[imports] (FixedGenerator) -- (Shape);
    \draw[imports] (FixedGenerator) -- (LengthThree);
    \draw[imports] (Coherent-Finite) -- (LengthThree);
    \draw[imports] (Coherent-Finite) -- (Coherent-Infinite);
    \draw[imports] (Pigeonhole) -- (Coherent-Infinite);
    \draw[imports] (LengthThree) -- (NoChain);
  \end{scope}
\end{tikzpicture}
\caption{This paper's modules, with the namespace \agda{DeansProblem} dropped; the examples and the index \agda{Paper} are not drawn.  An arrow runs from a module to each module that imports it, and an arrow implied by others is left out.  The library's modules are drawn as one box.}
\label{fig:modules}
\end{figure}
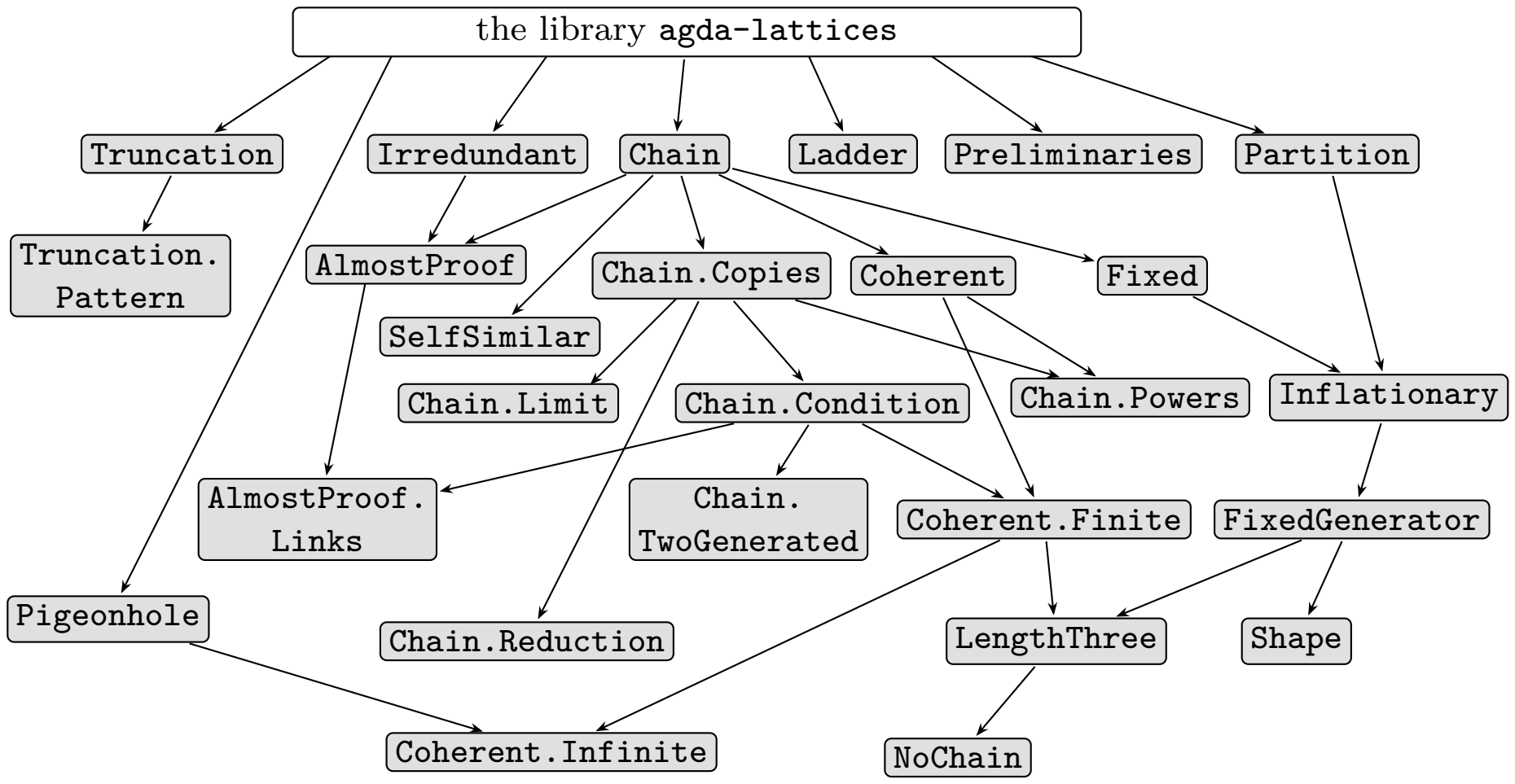

Figures~\ref{fig:dependencies} and~\ref{fig:dependencies-preliminaries} draw what the formal proof of Theorem~\ref{thm:H} uses, in the style of Figure~\ref{fig:map} but complete: each box is a result on the path of \agda{no-copy-chain}, with the definition of Table~\ref{tab:names} that states it, down to the results of Section~\ref{sec:preliminaries}, and an arrow runs from a result to each result whose formal proof uses it through definitions that belong to no result, an arrow implied by two others left out.  Where Figure~\ref{fig:map}'s arrows are citations in the text, these are uses in the development, traced from its sources at the commits above and checked against them by a script.  Twenty-four of the forty-eight results are on the path.  The others are not, among them the coherent form of Theorem~\formofC{coherent} and Corollary~\ref{cor:dean-iff-coherent}, Theorem~\ref{thm:D}, Theorem~\ref{thm:F}, Corollary~\ref{cor:E1}, and the results of Sections~\ref{sec:chains}, \ref{sec:theorem-a} and~\ref{sec:almost-proof}.  On the results that Figure~\ref{fig:map} draws, the development's uses are the text's citations, arrow for arrow.  Below them, in Section~\ref{sec:preliminaries}, they differ, and three of the differences bear on what the checked proof rests on.  The development proves Lemma~\ref{lem:9.13}, which the text cites from FJN, from Theorem~\ref{thm:refine}, which the text cites only for the truncation theorem.  Where the proof of Proposition~\ref{prop:fourteen} invokes $\W$, the development decides by Whitman's procedure, so Theorem~\ref{thm:W} is not on the path.  And the canonical form enters through its construction and the uniqueness of canonical terms, never through the rank by which Theorem~\ref{thm:canonical} defines it.

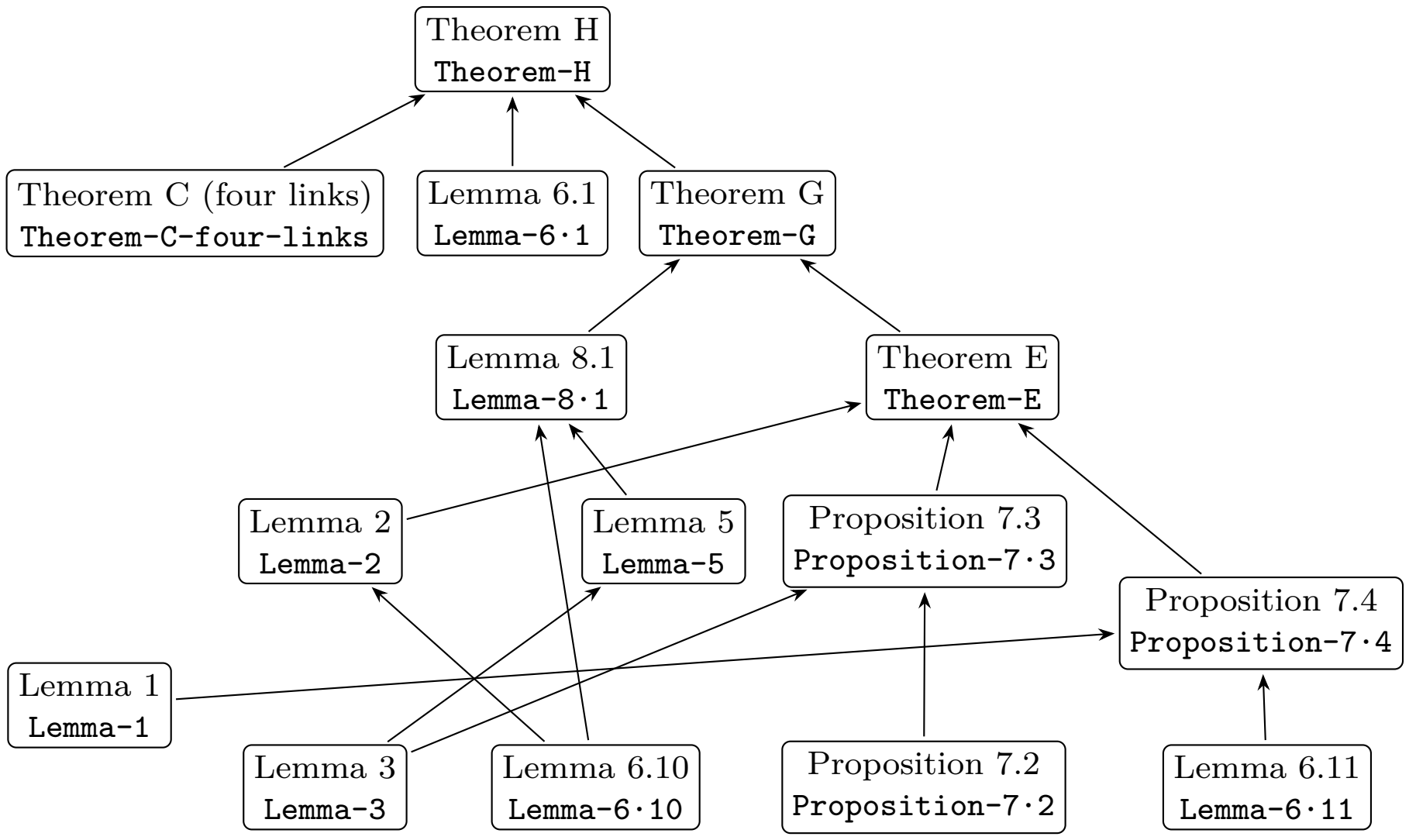
\begin{figure}[!tbp]
\centering
\begin{tikzpicture}[x=10mm, y=10mm, every node/.append style={font=\scriptsize}]
  \node[result] (TH) at (-1.65,-0.00) {Theorem H\\\agda{Theorem-H}};
  \node[result] (L61) at (-1.65,-1.40) {Lemma 6.1\\\agda{Lemma-6·1}};
  \node[result] (C4) at (-4.36,-1.40) {Theorem C (four links)\\\agda{Theorem-C-four-links}};
  \node[result] (TG) at (0.27,-1.40) {Theorem G\\\agda{Theorem-G}};
  \node[result] (L81) at (-1.49,-2.80) {Lemma 8.1\\\agda{Lemma-8·1}};
  \node[result] (TE) at (2.19,-2.80) {Theorem E\\\agda{Theorem-E}};
  \node[result] (L2) at (-3.29,-4.20) {Lemma 2\\\agda{Lemma-2}};
  \node[result] (L5) at (-0.36,-4.20) {Lemma 5\\\agda{Lemma-5}};
  \node[result] (P73) at (1.87,-4.20) {Proposition 7.3\\\agda{Proposition-7·3}};
  \node[result] (P74) at (4.74,-4.90) {Proposition 7.4\\\agda{Proposition-7·4}};
  \node[result] (L1) at (-5.26,-5.60) {Lemma 1\\\agda{Lemma-1}};
  \node[result] (L3) at (-3.25,-6.30) {Lemma 3\\\agda{Lemma-3}};
  \node[result] (L610) at (-0.94,-6.30) {Lemma 6.10\\\agda{Lemma-6·10}};
  \node[result] (L611) at (4.79,-6.30) {Lemma 6.11\\\agda{Lemma-6·11}};
  \node[result] (P72) at (1.86,-6.30) {Proposition 7.2\\\agda{Proposition-7·2}};
  \begin{scope}[on background layer]
    \draw[uses] (L61) -- (TH);
    \draw[uses] (C4) -- (TH);
    \draw[uses] (TG) -- (TH);
    \draw[uses] (L81) -- (TG);
    \draw[uses] (TE) -- (TG);
    \draw[uses] (L2) -- (TE);
    \draw[uses] (L5) -- (L81);
    \draw[uses] (P73) -- (TE);
    \draw[uses] (P74) -- (TE);
    \draw[uses] (L610) -- (L81);
    \draw[uses] (L1) -- (P74);
    \draw[uses] (L3) -- (L5);
    \draw[uses] (L3) -- (P73);
    \draw[uses] (L610) -- (L2);
    \draw[uses] (L611) -- (P74);
    \draw[uses] (P72) -- (P73);
  \end{scope}
\end{tikzpicture}
\caption{What the formal proof of Theorem~\ref{thm:H} uses, above the results of Section~\ref{sec:preliminaries}: a box for each result on its path, with the Agda definition whose type is its statement, and an arrow from a result to each result whose formal proof uses it, at the release's tag \texttt{v0.1.0} with the library at the commit it pins.  An arrow implied by two others is left out.  The results of Section~\ref{sec:preliminaries} enter as Figure~\ref{fig:dependencies-preliminaries} shows.}
\label{fig:dependencies}
\end{figure}

\begin{figure}[!tbp]
\centering
\begin{tikzpicture}[x=10mm, y=10mm, every node/.append style={font=\scriptsize}]
  \node[result] (L1) at (4.60,1.95) {Lemma 1\\\agda{Lemma-1}};
  \node[result] (L2) at (4.60,0.81) {Lemma 2\\\agda{Lemma-2}};
  \node[result] (L3) at (4.60,-2.07) {Lemma 3\\\agda{Lemma-3}};
  \node[result] (L61) at (4.60,-0.17) {Lemma 6.1\\\agda{Lemma-6·1}};
  \node[result] (L610) at (4.60,4.40) {Lemma 6.10\\\agda{Lemma-6·10}};
  \node[result] (L611) at (4.60,-1.22) {Lemma 6.11\\\agda{Lemma-6·11}};
  \node[result] (P72) at (4.60,-4.54) {Proposition 7.2\\\agda{Proposition-7·2}};
  \node[result] (P74) at (4.60,-3.24) {Proposition 7.4\\\agda{Proposition-7·4}};
  \node[result] (C4) at (4.60,3.42) {Theorem C (four links)\\\agda{Theorem-C-four-links}};
  \node[result] (C25) at (0.20,-1.53) {Corollary 2.5\\\agda{Corollary-2·5}};
  \node[result] (L22) at (1.40,3.65) {Lemma 2.2\\\agda{Lemma-2·2}};
  \node[result] (L29) at (0.20,-2.51) {Lemma 2.9\\\agda{Lemma-2·9}};
  \node[result] (TB) at (0.4,0.34) {Theorem B\\\agda{Theorem-B}};
  \node[result] (C24) at (-0.2,1.76) {Corollary 2.4\\\agda{Corollary-2·4}};
  \node[result] (L28) at (-1.60,-0.06) {Lemma 2.8\\\agda{Lemma-2·8}};
  \node[result] (T27) at (-2.80,0.90) {Theorem 2.7\\\agda{Theorem-2·7}};
  \node[result] (T23) at (-4.60,-1.88) {Theorem 2.3\\\agda{Theorem-2·3}};
  \node[result] (T26) at (-3.40,-3.33) {Theorem 2.6\\\agda{Theorem-2·6}};
  \begin{scope}[on background layer]
    \draw[uses] (C25) -- (L1);
    \draw[uses] (C25) -- (L2);
    \draw[uses] (C25) -- (L3);
    \draw[uses] (C25) -- (L61);
    \draw[uses] (C25) -- (L611);
    \draw[uses] (C25) -- (C4);
    \draw[uses] (L22) -- (L1);
    \draw[uses] (L22) -- (L610);
    \draw[uses] (L29) -- (L611);
    \draw[uses] (TB) -- (L61);
    \draw[uses] (TB) -- (P74);
    \draw[uses] (C24) -- (L2);
    \draw[uses] (C24) -- (L611);
    \draw[uses] (C24) -- (C4);
    \draw[uses] (L28) -- (L3);
    \draw[uses] (L28) -- (L611);
    \draw[uses] (T23) -- (L2);
    \draw[uses] (T23) -- (P72);
    \draw[uses] (T26) -- (P72);
    \draw[uses] (C24) -- (TB);
    \draw[uses] (L28) -- (TB);
    \draw[uses] (T23) -- (L29);
    \draw[uses] (T26) -- (L29);
    \draw[uses] (T27) -- (L28);
    \draw[uses] (T23) -- (T27);
    \draw[uses] (T26) -- (T27);
  \end{scope}
\end{tikzpicture}
\caption{How the results of Section~\ref{sec:preliminaries} on the path of Figure~\ref{fig:dependencies} enter it: an arrow from a result to each result whose formal proof uses it, the column on the right being the results of Figure~\ref{fig:dependencies} that use them.  Lemma~\ref{lem:prime} enters as the partition lemma at~$z$, Corollary~\ref{cor:1.6} as its form for a generator in the image of a homomorphism, Corollary~\ref{cor:1.13} through its two directions, and Theorem~\ref{thm:canonical} through the canonicalization and the uniqueness of canonical terms.}
\label{fig:dependencies-preliminaries}
\end{figure}
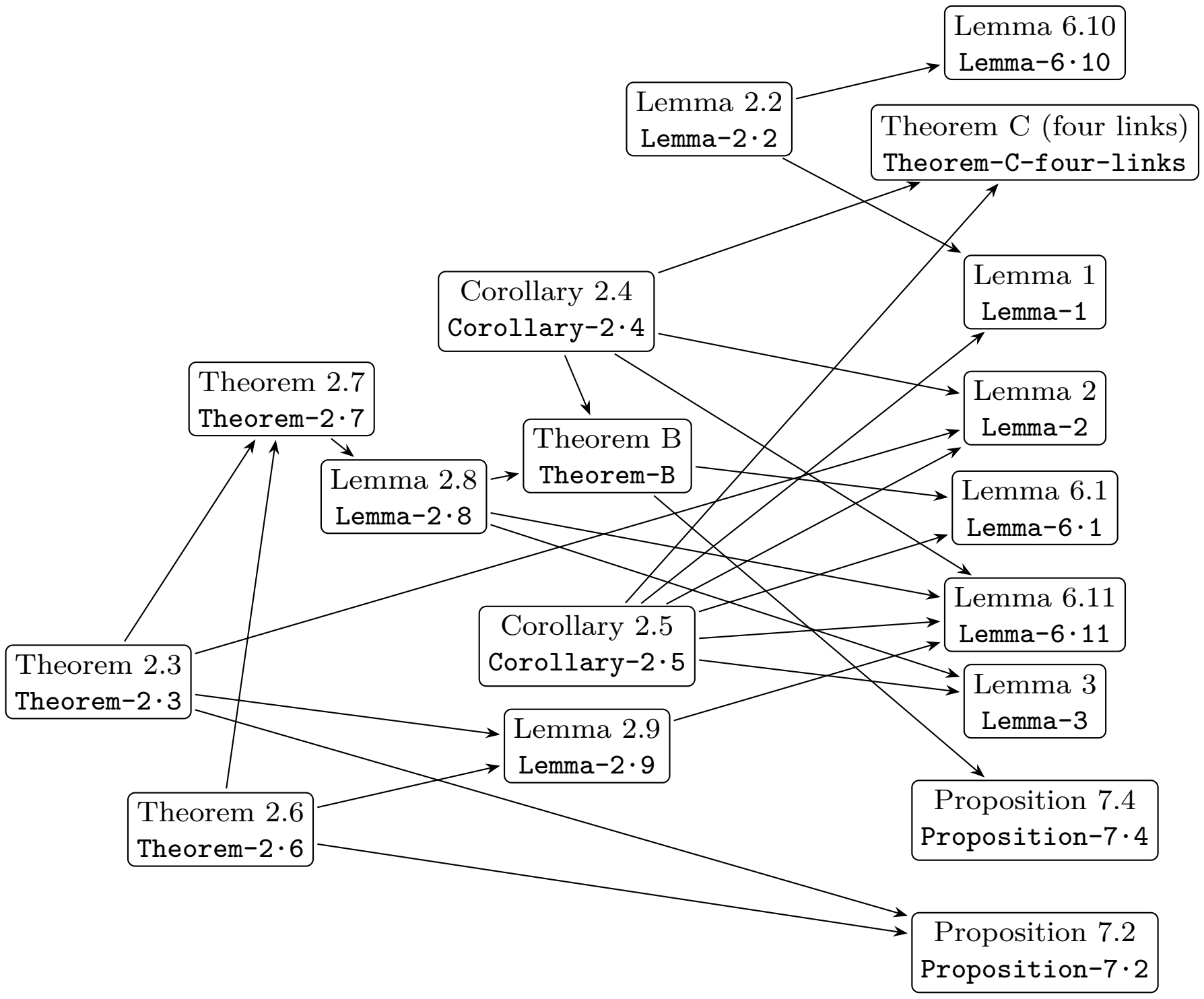

\begingroup
\footnotesize
\setlength{\LTleft}{0pt}
\setlength{\LTright}{0pt}
\begin{longtable}{@{}>{\raggedright\arraybackslash}p{0.18\textwidth}>{\raggedright\arraybackslash}p{0.30\textwidth}>{\raggedright\arraybackslash}p{0.45\textwidth}@{}}
\caption{Each result of this paper, the name of the Agda definition that states it, and the module that holds it; the last rows give the definitions of Sections~\ref{sec:chains}, \ref{sec:coherent} and~\ref{sec:structure}.  Every module is a submodule of \agda{DeansProblem} in \agda{ualib/agda-deans-problem}.  A statement marked $\dagger$ takes decidable equality on the generators as a hypothesis, where a decision needs it; one marked $\ast$ takes the infinite pigeonhole principle, which the checked proof of Theorem~\ref{thm:H} does not use.}
\label{tab:names}\\
\toprule
Result & Statement & Module \\
\midrule
\endfirsthead
\toprule
Result & Statement & Module \\
\midrule
\endhead
\bottomrule
\endfoot
Theorem~\ref{thm:H} & \agda{Theorem-H}\textsuperscript{$\dagger$} & \agda{NoChain} \\
Theorem~\ref{thm:W} & \agda{Theorem-2·1} & \agda{Preliminaries} \\
Lemma~\ref{lem:prime} & \agda{Lemma-2·2} & \agda{Preliminaries} \\
Theorem~\ref{thm:whitman} & \agda{Theorem-2·3}\textsuperscript{$\dagger$} & \agda{Preliminaries} \\
Corollary~\ref{cor:1.6} & \agda{Corollary-2·4} & \agda{Preliminaries} \\
Corollary~\ref{cor:1.13} & \agda{Corollary-2·5} & \agda{Preliminaries} \\
Theorem~\ref{thm:canonical} & \agda{Theorem-2·6}\textsuperscript{$\dagger$} & \agda{Preliminaries} \\
Theorem~\ref{thm:refine} & \agda{Theorem-2·7}\textsuperscript{$\dagger$} & \agda{Preliminaries} \\
Lemma~\ref{lem:9.13} & \agda{Lemma-2·8}\textsuperscript{$\dagger$} & \agda{Preliminaries} \\
Theorem~\ref{thm:B} & \agda{Theorem-B}\textsuperscript{$\dagger$} & \agda{Preliminaries} \\
Lemma~\ref{lem:finite} & \agda{Lemma-2·9}\textsuperscript{$\dagger$} & \agda{Preliminaries} \\
Proposition~\ref{prop:reformulation} & \agda{Proposition-4·2} & \agda{Chain.\allowbreak Limit} \\
Proposition~\ref{prop:self-similar} & \agda{Proposition-4·3} & \agda{Chain} \\
Lemma~\ref{lem:ladder} & \agda{Lemma-4·4} & \agda{Ladder} \\
Theorem~\ref{thm:ladder} & \agda{Theorem-4·5} & \agda{Examples.\allowbreak ChainCondition} \\
Theorem~\ref{thm:A} & \agda{Theorem-A}\textsuperscript{$\dagger$} & \agda{Chain.\allowbreak Condition} \\
Corollary~\ref{cor:A1} & \agda{Corollary-A1}\textsuperscript{$\dagger$} & \agda{SelfSimilar} \\
Corollary~\ref{cor:A2} & \agda{Corollary-A2} & \agda{Examples.\allowbreak FreeseEquations} \\
Corollary~\ref{cor:A3} & \agda{Corollary-A3}\textsuperscript{$\dagger$} & \agda{Chain.\allowbreak Condition} \\
Corollary~\ref{cor:A4} & \agda{Corollary-A4} & \agda{Fixed} \\
Lemma~\ref{lem:bounded-transitions} & \agda{Lemma-6·1}\textsuperscript{$\dagger$} & \agda{Coherent} \\
Theorem~\ref{thm:C} (one step) & \agda{Theorem-C-one-step}\textsuperscript{$\dagger$} & \agda{Coherent} \\
Lemma~\ref{lem:coherent-images} & \agda{Lemma-6·3} & \agda{Coherent.\allowbreak Finite} \\
Theorem~\ref{thm:C-four} & \agda{Theorem-C-four-links} & \agda{Coherent.\allowbreak Finite} \\
Theorem~\ref{thm:C-three} & \agda{Theorem-C-length-three}\textsuperscript{$\dagger$} & \agda{Coherent.\allowbreak Finite} \\
Theorem~\ref{thm:C-coherent} & \agda{Theorem-C-coherent}\textsuperscript{$\ast$}\textsuperscript{$\dagger$} & \agda{Coherent.\allowbreak Infinite} \\
Corollary~\ref{cor:dean-iff-coherent} & \agda{Corollary-6·4}\textsuperscript{$\ast$}\textsuperscript{$\dagger$} & \agda{Coherent.\allowbreak Infinite} \\
Proposition~\ref{prop:unbounded} & \agda{Proposition-6·6} & \agda{Coherent.\allowbreak Finite} \\
Theorem~\ref{thm:D} & \agda{Theorem-D} & \agda{Examples.\allowbreak FixedTriple} \\
Proposition~\ref{prop:rigid} & \agda{Proposition-6·9} & \agda{Examples.\allowbreak FixedTriple} \\
Lemma~\ref{lem:two-generators} & \agda{Lemma-6·10} & \agda{Partition} \\
Lemma~\ref{lem:descent} & \agda{Lemma-6·11} & \agda{Fixed} \\
Lemma~\ref{lem:1} & \agda{Lemma-1} & \agda{Inflationary} \\
Lemma~\ref{lem:2} & \agda{Lemma-2} & \agda{Inflationary} \\
Lemma~\ref{lem:3} & \agda{Lemma-3} & \agda{Inflationary} \\
Lemma~\ref{lem:4} & \agda{Lemma-4} & \agda{Inflationary} \\
Proposition~\ref{prop:fourteen} & \agda{Proposition-7·2} & \agda{FixedGenerator} \\
Proposition~\ref{prop:eleven} & \agda{Proposition-7·3} & \agda{FixedGenerator} \\
Proposition~\ref{prop:inflationary-power} & \agda{Proposition-7·4} & \agda{FixedGenerator} \\
Theorem~\ref{thm:E} & \agda{Theorem-E} & \agda{FixedGenerator} \\
Corollary~\ref{cor:freese-eleven} & \agda{Corollary-7·5} & \agda{Examples.\allowbreak SelfEmbedding} \\
Corollary~\ref{cor:E1} & \agda{Corollary-E1} & \agda{Shape} \\
Theorem~\ref{thm:F} & \agda{Theorem-F} & \agda{Inflationary} \\
Lemma~\ref{lem:5} & \agda{Lemma-5} & \agda{LengthThree} \\
Lemma~\ref{lem:core} & \agda{Lemma-8·1} & \agda{LengthThree} \\
Theorem~\ref{thm:G} & \agda{Theorem-G} & \agda{LengthThree} \\
Theorem~\ref{thm:truncation} & \agda{Theorem-9·1} & \agda{Truncation} \\
Theorem~\ref{thm:box-agree} & \agda{Theorem-9·2} & \agda{AlmostProof.\allowbreak Links} \\
\midrule
Definition~\ref{def:chain} & \agda{CopyChain}, \agda{Chain}, \agda{Step} & \agda{Chain.\allowbreak Copies}, \agda{Chain} \\
Definition~\ref{def:coherent} & \agda{Coherent}, \agda{Coherentω} & \agda{Coherent.\allowbreak Finite} \\
Definition~\ref{def:inflationary} & \agda{Inflationary}, \agda{Deflationary} & \agda{Inflationary} \\
\end{longtable}
\endgroup

\subsection{The cost}
\label{sec:cost}

Checking the release from source (its $38$ modules and the $147$ modules of the two libraries that they import, $42$ of \agda{agda-lattices} and $105$ of \agda{agda-algebras}), through a generated aggregator, takes about $35$ seconds, with a peak resident memory of $1.95$ GiB, on a 20-core machine.\footnote{The machine resources were shared with other work at a load average of $1.0$ to $2.4$; we used Agda 2.9.0 with a heap limit of $6$ GiB.}  With the libraries' interfaces in place, the modules of \agda{agda-deans-problem} check in about $15$ seconds.  The peak follows the garbage collector's schedule: each run from source peaked at $1.95$ GiB, and each with the libraries' interfaces in place at $2.10$ GiB.

The greatest cost is due to Whitman decisions.  Every decision of the release runs through the recursion except those that Section~\ref{sec:what-checked} names, on the terms of the one-step example, which run through the table.  On the terms of the one-step example the recursion is expensive; in modules outside the release, measured at its earlier commit \texttt{227ffdc} of 24 October 2026, free generation of the defining triple $(a\pr, b\pr, c\pr)$ exhausts the build's heap, and the equation of item (iii) alone takes $11.8$ seconds, at a peak of $5.7$ GiB; with the memoized table approach each of these takes at most $3.3$ seconds, about $1.6$ of them the loading of the module's imports (Section~\ref{sec:verification}).

Every other decision on which a statement of this paper rests is far cheaper: \agda{DeansProblem.\allowbreak Examples.\allowbreak FreeseBox}, which decides a few inequalities on the canonical form of $a$ and the canonicity of the canonical forms of $a$, $b$ and $c$ (540 Whitman decisions each), takes $1.9$ seconds in all, and every decision behind Theorems~\ref{thm:D}, \ref{thm:E} and~\ref{thm:G} is in $\FL(3)$, on terms with at most 30 generator occurrences, and takes milliseconds.  The proof of Theorem~\ref{thm:H} runs no decision while it is checked: its searches and bounds are functions of the chain, which is a variable.

\subsection{The searches}
\label{sec:searches}

\begin{computation}
\label{comp:searches}
The searches that led to Theorems~\ref{thm:D}, \ref{thm:E} and~\ref{thm:G} ran in Python, with a memoized implementation of Whitman's procedure and an enumeration of $\FL(3)$ by complexity (canonical forms as antichains of generators and meets satisfying the conditions of Theorem~\ref{thm:canonical}, dually for meets); the twelve symmetries of $\FL(3)$ (permutations of the generators and duality) reduce the spaces.  Every positive fact they found is decided in Agda; the counts are computations, not theorems, and no proof cites one.  ``Within the bounds'' means: fixed points of complexity at most $4$, the sublattice they generate capped at $150$ elements, and the free triples among those.
\begin{itemize}
\item Forward, over the $3318$ proper self-embeddings whose generator images have complexity at most $3$ ($305$ classes under the symmetries, of which $90$ have a periodic element of complexity $1$, as the descent requires of any embedding with a fixed point of positive complexity): one class with a fixed free triple within the bounds, that of $\theta_1$.
\item The class fixing $z$ with $\theta(x) \in [x, x \vee c]$, $\theta(y) \in [y, y \vee c]$, where $c = z \wedge (x \vee y)$, complexity at most $4$: $25$ candidates, $7$ proper members that freely generate, all with fixed free triples within the bounds; $6$ have a meet of two generators as a node of an image and $\theta_2$ the node $y \wedge (x \vee z)$; all $7$ are rigid by Theorem~\ref{thm:F}.
\item The join-fixing family moving every generator ($x < \theta(x) \leq x \vee y$, $y < \theta(y) \leq x \vee y$, $z < \theta(z)$, complexity at most $4$): $6535$ free triples, $1159$ up to symmetry, none with a fixed free triple within the bounds; now Theorem~\ref{thm:E}.
\item Backward, for each of the $7$ class members, every pair $(a\pr, b\pr)$ of elements of complexity at most $4$ with $(a\pr, b\pr, z)$ fixing the member's images: only $(x, y)$, except two members with one non-injective extra solution each; now Theorem~\ref{thm:G}.
\item Top-down, free triples $(a, b, z)$ inside $\Fix(\theta_1)$ and $\Fix(\theta_2)$, each the sublattice generated by its $18$ fixed points of complexity at most $4$, capped at $1200$ elements: none; now Theorem~\ref{thm:G}.
\item Every self-embedding fixing $z$ with images of complexity at most $4$ ($107\,826$ free triples $(a, b, z)$ with $a, b$ in the ideal $(x \vee y]$), in $28$ minutes: $352$ with a periodic element of complexity $1$, and $14$ with a fixed free triple within the bounds, the $7$ class members above and their duals, all of the shape of Corollary~\ref{cor:E1}.
\end{itemize}
\end{computation}
\end{fullonly}

\section{Remarks and questions}
\label{sec:questions}

\begin{enumerate}[label=(\arabic*), leftmargin=*]
\item \textbf{The chain condition without freeness}.  Theorem~\ref{thm:ladder} answers the variant of Section~\ref{sec:takahasi} in the negative: an ascending chain of $3$-generated sublattices of a free lattice need not stabilize.  Its links appear to be finite, of $10 + 6k$ and $12 + 6k$ elements (Computation~\ref{comp:ladder}).  Are they?  If so, the finite $3$-generated sublattices of a free lattice have unbounded size, and the chain condition fails even for chains of finite sublattices.  And for which $3$-generated lattices $K$ does a free lattice contain an infinite strictly ascending chain of sublattices isomorphic to $K$?  Theorem~\ref{thm:H} says that $K = \FL(3)$ is not one; the links of Theorem~\ref{thm:ladder}, if finite, are pairwise non-isomorphic, so they do not settle the question for any $K$.

\item \textbf{The self-embeddings with a fixed free triple}.  Corollary~\ref{cor:E1} describes their shape: a fixed generator $z$, and, up to duality and a power, $x$ and $y$ moved up inside $(x \vee y) \wedge (x \vee z)$ and $(x \vee y) \wedge (y \vee z)$.  The search of Computation~\ref{comp:searches} found seven such self-embeddings within complexity $4$ in the class fixing $z$ with $\theta(x) \in [x, x \vee c]$ and $\theta(y) \in [y, y \vee c]$, where $c = z \wedge (x \vee y)$, and, among all those fixing $z$ within complexity $4$ and its bounds, no others but their duals.  Which self-embeddings of $\FL(3)$ have a fixed free triple?  Is $\Fix(\theta_1)$ exactly the sublattice $\gen{p, q, z}$ of Theorem~\ref{thm:D}, as the eighteen fixed points of complexity at most $4$ suggest?\footnote{See \cite[Theorem~6.1]{Lupini2026} and the footnote{} on that preprint in Section~\ref{sec:problem}.}

\item \textbf{Fixed sublattices in general}.  The fixed subgroup of an injective endomorphism of a free group of rank $n$ has rank at most $n$ \cite{ImrichTurner1989, BestvinaHandel1992}.  Is the fixed sublattice of a proper self-embedding of $\FL(3)$ generated by three elements, or at least finitely generated?  Theorem~\ref{thm:D} does not decide it: $\Fix(\theta_1)$ contains a copy of $\FL(\omega)$, but so does $\FL(3)$, and the part of $\Fix(\theta_1)$ that is known, $\gen{p, q, z}$, has three generators.  Corollary~\ref{cor:A4} identifies the periodic points with the intersection of the images, and Lemma~\ref{lem:descent} says every fixed point of positive complexity descends to a join or meet of generators fixed by a power of the self-embedding.  What else can be said about $\Fix(\theta)$ and $\Per(\theta)$ for a proper self-embedding $\theta$ of $\FL(3)$, for instance whether $\Fix(\theta)$ is always finitely generated?
\end{enumerate}

\section*{Acknowledgments}

The first author thanks the second for the 1997 notes that started this work and for the problem's history.  The formal development described in Section~\ref{sec:formal}, the computer searches and the first draft of this paper were produced by large language models (Claude Fable~5.1 and Claude Opus~5.5, Anthropic; GPT-6-Astra, OpenAI; Kimi~K3, Moonshot AI) working as agents under the first author's direction, in sessions in which the proofs of the new results were also worked out, as Section~\ref{sec:how-produced} describes; the authors take full responsibility for every statement and proof.

\bibliographystyle{amsplain}
\bibliography{deans-problem}

\providecommand{\bysame}{\leavevmode\hbox to3em{\hrulefill}\thinspace}
\providecommand{\MR}{\relax\ifhmode\unskip\space\fi MR }
\providecommand{\MRhref}[2]{%
  \href{http://www.ams.org/mathscinet-getitem?mr=#1}{#2}
}
\providecommand{\href}[2]{#2}
\begin{thebibliography}{10}

\bibitem{Agda}
{Agda Developers}, \emph{Agda}, 2026, version 2.9.0, not yet released: the
  development version at commit \texttt{da66a8c} (5~October 2026) of
  \url{https://github.com/agda/agda}; documentation at
  \url{https://agda.readthedocs.io/}.

\bibitem{agda-mcp-field-report}
{Agda-native AIR}, \emph{2026-10-20, a private repository: agda-mcp across the
  seventeen sessions that wrote one formal development}, 2026, a field report
  of the Agda-native AIR project, written by Claude Opus~5.5 (Anthropic) at the
  first author's request from the sessions' own records; the file
  \texttt{docs/field-reports/2026-10-20-private-development-agda-mcp-use.md} at
  commit \texttt{137d9fb} (21~October 2026),
  \url{https://github.com/formalverification/agda-native-air/blob/137d9fb/docs/field-reports/2026-10-20-private-development-agda-mcp-use.md}.

\bibitem{BestvinaHandel1992}
Mladen Bestvina and Michael Handel, \emph{Train tracks and automorphisms of
  free groups}, Ann. of Math. (2) \textbf{135} (1992), no.~1, 1--51.
  \MR{1147956}

\bibitem{agda-stdlib}
Matthew~L. Daggitt, Guillaume Allais, James McKinna, Andreas Abel, Nathan van
  Doorn, James Wood, Ulf Norell, Donnacha~Ois{\'\i}n Kidney, Sergei
  Meshveliani, Sandro Stucki, Jacques Carette, Alex Rice, Jason Z.~S. Hu,
  Li-yao Xia, Shu-Hung You, Reed Mullanix, and Wen Kokke, \emph{The {A}gda
  standard library: version 2.0}, J. Open Source Softw. \textbf{10} (2025),
  no.~116, 9241, version~2.3 with the five changes it needs under Agda~2.9.0 is
  the tag \texttt{v2.3-agda-2.9.0} of
  \url{https://github.com/formalverification/agda-stdlib}.

\bibitem{Dean1961}
R.~A. Dean, \emph{Sublattices of free lattices}, Lattice Theory (R.~P.
  Dilworth, ed.), Proceedings of Symposia in Pure Mathematics, vol.~II,
  American Mathematical Society, Providence, RI, 1961, pp.~31--42. \MR{0133263}

\bibitem{agda-deans-problem}
William DeMeo, \emph{agda-deans-problem: the formal development of ``{N}o free
  lattice contains an infinite ascending chain of copies of {FL}(3)''}, 2026,
  the tag \texttt{v0.1.0} (commit \texttt{d562c05}, 19~November 2026),
  Apache-2.0, released with this paper and archived on Zenodo,
  \url{https://doi.org/10.5281/zenodo.24375155}; source at
  \url{https://github.com/ualib/agda-deans-problem/tree/v0.1.0}; it pins
  agda-lattices at commit \texttt{41fb55e}.

\bibitem{agda-lattices}
\bysame, \emph{agda-lattices: lattice theory in {A}gda}, 2026, Apache-2.0, at
  commit \texttt{41fb55e} (30~October 2026, its release v0.2.0) of
  \url{https://github.com/ualib/agda-lattices}.

\bibitem{agda-native-air}
\bysame, \emph{{A}gda-native {AIR}: an {A}gda-native artificial intelligence
  reasoning environment}, 2026, version 0.2.0, which includes the agda-mcp
  server, Apache-2.0, archived on Zenodo,
  \url{https://doi.org/10.5281/zenodo.23378575}; source at
  \url{https://github.com/formalverification/agda-native-air/tree/v0.2.0}.

\bibitem{agda-algebras}
William DeMeo and Jacques Carette, \emph{The {A}gda {U}niversal {A}lgebra
  {L}ibrary (agda-algebras)}, 2026, the development version toward version 3.0,
  at commit \texttt{b4d0ab6} (24~October 2026) of
  \url{https://github.com/ualib/agda-algebras}.

\bibitem{Freese1997}
Ralph Freese, \emph{Notes on {D}ean's problem}, manuscript dated 29 July 1997,
  four pages, \url{https://math.hawaii.edu/~ralph/Notes/DeansProb/dean.pdf},
  1997.

\bibitem{FreeseAppendixG}
\bysame, \emph{Free lattices}, General Lattice Theory, Birkh{\"a}user, Basel,
  second ed., 1998, Appendix~G of G.~Gr{\"a}tzer's book, pp.~575--590.
  \MR{1670580}

\bibitem{FreeseAppendixGRevised}
\bysame, \emph{Free lattices}, 2002, revised version of Appendix~G of
  Gr{\"a}tzer's \emph{General Lattice Theory},
  \url{https://math.hawaii.edu/~ralph/Preprints/appendix.pdf}.

\bibitem{FJN}
Ralph Freese, Jaroslav Je{\v{z}}ek, and J.~B. Nation, \emph{Free lattices},
  Mathematical Surveys and Monographs, vol.~42, American Mathematical Society,
  Providence, RI, 1995. \MR{1319815}

\bibitem{GuilloudPitClaudel2025}
Simon Guilloud and Cl{\'e}ment Pit-Claudel, \emph{Verified and optimized
  implementation of orthologic proof search}, Computer Aided Verification
  ({CAV} 2025), Part~I (Ruzica Piskac and Zvonimir Rakamari{\'c}, eds.),
  Lecture Notes in Computer Science, vol. 15933, Springer, Cham, 2025,
  arXiv:2501.09418, pp.~130--152.

\bibitem{Higman1951}
Graham Higman, \emph{A finitely related group with an isomorphic proper factor
  group}, J. London Math. Soc. \textbf{26} (1951), 59--61. \MR{0038347}

\bibitem{ImrichTurner1989}
W.~Imrich and E.~C. Turner, \emph{Endomorphisms of free groups and their fixed
  points}, Math. Proc. Cambridge Philos. Soc. \textbf{105} (1989), no.~3,
  421--422. \MR{0985677}

\bibitem{JamesHinze2011}
Daniel W.~H. James and Ralf Hinze, \emph{A reflection-based proof tactic for
  lattices in {C}oq}, Trends in Functional Programming, Volume 10 (Zolt{\'a}n
  Horv{\'a}th, Vikt{\'o}ria Zs{\'o}k, Peter Achten, and Pieter Koopman, eds.),
  Intellect, Bristol, 2011, TFP 2009, Kom{\'a}rno;
  \url{https://www.cs.ox.ac.uk/people/Ralf.Hinze/publications/TFP09.pdf},
  pp.~97--112.

\bibitem{JonssonKiefer1962}
Bjarni J{\'o}nsson and James~E. Kiefer, \emph{Finite sublattices of a free
  lattice}, Canad. J. Math. \textbf{14} (1962), 487--497. \MR{0137667}

\bibitem{KearnesKiss}
Keith~A. Kearnes and Emil~W. Kiss, \emph{The shape of congruence lattices},
  vol. 222, Memoirs of the American Mathematical Society, no. 1046, American
  Mathematical Society, Providence, RI, 2013.

\bibitem{Lupini2026}
Martino Lupini, \emph{Ascending chains of sublattices in free lattices}, 2026,
  preprint dated 19 October 2026, Zenodo,
  \url{https://doi.org/10.5281/zenodo.23477016}.

\bibitem{Nation1982}
J.~B. Nation, \emph{Finite sublattices of a free lattice}, Trans. Amer. Math.
  Soc. \textbf{269} (1982), no.~1, 311--337. \MR{0637041}

\bibitem{Takahasi1950}
Mutuo Takahasi, \emph{Note on locally free groups}, J. Inst. Polytech. Osaka
  City Univ. Ser. A \textbf{1} (1950), 65--70. \MR{0041845}

\bibitem{Takahasi1951}
\bysame, \emph{Note on chain conditions in free groups}, Osaka Math. J.
  \textbf{3} (1951), 221--225. \MR{0046362}

\bibitem{Whitman1941}
Philip~M. Whitman, \emph{Free lattices}, Ann. of Math. (2) \textbf{42} (1941),
  no.~1, 325--330. \MR{0003614}

\bibitem{Whitman1942}
\bysame, \emph{Free lattices {II}}, Ann. of Math. (2) \textbf{43} (1942),
  no.~1, 104--115. \MR{0006143}

\end{thebibliography}

\end{document}